\documentclass[preprint,11pt]{elsarticle}
\makeatletter
\def\ps@pprintTitle{%
 \let\@oddhead\@empty
 \let\@evenhead\@empty
 \def\@oddfoot{\centerline{\thepage}}%
 \let\@evenfoot\@oddfoot}
\makeatother

\usepackage{lineno}
\modulolinenumbers[5]
\usepackage[normalem]{ulem}

\journal{Paper accepted to European Journal of Operational Research}

\usepackage{graphicx}
\usepackage[english]{babel,varioref}
\usepackage[utf8]{inputenc}
\usepackage[T1]{fontenc}

 \usepackage[table,xcdraw]{xcolor}			
\usepackage{float}
\usepackage[numbers]{natbib}
\usepackage[hmargin=1.75cm,vmargin=2cm,bmargin=2cm]{geometry}
\usepackage{amssymb,amsthm,amsmath,dsfont}
\usepackage{setspace}
\usepackage{verbatim} 
\usepackage{fancyhdr}
\usepackage{multirow}
\usepackage{eurosym}
\usepackage{epstopdf}

\usepackage{caption}
\usepackage{subfigure}
\usepackage[section]{placeins}
\usepackage{lineno}
\usepackage {tikz}
\usetikzlibrary {positioning}
\definecolor {processblack}{cmyk}{0.96,0,0,0}
\usepackage[lined,algoruled,commentsnumbered, linesnumbered]{algorithm2e}

\usepackage{soul}

\usepackage{float}

\newcommand{\minor}[1]{{\color{black}#1}}
\newcommand{\typo}[1]{{\color{black}#1}}

\usepackage{xr-hyper}
\usepackage{hyperref}

\makeatletter
\newcommand*{\addFileDependency}[1]{
 \typeout{(#1)}
 \@addtofilelist{#1}
  \IfFileExists{#1}{}{\typeout{No file #1.}}
}
\makeatother

\newcommand*{\myexternaldocument}[1]{
    \externaldocument{#1}
    \addFileDependency{#1.tex}
    \addFileDependency{#1.aux}
}

\myexternaldocument{table_related}
\myexternaldocument{table_labels}
\myexternaldocument{table_params}

\usepackage{pdfpages} 

\begin{document}
\begin{frontmatter}



\title{An efficient Lagrangian-based heuristic to solve a multi-objective sustainable supply chain  problem}


\author[l1]{Camila P.S. Tautenhain}
\ead{santos.camila@unifesp.br}

\author[l2]{Ana Paula Barbosa-Povoa}
\ead{apovoa@tecnico.ulisboa.pt}

\author[l2]{Bruna Mota}
\ead{bruna.mota@tecnico.ulisboa.pt}

\author[l1]{Mariá C.V. Nascimento\corref{cor1}}
\ead{mcv.nascimento@unifesp.br}

\cortext[cor1]{Corresponding author.}
\address[l1]{Instituto de Ciência e Tecnologia, Universidade Federal de São Paulo, 
São José dos Campos, Brasil}

\address[l2]{Centro de Estudos de Gestão, Instituto Superior Técnico, Universidade Técnica de Lisboa, Lisboa 1049-101, Portugal}

\begin{abstract}
Sustainable Supply Chain (SSC) management aims at integrating economic, environmental and social goals to assist in the long-term planning of a company and its supply chains. There is no consensus in the literature as to whether social and environmental responsibilities are profit-compatible. However, the conflicting nature of these goals is explicit when considering specific assessment measures and, in this scenario, multi-objective optimization is a way to represent problems that simultaneously optimize the goals.  This paper proposes a Lagrangian matheuristic method, called \emph{AugMathLagr}, to solve a hard and relevant multi-objective problem found in the literature.
\emph{AugMathLagr} was extensively tested using artificial instances defined by a generator presented in this paper. 
The results show a competitive performance of \emph{AugMathLagr} when  compared with an exact multi-objective method limited by time and a matheuristic recently proposed in the literature and adapted here to address the studied problem. In addition, computational results on a case study are presented and analyzed, and demonstrate the outstanding performance of \emph{AugMathLagr}.
\end{abstract}

  \nonumnote{Accepted to the European Journal of Operational Research. ©2021. This manuscript version is made available under the CC-BY-NC-ND 4.0 license http://creativecommons.org/licenses/by-nc-nd/4.0/}
  
\end{frontmatter}


\section{Introduction}

The change in mentality of organizations to go beyond profit maximization and consider the long-term economic  success  by taking into account social and environmental responsibilities has raised the interest of researchers in the subject of Sustainable Supply Chain (SSC) management \citep{Seuring2013, Eskandarpour2015,Povoa2018}.  In particular,  multi-objective SSC management  models have relied on the  integration  of social, environmental and economic issues -- the so-called triple bottom line -- to explicitly define the long-term planning of sustainable companies  \citep{Carter2008}. 
Multi-objective optimization plays an important role in approaching SSC management since a compromise among conflicting goals related to profit, environment and social issues can be reached.

A measure commonly used as an economic criterion in SSC management optimization models is the Net Present Value (NPV) \citep{Copeland1983}.
There is a plethora of measures \citep{Ahi2015} to quantitatively assess the environmental impacts related to, for example, greenhouse gas emissions, pollution and resource usage. Additionally, there exist some methods specifically designed to quantify these impacts such as Eco-Indicator 99 \citep{Goedkoop2000} and ReCiPe 2008 \citep{Goedkoop2009}, which are based on the Life Cycle Analysis (LCA) methodology.

In spite of the relevance of the social impacts in the SSC 
management, a minority of models in the literature optimize this indicator along with environmental and economic indicators, e.g. in {\cite{Mota2015,Boukherroub2015,Cambero2016,Mota2018,Tautenhain2019}}. \citet{Eskandarpour2015} discuss various indicators to assess social responsibilities, most of which related to human rights and social justice laws.
In particular, the authors suggest metrics that encourage the creation of job vacancies \citep{Mota2015,Boukherroub2015,Mota2018},  employment stability \citep{Boukherroub2015} and working conditions. Recently, \citet{Mota2015} introduced the Social Benefit Indicator (SB) whose calculation relies on the total {number of jobs created} in each location across the supply chain. This  indicator weighs preferably the creation of entities in less developed countries. \citet{Mota2018} {have optimized} this social indicator in the multi-objective SSC management problem under investigation. This indicator is recommended \typo{in} situations that involve hiring and do not require layoffs \citep{Mota2018}.
In addition to optimizing the SB, the SSC management problem introduced in \cite{Mota2018} aims at maximizing the NPV and minimizing the environmental impacts quantified by ReCiPe 2008.

The problem in question models a generic SSC management problem that manufactures multiple products in a  planning horizon composed of multiple periods. The appeal of the model lies 
\typo{in} its applicability to a wide variety of industries, since it integrates a number of strategic and tactical decisions related to, for example, the use of technologies to manufacture and refurbish products and the shipment of items through different transport modes, among others.  The authors solved mono-objective problems to minimize the environmental objective function and to maximize the economic, i.e., the NPV, and social objective functions.
Moreover, they have also considered two scenarios to maximize the social objective function with additional constraints for the NPV value. In these scenarios, they required the NPV value to be at least $85\%$ and $95\%$ of the value found when optimizing the mono-objective problem that maximizes the economic objective function. The computational difficulty {to solve} the resulting problem was the main reason why the authors did not consider the three objective functions simultaneously when optimizing the problems in the case studies.

Many studies employ heuristic methods instead of exact solution strategies to solve \typo{SSC management problems within a reasonable time} \citep{Devika2014,Sahebjamnia2018,Govindan2019}. However, to design a heuristic to find feasible solutions for the SSC management problem introduced in \citep{Mota2018} is particularly challenging. This is due to the integration of several complicated constraints involving binary and integer variables with \typo{strategic} planning decisions. In line with this, a promising approach to this problem is the relaxation of complicated constraints and decision variables. Lagrangian decomposition allows \typo{us} to explore these aspects and, in specific, can provide efficient methods for this problem by relaxing complicated constraints while penalizing their violation.

In this context, the primary contributions of this paper are:
\begin{itemize}
\item The development of an efficient Lagrangian matheuristic for SSC management problems, here called \emph{AugMathLagr}, to tackle the SSC problem proposed in \cite{Mota2018}.
\item The adaptation of the \emph{AugMathFix} matheuristic introduced in \citep{Tautenhain2019} to the SSC management problem put forward in \cite{Mota2018}.
\item A test bed of artificial instances and an instance generator loosely based on the study of real data.

\end{itemize}

The present paper also reports the computational experience  with \emph{AugMathLagr} on  a case study and a set of medium-sized randomly generated instances.  Results from of  \emph{AugMathLagr} were compared with  results obtained with the Augmented $\epsilon$-Constraint method (AUGMECON2) \citep{Mavrotas2013}, {which} is an enhancement of AUGMECON \citep{Mavrotas2009}, and with the adapted \emph{AugMathFix}. The comparison involved {evaluating} the solution quality and time-to-solution. The proposed matheuristic achieved very good results in significantly lesser computational times than AUGMECON2 and \emph{AugMathFix}.

The remainder of this paper is organized as follows. 
Section~\ref{sec:relatedMethods} reviews the related methods for solving multi-objective {SSC management} problems{;}
Section \ref{sec:targetModel} presents the {SSC} model to which we propose a Lagrangian heuristic;
Section \ref{sec:method} describes the proposed multi-objective matheuristic, \emph{AugMathLagr}, and briefly discusses the adapted \emph{AugMathFix};
Section \ref{sec:instance_generator} introduces the random instance generator;
Section \ref{sec:experiments} presents the computational experiments {and} Section \ref{sec:conclusion} sums up the paper drawing some conclusions {and suggesting} future works.

\section{Related works}\label{sec:relatedMethods}

In this section, we review SSC management multi-objective optimization models and solution methods related to this paper. 
Due to the wide range of SSC problems, we narrowed the scope of our search to forward \typo{and} closed-loop multi-objective SSC models that aim at optimizing economic, environmental and social objectives and that integrate both \typo{strategic} and tactical planning. 
The integration of \typo{strategic} and tactical planning decisions prevents optimization models from finding sub-optimal solutions \citep{Lopes2008, Farrokhi2019}. Tactical decisions, such as the inventory and production quantities in warehouses and factories, may affect long-term decisions on where to locate these entities. They represent \typo{the} average values of these variables within the periods considered. For further details on  optimization models and solution methods for SSC management problems\typo{,} we refer to \citet{Povoa2018}.

Consider the notations for the objective functions presented in Table \ref{tab:related_labels}.
Table \ref{tab:related_models} shows SSC optimization models and solution methods published in the last 5 years.
In this table, columns ``CL'' to ``MH'' indicate whether or not the study is characterized by: closed-loop or forward supply chains (CL); inventory decisions (IQ); production quantity choices (PQ); multiple products (MP); multi-period planning horizon (PT); multiple transportation modes (MT) and multiple technologies to manufacture the products (MH).
Columns ``Eco'', ``Env'' and ``Soc'' show respectively the  economic, environmental and social indicators optimized by the model or solution method.
The references marked as ``-'' in column ``Env'' did not specify which environmental impact indicators they use in the optimization. 
Column ``Method'' presents the optimization approach to solve the corresponding problem. Column ``Type'' mentions the type of solution method the studies consider among exact, meta-heuristic and matheuristic methods. 
Finally, column \minor{``Solutions'' classifies the solution methods reported in the corresponding references as generative, \textit{a priori} or interactive.
Generative methods approximate the Pareto frontiers; \textit{a priori} methods assign preferences to the objectives; and interactive methods require the intervention of a  decision-maker.}

\begin{center}
\begin{table}[htb]
\begin{center}
\scriptsize
\caption{Notations for the economic, environmental and social functions.}
\label{tab:related_labels}
\begin{tabular}{ll}
\hline
\multicolumn{1}{c}{\textbf{Description}} & \textbf{Notation} \\ \hline
\multicolumn{2}{c}{\textbf{Economic function}}  \\ \hline
Max NPV                    & NPV               \\
Min costs                  & C                 \\
Max profit                 & P                 \\
Global after tax profit    & ATP               \\
\hline
\multicolumn{2}{c}{\textbf{Environmental function}} \\ \hline
ReCiPe 2008                        & ReCiPe            \\
Min GHG emissions or CO2 emissions & GHG               \\
Min energy consumption             & EC                \\
Min waste                          & W                 \\
Eco-indicator 99                   & EI99              \\
Environmental impact indicator     & IC                \\
Water consumption                  & WC                \\
\hline
\multicolumn{2}{c}{\textbf{Social function} } \\ \hline
Jobs created                                                  & JC                \\
Maximize justice based employment                             & EJ                \\
Workers' damage                                               & WD                \\
Workers' health                                               & WH                \\
Social Benefit Indicator                          & SB                \\
Min proximity of employees to work                            & EP                \\
Min hires and layoff to incentive stability                   & ES                \\
Max the distance between undesirable facilities and customers & UF                \\
Equity and workload in work                                   & WE                \\
Workers' safety                                               & WS                \\
Consumer risk                                                 & CR                \\
Maximizes humanitarian forces                                 & HF                \\
Min total dissatisfaction of served customer zones            & CD                \\
Maximize food security                                        & FS                \\
Charity donations or support                                  & Ch                \\
Training hours                                                & WT                \\
Community service hours                                       & CS                \\
Social responsibility                                         & SR                \\
Social impact to open plants                                  & SI                \\
Hazardous by-products                                         & HP                \\
Uncollected wasted cooking oil                                & UWCO              \\ \hline
\end{tabular}
\end{center}
\end{table}
\end{center}

\begin{table}[htb]
\scriptsize
\tabcolsep=0.05cm
\caption{{SSC optimization models and solution methods.}}
\label{tab:related_models}
\begin{tabular}{|l|l|l|l|l|l|l|l|l|l|l|l|l|l|l|}
\hline
\textbf{Ref.}                                  & \textbf{Problem}                                                                        & \textbf{CL}                 & \textbf{IQ}                 & \textbf{PQ}                  & \textbf{MP}                 & \textbf{PT}                 & \textbf{MT}                 & \textbf{MH}                 & \textbf{Eco} & \textbf{Env}                                           & \textbf{Soc}                                                & \textbf{Method}                                                                                     & \textbf{Type}                                                                       & \textbf{Solutions} \\ \hline
This             & \begin{tabular}[c]{@{}l@{}}Eletronic\\components\end{tabular}                                                                              & \cellcolor[HTML]{B7E1CD}yes & \cellcolor[HTML]{B7E1CD}yes & \cellcolor[HTML]{B7E1CD}yes  & \cellcolor[HTML]{B7E1CD}yes & \cellcolor[HTML]{B7E1CD}yes & \cellcolor[HTML]{B7E1CD}yes & \cellcolor[HTML]{B7E1CD}yes & NPV          & ReCiPe                                                 & SB                                                          & \begin{tabular}[c]{@{}l@{}}AUGMECON and\\ Lagrangian heuristic\end{tabular}                                                                                & Matheuristic                                                                        & \textcolor{black}{Generative}          \\ \hline
\cite{Yun2020}                 & Numerical                                                                               & \cellcolor[HTML]{B7E1CD}yes & \cellcolor[HTML]{F4C7C3}no  & \cellcolor[HTML]{F4C7C3}no   & \cellcolor[HTML]{F4C7C3}no  & \cellcolor[HTML]{F4C7C3}no  & \cellcolor[HTML]{B7E1CD}yes & \cellcolor[HTML]{F4C7C3}no  & C            & CO2                                                    & \begin{tabular}[c]{@{}l@{}}JC,WD,\\ WU\end{tabular}         & GA and CS                                                                                           & Meta-heuristic                                                                      & \textcolor{black}{Generative}          \\ \hline
\cite{Vafaeenezhad2019}        & Paper industry                                                                          & \cellcolor[HTML]{F4C7C3}no  & \cellcolor[HTML]{B7E1CD}yes & \cellcolor[HTML]{B7E1CD}yes  & \cellcolor[HTML]{B7E1CD}yes & \cellcolor[HTML]{B7E1CD}yes & \cellcolor[HTML]{F4C7C3}no  & \cellcolor[HTML]{B7E1CD}yes & C            & \begin{tabular}[c]{@{}l@{}}GHG,\\ EC, W\end{tabular}   & EP,ES                                                       & AUGMECON2                                                                                           & Exact                                                                               & \textcolor{black}{Generative}          \\ \hline
\cite{Tautenhain2019}          & \begin{tabular}[c]{@{}l@{}}Electronic\\components\end{tabular}                                                                              & \cellcolor[HTML]{B7E1CD}yes & \cellcolor[HTML]{B7E1CD}yes & \cellcolor[HTML]{B7E1CD}yes  & \cellcolor[HTML]{B7E1CD}yes & \cellcolor[HTML]{B7E1CD}yes & \cellcolor[HTML]{B7E1CD}yes & \cellcolor[HTML]{B7E1CD}yes & C            & ReCiPe                                                 & JC                                                          & \begin{tabular}[c]{@{}l@{}}AUGMECON and \\relaxation-based\\heuristic\end{tabular}           & Matheuristic                                                                        & \textcolor{black}{Generative}          \\ \hline
\cite{Razm2019}                & Bio-energy SC                                                                            & \cellcolor[HTML]{F4C7C3}no  & \cellcolor[HTML]{F4C7C3}no  & \cellcolor[HTML]{B7E1CD}yes  & \cellcolor[HTML]{B7E1CD}yes & \cellcolor[HTML]{B7E1CD}yes & \cellcolor[HTML]{F4C7C3}no  & \cellcolor[HTML]{B7E1CD}yes & ATP          & GHG                                                    & JC                                                          & AUGMECON                                                                                            & Exact                                                                               & \textcolor{black}{Generative}          \\ \hline
\cite{Rahimi2019}              & Numerical                                                                               & \cellcolor[HTML]{F4C7C3}no  & \cellcolor[HTML]{B7E1CD}yes & \cellcolor[HTML]{B7E1CD}yes  & \cellcolor[HTML]{B7E1CD}yes & \cellcolor[HTML]{B7E1CD}yes & \cellcolor[HTML]{F4C7C3}no  & \cellcolor[HTML]{B7E1CD}yes & P            & GHG                                                    & JC,WD                                                       & $\epsilon$-Constraint                                                                                  & Exact                                                                               & \textcolor{black}{Generative}          \\ \hline
\cite{Pourjavad2019}           & Numerical                                                                               & \cellcolor[HTML]{B7E1CD}yes & \cellcolor[HTML]{F4C7C3}no  & \cellcolor[HTML]{F4C7C3}no   & \cellcolor[HTML]{F4C7C3}no  & \cellcolor[HTML]{B7E1CD}yes & \cellcolor[HTML]{F4C7C3}no  & \cellcolor[HTML]{F4C7C3}no  & C            & -                                                      & JC                                                          & \begin{tabular}[c]{@{}l@{}}Fuzzy solution approach\\ and weighted sum\end{tabular}                  & Fuzzy                                                                               & \textit{A priori}          \\ \hline
\cite{Nobari2019}              & Numerical                                                                               & \cellcolor[HTML]{B7E1CD}yes & \cellcolor[HTML]{F4C7C3}no  & \cellcolor[HTML]{B7E1CD}yes  & \cellcolor[HTML]{B7E1CD}yes & \cellcolor[HTML]{B7E1CD}yes & \cellcolor[HTML]{F4C7C3}no  & \cellcolor[HTML]{F4C7C3}no  & P            & \begin{tabular}[c]{@{}l@{}}Green\\ degree\end{tabular} & CD                                                          & \begin{tabular}[c]{@{}l@{}}Game theory\\combined with\\NSGA-II and MOICA\end{tabular}               & Meta-heuristic                                                                      & \textcolor{black}{Generative}          \\ \hline
\cite{Martins2019}             & Food bank SC                                                                            & \cellcolor[HTML]{F4C7C3}no  & \cellcolor[HTML]{F4C7C3}no  & \cellcolor[HTML]{F4C7C3}no   & \cellcolor[HTML]{B7E1CD}yes & \cellcolor[HTML]{B7E1CD}yes & \cellcolor[HTML]{F4C7C3}no  & \cellcolor[HTML]{F4C7C3}no  & C            & FW, CO2                                                & Ch                                                     & Lexicographic ordering                                                                              & Exact                                                                               & \textit{A priori}       \\ \hline
\cite{Hajiaghaei2019} & Glass industry                                                                          & \cellcolor[HTML]{B7E1CD}yes & \cellcolor[HTML]{F4C7C3}no  & \cellcolor[HTML]{B7E1CD}yes  & \cellcolor[HTML]{F4C7C3}no  & \cellcolor[HTML]{F4C7C3}no  & \cellcolor[HTML]{F4C7C3}no  & \cellcolor[HTML]{F4C7C3}no  & C            & -                                                      & JC,WD                                                       & \begin{tabular}[c]{@{}l@{}}Several meta-heuristic,\\ among them ICA,\\ GA, SA and TS\end{tabular} & Meta-heuristic                                                                      & \textcolor{black}{Generative}          \\ \hline
\cite{Govindan2019}            & Numerical                                                                               & \cellcolor[HTML]{F4C7C3}no  & \cellcolor[HTML]{F4C7C3}no  & \cellcolor[HTML]{B7E1CD}yes  & \cellcolor[HTML]{B7E1CD}yes & \cellcolor[HTML]{F4C7C3}no  & \cellcolor[HTML]{B7E1CD}yes & \cellcolor[HTML]{B7E1CD}yes & C            & GHG                                                    & JC,WE                                                       & \begin{tabular}[c]{@{}l@{}}Swarm intelligence\\ algorithms and VNS\end{tabular}                     & Meta-heuristic                                                                      & \textcolor{black}{Generative}          \\ \hline
\cite{Gonela2019}              & \begin{tabular}[c]{@{}l@{}}Electricity\\ generation\end{tabular}                         & \cellcolor[HTML]{F4C7C3}no  & \cellcolor[HTML]{B7E1CD}yes & \cellcolor[HTML]{B7E1CD}yes  & \cellcolor[HTML]{F4C7C3}no  & \cellcolor[HTML]{B7E1CD}yes & \cellcolor[HTML]{F4C7C3}no  & \cellcolor[HTML]{F4C7C3}no  & C            & GHG                                                    & JC                                                          & AUGMECON                                                                                            & Exact                                                                               & \textcolor{black}{Generative}          \\ \hline
\cite{Farrokhi2019}        & Numerical                                                                               & \cellcolor[HTML]{B7E1CD}yes & \cellcolor[HTML]{F4C7C3}no  & \cellcolor[HTML]{F4C7C3}no   & \cellcolor[HTML]{F4C7C3}no  & \cellcolor[HTML]{F4C7C3}no  & \cellcolor[HTML]{B7E1CD}yes & \cellcolor[HTML]{B7E1CD}yes & C            & CO2                                                    & UF                                                          & NSGA-II and SA                                                                                      & Meta-heuristic                                                                      & \textcolor{black}{Generative}          \\ \hline
\cite{Darbari2019}             & Numerical                                                                               & \cellcolor[HTML]{B7E1CD}yes & \cellcolor[HTML]{F4C7C3}no  & \cellcolor[HTML]{B7E1CD}yes  & \cellcolor[HTML]{B7E1CD}yes & \cellcolor[HTML]{F4C7C3}no  & \cellcolor[HTML]{B7E1CD}yes & \cellcolor[HTML]{F4C7C3}no  & P            & CO2                                                    & \begin{tabular}[c]{@{}l@{}}JC,WT,\\CS,Ch\end{tabular} & \begin{tabular}[c]{@{}l@{}}GP and fuzzy\\solution approach\end{tabular}              & Fuzzy                                                                               & \textit{A priori}          \\ \hline
\cite{Tsao2018}                & Numerical                                                                               & \cellcolor[HTML]{F4C7C3}no  & \cellcolor[HTML]{F4C7C3}no  & \cellcolor[HTML]{B7E1CD}yes  & \cellcolor[HTML]{F4C7C3}no  & \cellcolor[HTML]{F4C7C3}no  & \cellcolor[HTML]{F4C7C3}no  & \cellcolor[HTML]{B7E1CD}yes & C            & CO2                                                    & JC,WD                                                       & Fuzzy solution approach                                                                             & Fuzzy                                                                               & Interactive      \\ \hline
\cite{Sahebjamnia2018}         & Numerical                                                                               & \cellcolor[HTML]{B7E1CD}yes & \cellcolor[HTML]{F4C7C3}no  & \cellcolor[HTML]{B7E1CD}yes  & \cellcolor[HTML]{B7E1CD}yes & \cellcolor[HTML]{F4C7C3}no  & \cellcolor[HTML]{F4C7C3}no  & \cellcolor[HTML]{B7E1CD}yes & C            & ReCiPe                                                 & \begin{tabular}[c]{@{}l@{}}JC,WS\\ HP\end{tabular}          & \begin{tabular}[c]{@{}l@{}}Meta-heuristics and\\ $\epsilon$ constraint\end{tabular}                    & Meta-heuristic                                                                      & \textcolor{black}{Generative}          \\ \hline
\cite{Rezaei2018}              & Numerical                                                                               & \cellcolor[HTML]{B7E1CD}yes & \cellcolor[HTML]{F4C7C3}no  & \cellcolor[HTML]{F4C7C3}no   & \cellcolor[HTML]{F4C7C3}no  & \cellcolor[HTML]{F4C7C3}no  & \cellcolor[HTML]{B7E1CD}yes & \cellcolor[HTML]{B7E1CD}yes & C            & -                                                      & JC                                                          & CS                                                                                                  & Meta-heuristic                                                                      & \textcolor{black}{Generative}          \\ \hline
\cite{Rabbani2018}             & \begin{tabular}[c]{@{}l@{}}Switchgrass-based\\ bioenergy SC\end{tabular}                & \cellcolor[HTML]{F4C7C3}no  & \cellcolor[HTML]{B7E1CD}yes & \cellcolor[HTML]{B7E1CD}yes & \cellcolor[HTML]{B7E1CD}yes & \cellcolor[HTML]{B7E1CD}yes & \cellcolor[HTML]{B7E1CD}yes & \cellcolor[HTML]{F4C7C3}no  & C            & GHG                                                    & JC                                                          & \begin{tabular}[c]{@{}l@{}}AUGMECON and\\ TOPSIS\end{tabular}                                       & Exact                                                                               & \textcolor{black}{Generative}          \\ \hline
\cite{Pourjavad2018}           & Numerical                                                                               & \cellcolor[HTML]{B7E1CD}yes & \cellcolor[HTML]{F4C7C3}no  & \cellcolor[HTML]{F4C7C3}no   & \cellcolor[HTML]{F4C7C3}no  & \cellcolor[HTML]{B7E1CD}yes & \cellcolor[HTML]{F4C7C3}no  & \cellcolor[HTML]{F4C7C3}no  & C            & -                                                      & JC                                                          & NSGA-II                                                                                             & Meta-heuristic                                                                      & \textit{A priori}          \\ \hline
\cite{Nobari2018}              & Numerical                                                                               & \cellcolor[HTML]{B7E1CD}yes & \cellcolor[HTML]{B7E1CD}yes & \cellcolor[HTML]{B7E1CD}yes  & \cellcolor[HTML]{B7E1CD}yes & \cellcolor[HTML]{B7E1CD}yes & \cellcolor[HTML]{F4C7C3}no  & \cellcolor[HTML]{F4C7C3}no  & P            & GHG                                                    & JC                                                          & MOICA                                                                                               & Meta-heuristic                                                                      & \textcolor{black}{Generative}          \\ \hline
\cite{Mota2018}                & \begin{tabular}[c]{@{}l@{}}Electronic\\components\end{tabular}                          & \cellcolor[HTML]{B7E1CD}yes & \cellcolor[HTML]{B7E1CD}yes & \cellcolor[HTML]{B7E1CD}yes  & \cellcolor[HTML]{B7E1CD}yes & \cellcolor[HTML]{B7E1CD}yes & \cellcolor[HTML]{B7E1CD}yes & \cellcolor[HTML]{B7E1CD}yes & NPV          & ReCiPe                                                 & SB                                                          & $\epsilon$-Constraint                                                                                  & Exact                                                                               & \textit{A priori}          \\ \hline
\cite{How2018}                 & Biomass SC                                                                              & \cellcolor[HTML]{F4C7C3}no  & \cellcolor[HTML]{F4C7C3}no  & \cellcolor[HTML]{B7E1CD}yes  & \cellcolor[HTML]{B7E1CD}yes & \cellcolor[HTML]{F4C7C3}no  & \cellcolor[HTML]{B7E1CD}yes & \cellcolor[HTML]{B7E1CD}yes & P            & IC                                                     & \begin{tabular}[c]{@{}l@{}}JC,WH,\\ WS\end{tabular}         & \begin{tabular}[c]{@{}l@{}}PCA and \\ AHP\end{tabular}                                              & Decomposition                                                                       & \textit{A priori}          \\ \hline
\cite{Ghaderi2018}             & \begin{tabular}[c]{@{}l@{}}Switchgrass-based \\ bioethanol SC\end{tabular}              & \cellcolor[HTML]{F4C7C3}no  & \cellcolor[HTML]{B7E1CD}yes & \cellcolor[HTML]{B7E1CD}yes  & \cellcolor[HTML]{B7E1CD}yes & \cellcolor[HTML]{B7E1CD}yes & \cellcolor[HTML]{B7E1CD}yes & \cellcolor[HTML]{B7E1CD}yes & C            & ReCiPe                                                 & SR        
& \begin{tabular}[c]{@{}l@{}}fuzzy solution method\end{tabular}               & \begin{tabular}[c]{@{}l@{}}Fuzzy\\Interactive\end{tabular}                                                                               & \textcolor{black}{Generative}         \\ \hline
\cite{Eskandari2018} & Blood bank SC                                                                           & \cellcolor[HTML]{F4C7C3}no  & \cellcolor[HTML]{B7E1CD}yes & \cellcolor[HTML]{F4C7C3}no   & \cellcolor[HTML]{B7E1CD}yes & \cellcolor[HTML]{B7E1CD}yes & \cellcolor[HTML]{B7E1CD}yes & \cellcolor[HTML]{F4C7C3}no  & C            & -                                                      & JC                                                          & $\epsilon$-Constraint                                                                                  & Exact                                                                               & \textcolor{black}{Generative}          \\ \hline
\cite{Chavez2018}              & \begin{tabular}[c]{@{}l@{}}Biofuel SC from\\ coffee crop residues\end{tabular}          & \cellcolor[HTML]{F4C7C3}no  & \cellcolor[HTML]{B7E1CD}yes & \cellcolor[HTML]{B7E1CD}yes  & \cellcolor[HTML]{B7E1CD}yes & \cellcolor[HTML]{B7E1CD}yes & \cellcolor[HTML]{F4C7C3}no  & \cellcolor[HTML]{F4C7C3}no  & NPV          & CO2                                                    & JC,FS                                                       & $\epsilon$-Constraint                                                                                  & Exact                                                                               & \textit{A priori}          \\ \hline
\cite{Varsei2017}              & Wine industry                                                                           & \cellcolor[HTML]{F4C7C3}no  & \cellcolor[HTML]{F4C7C3}no  & \cellcolor[HTML]{F4C7C3}no   & \cellcolor[HTML]{F4C7C3}no  & \cellcolor[HTML]{F4C7C3}no  & \cellcolor[HTML]{B7E1CD}yes & \cellcolor[HTML]{F4C7C3}no  & C            & CO2                                                    & SI                                                          & AUGMECON                                                                                            & Exact                                                                               & \textcolor{black}{Generative}          \\ \hline
\cite{Sabegh2017}              & \begin{tabular}[c]{@{}l@{}}Healthcare SC\\ for natural\\ disaster response\end{tabular} & \cellcolor[HTML]{F4C7C3}no  & \cellcolor[HTML]{F4C7C3}no  & \cellcolor[HTML]{B7E1CD}yes  & \cellcolor[HTML]{F4C7C3}no  & \cellcolor[HTML]{F4C7C3}no  & \cellcolor[HTML]{F4C7C3}no  & \cellcolor[HTML]{B7E1CD}yes & C            & -                                                      & HF                                                          & \begin{tabular}[c]{@{}l@{}}Neural network,\\ GA and PSO\end{tabular}                                & \begin{tabular}[c]{@{}l@{}}Neural network,\\Meta-heuristic\\and GA\end{tabular} & \textit{A priori}          \\ \hline
\cite{Osmani2017}              & Bioethanol SC                                                                           & \cellcolor[HTML]{F4C7C3}no  & \cellcolor[HTML]{F4C7C3}no  & \cellcolor[HTML]{B7E1CD}yes  & \cellcolor[HTML]{B7E1CD}yes & \cellcolor[HTML]{B7E1CD}yes & \cellcolor[HTML]{F4C7C3}no  & \cellcolor[HTML]{B7E1CD}yes & P            & GHG                                                    & JC                                                          & \begin{tabular}[c]{@{}l@{}}AUGMECON\\and Benders\\ decomposition\end{tabular}                      & Matheuristic                                                                        & \textcolor{black}{Generative}          \\ \hline
\cite{Jafari2017}              & Textile industry                                                                        & \cellcolor[HTML]{B7E1CD}yes & \cellcolor[HTML]{F4C7C3}no  & \cellcolor[HTML]{B7E1CD}yes  & \cellcolor[HTML]{B7E1CD}yes & \cellcolor[HTML]{F4C7C3}no  & \cellcolor[HTML]{B7E1CD}yes & \cellcolor[HTML]{F4C7C3}no  & C            & WC                                                     & EJ                                                          & MOVDO                                                                                               & Meta-heuristic                                                                      & \textcolor{black}{Generative}          \\ \hline
\cite{Arampantzi2017}          & \begin{tabular}[c]{@{}l@{}}Global manufacturer\end{tabular}                                                                     & \cellcolor[HTML]{F4C7C3}no  & \cellcolor[HTML]{B7E1CD}yes & \cellcolor[HTML]{B7E1CD}yes  & \cellcolor[HTML]{F4C7C3}no  & \cellcolor[HTML]{B7E1CD}yes & \cellcolor[HTML]{B7E1CD}yes & \cellcolor[HTML]{F4C7C3}no  & C            & GHG, W                                                 & JC*                                                         & \begin{tabular}[c]{@{}l@{}}GP and\\ $\epsilon$-Constraint\end{tabular}                                 & Exact                                                                               & \textcolor{black}{Generative}          \\ \hline
\cite{Zhang2016robust}               & Biodisel SC                                                                             & \cellcolor[HTML]{F4C7C3}no  & \cellcolor[HTML]{F4C7C3}no  & \cellcolor[HTML]{F4C7C3}no   & \cellcolor[HTML]{F4C7C3}no  & \cellcolor[HTML]{F4C7C3}no  & \cellcolor[HTML]{F4C7C3}no  & \cellcolor[HTML]{F4C7C3}no  & P            & CO2                                                    & UWCO                                                        & \begin{tabular}[c]{@{}l@{}}Decomposition\\ GA and\\ weighted sum\end{tabular}                       & \begin{tabular}[c]{@{}l@{}}Heuristic/\\ Decomposition\end{tabular}                  & \textit{A priori}       \\ \hline
\cite{Miret2016}               & Bioethanol SC                                                                           & \cellcolor[HTML]{F4C7C3}no  & \cellcolor[HTML]{B7E1CD}yes & \cellcolor[HTML]{B7E1CD}yes  & \cellcolor[HTML]{B7E1CD}yes & \cellcolor[HTML]{B7E1CD}yes & \cellcolor[HTML]{F4C7C3}no  & \cellcolor[HTML]{B7E1CD}yes & C            & Eco-costs                                              & JC                                                          & GP                                                                                                  & Exact                                                                               & \textit{A priori}          \\ \hline
\cite{Bairamzadeh2016}         & Bioethanol SC                                                                           & \cellcolor[HTML]{F4C7C3}no  & \cellcolor[HTML]{B7E1CD}yes & \cellcolor[HTML]{B7E1CD}yes  & \cellcolor[HTML]{B7E1CD}yes & \cellcolor[HTML]{B7E1CD}yes & \cellcolor[HTML]{B7E1CD}yes & \cellcolor[HTML]{B7E1CD}yes & P            & EI99                                                   & JC                                                          & \begin{tabular}[c]{@{}l@{}}Robust possibilistic\\ programming\end{tabular}                          & \begin{tabular}[c]{@{}l@{}}Exact/\\Robust\end{tabular}                                                                     & Interactive     \\ \hline
\cite{Mota2015}                & Battery industry                                                                        & \cellcolor[HTML]{B7E1CD}yes & \cellcolor[HTML]{B7E1CD}yes & \cellcolor[HTML]{F4C7C3}no  & \cellcolor[HTML]{B7E1CD}yes & \cellcolor[HTML]{B7E1CD}yes & \cellcolor[HTML]{F4C7C3}no  & \cellcolor[HTML]{F4C7C3}no  & C            & ReCiPe                                                 & SB                                                          & $\epsilon$-Constraint                                                                                  & Exact                                                                               & \textit{A priori}       \\ \hline
\cite{Boukherroub2015}         & Lumber industry                                                                         & \cellcolor[HTML]{F4C7C3}no  & \cellcolor[HTML]{B7E1CD}yes & \cellcolor[HTML]{B7E1CD}yes  & \cellcolor[HTML]{B7E1CD}yes & \cellcolor[HTML]{B7E1CD}yes & \cellcolor[HTML]{B7E1CD}yes & \cellcolor[HTML]{B7E1CD}yes & C            & GHG                                                    & EP                                                          & Weighted sum                                                                                        & Exact                                                                               & \textit{A priori}          \\ \hline
\cite{Pishvaee2014}            & \begin{tabular}[c]{@{}l@{}}Medical needle\\ and syringe SC\end{tabular}                 & \cellcolor[HTML]{B7E1CD}yes & \cellcolor[HTML]{F4C7C3}no  & \cellcolor[HTML]{B7E1CD}yes  & \cellcolor[HTML]{B7E1CD}yes & \cellcolor[HTML]{F4C7C3}no  & \cellcolor[HTML]{F4C7C3}no  & \cellcolor[HTML]{B7E1CD}yes & C            & ReCiPe                                                 & JC,CR                                                       &  \begin{tabular}[c]{@{}l@{}}Benders decomposition\\and fuzzy solution\\
approach\end{tabular}                                                                              & Matheuristic                                                                        & \textcolor{black}{Generative}          \\ \hline
\cite{Devika2014}              & Glass industry                                                                          & \cellcolor[HTML]{B7E1CD}yes & \cellcolor[HTML]{F4C7C3}no  & \cellcolor[HTML]{B7E1CD}yes  & \cellcolor[HTML]{F4C7C3}no  & \cellcolor[HTML]{F4C7C3}no  & \cellcolor[HTML]{F4C7C3}no  & \cellcolor[HTML]{B7E1CD}yes & C            & GHG, W                                                 & JC,WD                                                       & ICA based and VNS                                                                                   & Meta-heuristic                                                                      & \textcolor{black}{Generative}          \\ \hline
\end{tabular}
\end{table}

One may observe that only the closed-loop models presented in \cite{Tautenhain2019} and \cite{Mota2018} include all the \typo{strategic} and tactical decisions listed in Table \ref{tab:related_models}. 
The model introduced in \citep{Tautenhain2019} also considers \typo{the} carbon market. \citet{Mota2018}, on the other hand, model cross-docking transportation through airports and seaports and further explore transportation decisions such as the number of vehicles and trips required to deliver products. 
Despite including all the \typo{strategic} and tactical decisions presented in Table \ref{tab:related_models}, \citet{Ghaderi2018, Bairamzadeh2016, Boukherroub2015} do not address the reverse flow in their models. The variety of decisions each model takes into account suggests that generic models should include as many decisions as possible in order to represent a wider range of supply chains.

The minimization of costs and profit maximization are the most commonly employed economic objective functions in the reviewed SSC models. \minor{Most studies minimize the environmental indicators related to Greenhouse Gases (GHG) or specific CO$_2$ emissions, despite the European Commission putting forward the recommendation that Lyfe Cycle Approaches (LCA) are the best available methodology to evaluate the environmental impact \citep{EC2003}.} 
\typo{Other} case-oriented indicators are the minimization of energy and water consumption besides the waste resulting from SC operations.

Social indicators are difficult to quantify \citep{Eskandarpour2015}, being commonly linked to the planning situation. The majority of references in Table~\ref{tab:related_models}  regard  the socio-economic development of regions by optimizing job creation indicators.
\citet{Mota2015,Mota2018} designed the Social Benefit Indicator (SB) to favor the job creation in less developed countries. 
\minor{SB assesses the development of a country according to its Gross Domestic Product (GDP) \citep{Brezina2011}, which is a widely accepted index to measure the economic activity of a nation. In this context, countries with less industrial development have a lower GDP and, according to the defined SB, are less developed.}
The authors recommend, however, its use only in planning situations that preferably involve hiring and do not require lay-offs. Another concern about job creation is related to working conditions. To this end, other widely employed
indicators take into account \typo{workers'} safety and health. When combined, these indicators can help \typo{to prevent} the creation of jobs in poor working conditions.

A vast body of literature on {SSC management} multi-objective solution methods is based on scenario analysis, the $\epsilon$-Constraint or weighting aggregation methods and heuristics or matheuristics.

To approach a multi-objective minimization problem, the $\epsilon$-Constraint method \citep{Marglin1967,Haimes1971} optimizes one of its objective functions and adds constraints to assign upper bounds for the other objective functions.
In particular, \citet{Mavrotas2009} introduced the Augmented $\epsilon$-Constraint method (AUGMECON) that solves a sequence of mono-objective problems. Later, \citet{Mavrotas2013} proposed AUGMECON2 as an improved version of AUGMECON.
Examples of works that employ AUGMECON to approach multi-objective {SSC management problems} are found in \citep{Rabbani2018,Osmani2017,Vafaeenezhad2019}.

 Exact-based methods often employ mixed-integer linear solvers \typo{such as} CPLEX
 to solve the resulting mono-objective problems. In practice, these methods are computationally costly to approach most SSC problems. In such cases, AUGMECON and AUGMECON2 are impracticable to approximate the Pareto frontier. To overcome this limitation, \citet{Mota2018} optimized the SSC problem using \typo{scenario} analysis.

 Goal Programming (GP) is an exact-based method used in a number of studies presented in Table \ref{tab:related_models}.
 GP methods aim at finding solutions by minimizing the deviance from pre-defined values for  each objective function, also known as goals.
 In the references indicated in Table \ref{tab:related_models},  there are both generative and \textit{a priori} versions of GP methods. 
 When dealing with uncertainty in SSC, fuzzy-based solution approaches are commonly employed as \textit{a priori} or interactive methods.

Several multi-objective algorithms to heuristically solve SSCs  have  been introduced.  Among them, one can find Non-dominated Sorting Genetic Algorithm (NSGA-II) \citep{Pourjavad2018,Farrokhi2019,Nobari2019}, Genetic Algorithms (GA) \citep{Zhang2016robust,Sabegh2017,Hajiaghaei2019,Yun2020},  Variable Neighborhood Search (VNS) \citep{Devika2014, Govindan2019}, Imperative Competitive Algorithm (ICA) \citep{Devika2014, Hajiaghaei2019}, multi-objective ICA (MOICA) \citep{Nobari2018, Nobari2019}, Cuckoo Search (CS) \citep{Rezaei2018, Yun2020},   Multi-Objective Vibration Damping Optimization (MOVDO) \cite{Jafari2017}, Particle Swarm Optimization (PSO) \citep{Sabegh2017, Govindan2019}, Simulated Annealing (SA) \citep{Hajiaghaei2019,Farrokhi2019} and Tabu Search (TS) \citep{Hajiaghaei2019}. Other strategies such as Neural Networks \citep{Sabegh2017}, Principal Component Analysis (PCA) and Analytical Hierarchical Process (AHP) \citep{How2018} can also be found in the literature.

Matheuristics are powerful solution methods that hybridize heuristic methods with mathematical programming techniques.  A few matheuristics can be found in the literature to approach multi-objective SSC problems. The matheuristic proposed by \citet{Osmani2017} combined AUGMECON with a Benders decomposition  to find Pareto solutions for a SSC management problem. \citet{Zhang2016robust} also introduced a matheuristic that employs a Benders decomposition to split a SSC problem into two subproblems associated with the \typo{strategic} and tactical decisions.
\citet{Tautenhain2019} introduced a multi-objective matheuristic founded on AUGMECON2, called \emph{AugMathFix}, which iteratively solves mono-objective SSC management  problems with strategic constraints relaxed. In this case, for the method to find a feasible solution to the original problem in a reduced computational time,  the variables related to the tactical decisions are fixed at the values of the solution of the relaxed problem. Moreover, {at each iteration} a local search strategy is applied to the feasible solutions to enhance the quality of the solutions found.

As presented in Table \ref{tab:related_models}, with the exception of \emph{AugMathFix}, the existing SSC heuristics do not account for all the decisions considered in the model introduced in \citep{Mota2018}.
However, the resulting integrated SSC problem has a significant amount of integer variables, forcing the authors to simplify its solution through scenario analysis \cite{Mota2018}.
Optimization solvers like CPLEX  have a better performance in solving  mixed-integer programs when they have a reduced number of integer variables. In this context, matheuristics that decompose the original problem into simpler ones or that relax binary and integer variables, for example, may speed up mixed-integer solvers in finding solutions.

Lagrangian decomposition methods, in particular, are potential candidates to solve SSC optimization problems since they allow us to simplify complicated constraints, such as those including binary variables, for example. These methods presented good results in finding solutions for related SC problems. \citet{Heidari2018} and \citet{Yousefi2017} studied  the Lagrangian decomposition of multi-objective supply chain {problems.} \citet{Heidari2018} employed the $\epsilon$-Constraint method to transform the multi-objective SSC management problem into a sequence of mono-objective problems. 
The resulting mono-objective problems were then solved by a Lagrangian heuristic. The pair of constraints to be relaxed was the one that, when removed from the problems, resulted in the lowest computational running times.

In the introduced mathematical formulation, \citet{Yousefi2017} did not consider environmental issues in any of the objective functions, however, they included recycling centers. The authors employed a weighting  strategy to aggregate the objective functions into a mono-objective function and then relaxed the most complex constraints of the problem to then solve it using an optimization solver.

Lagrangian-based heuristics {are extensively studied to approach} mono-objective supply {chain} problems \citep{Eskigun2005,Lidestam2011,Elhedhli2012, Zhang2014,Rafie2018}.
\citet{Eskigun2005} reduced the original capacitated network design problem into simpler and independent subproblems {to approach the {design of vehicle distribution centers of several instances constructed from industrial data} by a Lagrangian heuristic}.  
{In the Lagrangian heuristic proposed by} \citet{Lidestam2011}{, the authors} decomposed the supply chain planning problem into two different subproblems associated with different stages of  cellulose production.
\citet{Elhedhli2012} studied a supply chain network design problem that aims to minimize environmental costs due to carbon dioxide emissions. The authors employed a Langragian relaxation by decomposing the original problem according to entity types and  warehouse site.
\citet{Zhang2014} approached a supply chain problem
that required  distribution centers between (i) suppliers and factories, and (ii) factories and customers.
 The authors relaxed the constraints which ensured that each customer and supplier was assigned to only one distribution center in the Lagrangian heuristic they suggested. \citet{Rafie2018} studied a supply chain of perishable products that takes into consideration fuel consumption and product wastage. The Lagrangian heuristic presented by the authors considers the relaxation of the vehicle capacity constraints, a nonlinear constraint and a constraint related to the allocation of distribution centers to customers.

To the best of our knowledge, this is the first Lagrangian-based method to solve an integrated SSC problem that considers all the dimensions of sustainability, as presented in Table~\ref{tab:related_models}.

 The next section describes the SSC management problem, the focus of this paper.

\section{Problem description} \label{sec:targetModel}

The studied {(SSC)} is composed of suppliers, factories, warehouses, customers, airports and seaports, and is tailored to support multi-period planning. 
Moreover, there are three types of items in the SSC: raw materials, final products and recovered products.

\begin{figure}[htb] 
 \includegraphics[width=\linewidth]{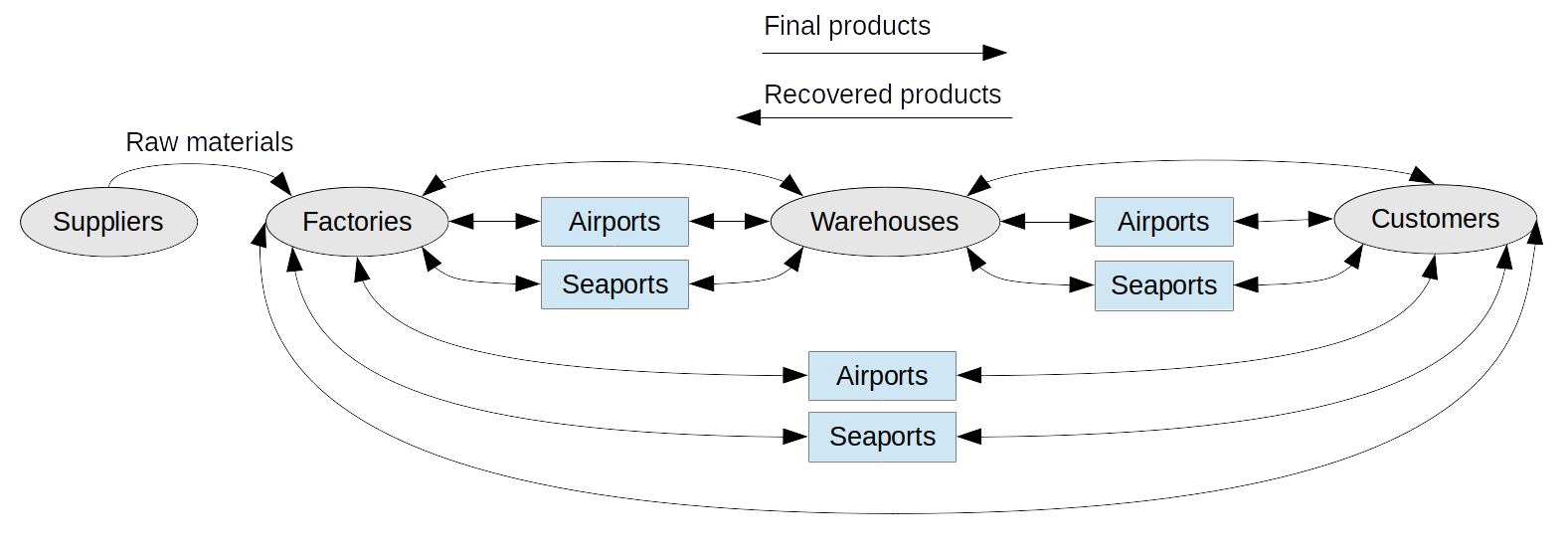}
 \caption{An illustration with a general representation of the studied SSC.}
 \label{fig:ssc_mota}
\end{figure}

Figure \ref{fig:ssc_mota} illustrates a flow network whose nodes are the entities of the SSC, and the arcs represent the flow of goods between the sites. 
In this figure, the rectangles represent airports and seaports whereas the ellipses represent the remaining entities.
The double-sided arrows indicate that the corresponding arcs represent the flow of both final products and recovered products between the sites. The rightwards arrow represents the arc indicating the flow of raw materials from suppliers to factories.

Factories employ production technologies to manufacture final products from a bill of raw materials and {may} also use remanufacturing technologies to reuse products. The problem allows {the storage} of final products sent from factories to warehouses. 
Factories and warehouses can ship final products to meet  customer demands, indicated by the arc linking these entities. Customers return used products to factories and warehouses after the end of their lifetime, which is defined as a parameter measured in \typo{the} number of periods of the planning horizon in the model. Warehouses can then return the recovered products {to the} factories. 

Land transportation is responsible for transporting raw materials from suppliers to factories.
The possible transportation modes between factories, warehouses and customers are by land, air or sea. 
Trucks carry items from and to entities through land transportation, {whereas} air and sea transportation are only allowed from/to airports or seaport hubs.

The goals of the studied SSC management are related to the triple bottom line:  (i) the  maximization of the Net Present Value (NPV) as the economic function; (ii) the minimization of the environmental impact evaluated by ReCiPe 2008 \citep{Goedkoop2009} as the environmental function; and (iii) the maximization of the Social Benefit Indicator (SB) \citep{Mota2018} as the social function. This social indicator favors the creation of jobs in less developed countries.
\minor{The Gross Domestic Product (GDP) \citep{Brezina2011} was used to assess and compare the country's industrial base development.}

\minor{SB  can also be adjusted according to the context of the supply chain under study. In analyzing supply chains at the regional level, for example, the development of regions might be assessed through unemployment rates or population densities \citep{Mota2015}. }

In this paper, in order to prevent job creation in poor working environments, we also considered an additional indicator linked to the labor conditions in entities in the social objective. 
Moreover, we included the impacts of suppliers in the environment function as well.

In the studied SSC, the set of periods is given by $T=\{1,2,\dots,|T|\}$ {and} the set of entities by $I= I_{sup} \cup I_{f} \cup I_{w} \cup I_{c} \cup I_{air} \cup I_{port}$, where $I_{sup} , I_{f} , I_w, I_{c} , I_{air}$ and $I_{port}$ are, respectively, the sets of suppliers, factories, warehouses, customers, airports and seaports. The set of items is given by $M = M_{rm} \cup M_{fp} \cup M_{rp}$, where $M_{rm}, M_{fp}$ and $M_{rp}$ are the sets of raw materials, final products and recovered products, respectively. 
In this SSC, intermediate products required in the bill to manufacture other products are also defined as final products.
The set of transportation modes is given by  $A = A_{truck} \cup A_{plane} \cup A_{boat}$, where $A_{truck} , A_{plane}$ and $A_{boat}$ are the set of trucks, airplanes and ships, respectively. The {technologies set} is given by {$G=G_{prod} \cup G_{rem}$, where $G_{prod}$ and $G_{rem}$} are, respectively,  the production and the \typo{remanufacturing} technologies sets. Moreover, consider $H_{prod}=\{(m,g) \mbox{ such that } m \in M_{fp} \mbox{ can be manufactured using technology } g \in G_{prod}\}$ and $H_{rem}=\{(m,g) \mbox{ such that }  m \in M_{fp} \mbox{ can be remanufactured using technology } g \in G_{rem}\}$. 
In Section 1 of the Supplementary Material, we present a brief discussion on how this formulation can be adapted to a more general network.

  \minor{This model assumes data regarding economic, environmental and social criteria are available. We acknowledge, however, that the access  of information in the supply chains is a problem, not often because it does not exist but because it is not treated. So, from our experience, a large period of time must be assigned to the identification of the required data. Also and in particular in what concerns the environmental and social data, assumptions must be often taken. For instance, in the environmental component, similar operations may be taken as the basis to feed certain supply chain operations that do not have available data. This need implies that a sensitivity analysis should be performed on the parameters more subject to uncertainty. This type of work allows addressing the lack of information by looking at which information is more important, that is, which one has a more significant impact on strategic and tactical decisions.
 }

The decisions associated with the management of the studied SSC relate to the amount of raw materials acquired from each supplier; the opening of factories and warehouses and their capacities; production and remanufacturing technologies assigned to factories and recycling centers, respectively; 
amount of items produced and remanufactured in factories; storage levels in warehouses; and shipment of items between entities using different transportation modes. 
Table \ref{tab:decision_vars} presents {a} complete list of decision variables. Tables \ref{tab:params_constraints} and \ref{tab:params_objectives} present the parameters \typo{associated with} the model constraints and objective functions, respectively.

\begin{table}[htb]
\centering
\caption{Decision variables of the  SSC formulation introduced in \citep{Mota2018}.}
\tabcolsep=0.05cm
\small
\label{tab:decision_vars}
\begin{tabular}{cl}
\hline
\multicolumn{1}{c}{\textbf{Variable}} & \multicolumn{1}{c}{\textbf{Description}} \\[0.5ex] \hline
\multicolumn{2}{c}{Continuous decision variables} \\ \hline
 $S_{mit}$ 			&  Amount of product $m$ stocked in entity $i \in I_f \cup I_w$ in time period $t \in T$. \\ [0.5ex]
 \multirow{2}{*}{$X_{maijt}$} &  {Amount} of item  $m \in M$ transported from entity $i \in I$ to entity $j \in I$. \\
				& 	by transportation mode $a$ in time period $t \in T$. \\ [0.5ex]
 $P_{mgit}$ 			&  Amount of product $m \in M$ produced by technology $g \in G_{prod}$ in factory $i \in I_f$ in time period $t \in T$. \\ [0.5ex]
 $R_{mgit}$ 			&  Amount of product $m$ remanufactured by technology $g \in G_{rem}$ in factory $i \in I_f$ in time period $t \in T$. \\ [0.5ex]
 $YC_{i}$ 			&  Capacity of entity $i \in I_f \cup I_w$. \\ [0.5ex]
 $YCT_{it}$ 			&  Effective use of capacity in entity $i \in I_f \cup I_w$ in time period $t \in T$. \\ [0.5ex]
 \multirow{1}{*}{$KT_{ait}$} 	&  Upper bound to the number of transportation modes $a$ from entity $i \in I$ 
				to another in time period $t \in T$. \\ [0.5ex] \hline
 \multicolumn{2}{c}{Integer decision variables} \\ \hline
 $K_{ai}$			& Number of transportation modes $a \in A$ {used to transport products from entity $i \in I$.} \\
 $Q_{aijt}$			& Number of trips from entity $i \in I$ to $j \in I$ by transportation mode $a$ in time period $t \in T$. \\[0.5ex] \hline
 \multicolumn{2}{c}{Binary decision variables} \\ \hline
 $Y_i$				& value $1$ indicates that entity $i \in I$ is installed and $0$, otherwise. \\ [0.5ex]
 \multirow{2}{*}{$Z_{gmi}$}	& value $1$ indicates that the technology $g$ is selected to produce product $m$ in factory $i \in I_f$ \\
				& and $0$, otherwise.\\ \hline
 \end{tabular}
\end{table}

\begin{table}[htb]
\centering
\caption{{Data parameters for the constraints of the SSC formulation introduced in \citep{Mota2018}.}}
\footnotesize
\label{tab:params_constraints}
\begin{tabular}{cl}
\hline
\textbf{Parameter} &  \multicolumn{1}{c}{\textbf{Description}}	\\ [0.5ex]
\hline 
\multicolumn{2}{c}{Demands and specification of items} \\ \hline
$dmd_{mit}$ & {Demand of customer $i\in I_c$ for product $m \in M_{fp}$ in a period $t \in T$.}  \label{eq:ger_dc} \\   [0.5ex]
$BOM^{prod}_{mng}$  & {Bill of materials $m\in M_{rm}$ to produce a product $n \in M_{fp}$} \label{eq:bomprod} \\
  & {using technology $g \in G_{prod}$.}  \\  [0.5ex]
$BOM^{rem}_{mn}$ & {Bill of recovery products $m \in M_{rp}$ to remanufacture a final product $n \in M_{fp}$.} \label{eq:bomrem} \\  [0.5ex]
$pw_m$ & {Weight of a raw material, final product or} \label{eq:pw} recovery product $m \in M_{rm} \cup M_{fp} \cup M_{rp}$ ($kg$).  \\ [0.5ex]
$apu_m$ & {Necessary area per unit of $m \in M_{rm} \cup M_{fp} \cup M_{rp}$ ($m^2$).}	 \label{eq:apu} 	\\ [0.5ex]
$retF_m$ & {Return rate of a recovery product $m \in M_{rp}$.}		\\ [0.5ex]
{$LF$} & {Lifetime  in time periods of final products.}  \\ [0.5ex]
%
%
\multicolumn{2}{c}{Maximum and minimum capacities of entities} \\ \hline
$ic^{max}_{mi}$ &  {Storage capacity for a product $M_{fp}$ in a entity $i \in I_f \cup I_w$.}	 \label{eq:gen_ic} \\ [0.5ex]
$ic^{min}_{mi}$ & {Minimum stock for a product $m \in M_{fp}$ in a entity $i \in I_f \cup I_w$}	\\ [0.5ex]
$pc^{max}_{g}$ & {Maximum capacity of a technology $g \in G$.} \label{eq:gen_capMax}	\\ [0.5ex]
$pc^{min}_{g}$ & {Minimum use required {by} a technology $g \in G$.}	\\ [0.5ex]
$ea^{max}_{i}$ & {Maximum installation area of a entity $i \in I_f \cup I_w$ $(m^2)$.} \label{eq:icmax}\\ [0.5ex]
$ea^{min}_{i}$ & {Minimum installation area of a entity $i \in I_f \cup I_w$ $(m^2)$.}	 \label{eq:icmin}\\ [0.5ex]
$sc^{max}_{mi}$  & {Maximum amount of a raw material $m \in M_{rm}$}  that can be supplied by the supplier $i \in I_p$ in each time period.\label{eq:gen_supMax}	\\ [0.5ex]
$sc^{min}_{mi}$ & {Minimum order of a raw material $m \in M_{rm}$ from a supplier $i \in I_p$ in each time period.} \label{eq:gen_scmin} 	\\ [0.5ex]
$dist_{ij}$ & {Distance between two entities $i {\in I}$ and $j {\in I}$ ($km$)}	\\ [0.5ex]
%
 \hline
\end{tabular}
\end{table}

\begin{table}[htb]
\centering
\caption{{Data parameters for the objective functions of the SSC formulation introduced in \citep{Mota2018}.}}
\footnotesize
\label{tab:params_objectives}
\begin{tabular}{cl}
\hline
\textbf{Parameter} &  \multicolumn{1}{c}{\textbf{Description}}	\\ [0.5ex]
\hline 
\multicolumn{2}{c}{Costs and related parameters} \\ \hline
$tec_g$ & {Cost to install a technology $g \in G$ (\euro).} \label{eq:gen_tec} \\ [0.5ex]
$opc_g$ & {The operating cost of a technology $g \in G$ (\euro).} \label{eq:gen_opc} \\ [0.5ex]
$psu_m$ & {Selling price of a final product $m \in M_{fp}$ (\euro).} \label{eq:psu} \\ [0.5ex]
$sc_m$ & {Inventory price of a final product $m \in M_{fp}$ (\euro).} \label{eq:sc}		\\ [0.5ex]
$rpc_m$ & {Cost of recovery product $m \in M_{rp}$ (\euro).}	 \label{eq:rpc}	\\ [0.5ex]
$rmc_m$ & {Cost of raw material $m \in M_{rm}$ (\euro).} \label{eq:rmc}		\\ [0.5ex]
$avc_a$ & {Average vehicle consumption for $a \in \mathcal{A}_{truck}$ ($l$ per $100km$).}		\\ [0.5ex] \hline
%
\multicolumn{2}{c}{Social parameters} \\ \hline
$sqmc_i$  & {Construction costs for each entity $i \in I_f \cup I_w$ (\euro).} \label{eq:sqmc} \\ [0.5ex]
$lc_i$ & {Labor cost in an entity $i \in I$ (\euro).}	\\ [0.5ex]
$\mu^{gdp}_i$ & {Inverse of the country GDP per capita  (GDPPC) where entity $i \in I$ is located.} \\ [0.5ex]
%
{$pWork_i$} & {Working condition index in an entity $i \in I$.} 	 \\ [0.5ex]
{$wi_i$} & {Minimum number of workers in entity $i \in I_f \cup I_w$.} \label{eq:wi} \\ [0.5ex]
{$wa_a$} & {Number of workers in transportation mode $a \in A$.} \label{eq:wa} \\ [0.5ex]
{$wg_g$} & {Number of workers to operate technology $g \in G$.} \label{eq:wg} \\ [0.5ex]
{$wpsq_i$} & {Minimum number of workers per square feet in entity $i \in I_f \cup I_w$.} 	 \label{eq:wpsq}	\\ [0.5ex]
%
 \hline
\multicolumn{2}{c}{Environmental impact parameters} \\ \hline
$ei_{c}$ &  {Entity installation impact for category $c \in C$ (per $m^2$). }	 \\ [0.5ex]
{$es_{ic}$} &  {Impact from supplier $i \in I_p$ for category $c \in C$. }
	 \\ [0.5ex]
\multirow{2}{*}{${eg_{mgc}}$} & {Production impact for manufacturing $m \in M_{fp}$ using technology $g \in G$}\\
& {for category $c \in C$ (per product).}	 \\  [0.5ex]
${et_{ac}}$ & {Impact of the transportation mode $k \in K$ for the category $c \in C$ (per $kg$).}	
\\[0.25cm] \hline
\end{tabular}
\end{table}

The following section presents the multi-objective formulation of the {SSC management} problem studied in \cite{Mota2018}.

\subsection{Multi-objective formulation}

To fully understand the proposed method, we show the constraints associated with the Lagrangian relaxation we perform in our method in \eqref{eq:mota_10}-\eqref{eq:mota_39}.

Constraints \eqref{eq:mota_10} and \eqref{eq:mota_11} ensure that the amount of raw materials acquired by the factories from the selected suppliers {are} within the interval {$[sc^{min}_{mi},sc^{max}_{mi}]$}.
\begin{align} 
& \sum_{\substack{a \in A,  j \in I_f}} X_{maijt} \leq sc^{max}_{mi} Y_i, \quad i \in I_{sup}, m \in M_{rm}, t \in T \label{eq:mota_10} \\
& \sum_{\substack{a \in A, j \in I_f}} X_{maijt} \geq sc^{min}_{mi} Y_i, \quad i \in I_{sup}, m \in M_{rm}, t \in T \label{eq:mota_11}
\end{align}

Constraints \eqref{eq:mota_12} and \eqref{eq:mota_13} define{, respectively,} $ec_i^{max}$ as the maximum in- and out-flow {of products between a pair of} installed entities {$i,j \in I$.}
\begin{align} 
 \sum_{m \in M, a \in A, j \in I} X_{maijt} \leq ec^{max}_{i} Y_i, \quad i \in I, t \in T \label{eq:mota_12} \\
 \sum_{m \in M, a \in A,  j \in I} X_{majit} \leq ec^{max}_{i} Y_i, \quad i \in I, t \in T \label{eq:mota_13}
\end{align}

Constraints \eqref{eq:mota_14} and \eqref{eq:mota_15} ensure that the amount of final products $m$ stored at installed factory or warehouse $i$ {is} within the interval $[ic_{mi}^{min},ic_{mi}^{max}]$.
\begin{align} 
 S_{mit} \leq ic^{max}_{mi} Y_i, \quad m \in M_{fp}, i \in I_f \cup I_w, t \in {T} \label{eq:mota_14} \\
 S_{mit} \geq ic^{min}_{mi} Y_i, \quad m \in M_{fp}, i \in I_f \cup I_w, t \in {T} \label{eq:mota_15} 
\end{align}

Constraints \eqref{eq:mota_18} and \eqref{eq:mota_19} ensure that the installation area  of each factory or warehouse $i$ {is} within {the interval} $[ea_i^{min},ea_i^{max}]$, where $ea_i^{min}$ and $ea_i^{max}$ are non-negative scalars.
\begin{align} 
& YC_i \leq ea^{max}_i Y_{i}, \quad i \in I_f \cup I_w \label{eq:mota_18} \\
& YC_i \geq ea^{min}_i Y_{i}, \quad i \in I_f \cup I_w \label{eq:mota_19}
\end{align}

Constraints \eqref{eq:mota_20} and \eqref{eq:mota_21} guarantee  that only entities selected to be installed can receive or send items.
\begin{align} 
& \sum_{m \in M, a \in A, i \in I, t \in T} X_{maijt} \geq Y_j, \quad j \in I \label{eq:mota_20} \\ 
& \sum_{m \in M, a \in A, j \in I, t \in T} X_{maijt} \geq Y_i, \quad i \in I \label{eq:mota_21}  
\end{align}

Constraints \eqref{eq:mota_26} and \eqref{eq:mota_27} guarantee that if the number of trips to transport items from/to an entity is higher than $0$, the entity must be installed.
\begin{align} 
&  Q_{aijt} \leq BigM Y_i, \quad a \in A, i,j \in I, t \in T  \label{eq:mota_26} \\ 
&  Q_{aijt} \leq BigM Y_j, \quad a \in A, i,j \in I, t \in T  \label{eq:mota_27}  
\end{align}

Constraints \eqref{eq:mota_33} {restrict} the purchase of trucks $a \in A_{truck}$ only at installed entities $i \in I$.
\begin{align} 
& K_{ai} \leq BigM Y_i, \quad a \in A_{truck}, i \in I \label{eq:mota_33} 
\end{align}

Let $pc^{max}_g$ and $pc^{min}_g$ be, respectively, the maximum and minimum {amounts} of products that technology $g \in {G}$ can produce.
For each factory $i$ and time period $t$, constraints \eqref{eq:mota_36} and \eqref{eq:mota_34} ensure that the production levels of the final product $m$ using technology $g$, ${\forall} (m,g) \in H_{prod}$, {are} within the {interval $[pc^{min}_g,pc^{max}_g]$}. Analogously, constraints \eqref{eq:mota_37} and \eqref{eq:mota_35} ensure that the remanufacturing levels {of technology $g$ {are} within the interval $[pc^{min}_g,pc^{max}_g]$}.
\begin{align}
& P_{mgit} \geq pc^{min}_g Z_{gmi}, \quad i \in I_f, (m,g) \in H_{prod}, t \in T \label{eq:mota_36} \\
& P_{mgit} \leq pc^{max}_g Z_{gmi}, \quad i \in I_f, (m,g) \in H_{prod}, t \in T \label{eq:mota_34} \\
& R_{mgit} \geq pc^{min}_g Z_{gmi}, \quad i \in I_f, (m,g) \in H_{rem}, t \in T \label{eq:mota_37} \\
& R_{mgit} \leq pc^{max}_g Z_{gmi}, \quad i \in I_f, (m,g) \in H_{rem}, t \in T \label{eq:mota_35} 
\end{align}

Constraints \eqref{eq:mota_38} and \eqref{eq:mota_39} define that production or remanufacturing technologies {can only} be selected { in installed factories}.
\begin{align}
& \sum_{g:(m,g) \in H_{prod}} Z_{gmi} \leq Y_i, \quad m \in M_{fp}, i \in I_f \label{eq:mota_38}\\
& \sum_{g:(m,g) \in H_{rem}} Z_{gmi} \leq Y_i, \quad m \in M_{fp}, i \in I_f \label{eq:mota_39}
\end{align}

Let $n^{'}, n^{''}, b^{'}$ and $b^{''}$ be natural numbers that define the dimension of the following decision variables: $u \in \mathds{R}^{n^{'}}_{{\geq 0}}, v \in  \mathds{Z}_{{\geq0}}^{n^{''}}$, $w' \in \{0,1\}^{b'}$ and $w'' \in \{0,1\}^{b''}$. Moreover, consider ${f}_{eco}$, ${f}_{env}$, ${f}_{soc}: \mathds{R}^{n{'}}\times \mathds{Z}_{{\geq0}}^{n{''}}\times  \{0,1\}^{b'} \times \{0,1\}^{b''}  \rightarrow \mathds{R}$ as the economic, environmental and social  functions of the problem, respectively.
The economic function ($f_{eco}$) and social function ($f_{soc}$) must be maximized to achieve the best  values of  NPV and of the SB, respectively. Thereby, without loss of generality, we minimize $-f_{eco}, f_{env}$ and $-f_{soc}$ to describe the multi-objective minimization problem \eqref{eq:mo_fo}-\eqref{eq:mo_ctr3}.

\begin{align}
   {\min} \quad & -f_{eco}(u,v,w',w''), f_{env}(u,v,w',w''), -f_{soc}(u,v,w',w'') \label{eq:mo_fo}\\
 s.t. \quad & \mathcal{A}^{'} u +  \mathcal{A}^{''} v \leq \beta' w' + \beta'' w''  \label{eq:mo_ctr}\\ 
   \quad & E^{'}u + E^{''} v \leq d \label{eq:mo_ctr2}\\
   \quad & u \in \mathds{R}_{\geq 0}, v \in \mathds{Z}_{\geq 0}^{n{''}}, w' \in \{0,1\}^{b'}, w'' \in \{0,1\}^{b''}  \label{eq:mo_ctr3}
   \end{align}

\noindent where $p$ and $q$ are natural numbers, 
$ \mathcal{A}^{'} \in \mathds{R}^{p \times n^{'}}$ and $E^{'} \in \mathds{R}^{q\times n^{'}}$ are  parameters associated with the real variables;
$\mathcal{A}^{''} \in \mathds{R}^{p \times n^{''}}$ and $E^{''} \in \mathds{R}^{q\times n^{''}}$ are parameters associated with the integer variables;
$\beta' \in \mathds{R}^{p \times b'}$ and $\beta'' \in \mathds{R}^{p \times b''}$ are parameters associated with the binary  variables $Y$ and $Z$ {of the problem, respectively}; $d \in \mathds{R}^q$ are the  parameters not linked to  decision variables;
$ \mathcal{A}^{'} u +  \mathcal{A}^{''} v \leq \beta' w' + \beta'' w''$ are the constraints associated with the binary variables and presented in \eqref{eq:mota_10}-\eqref{eq:mota_39};
and $E^{'}u + E^{''} v \leq d$ are the remaining constraints.

Constraints $E^{'}u + E^{''} v \leq d$ are related to the tactical planning.
They model 
(i) the material balance at entities to ensure that their in-flow  is equal to the out-flow, considering item manufacturing, remanufacturing and stock; 
(ii) cross-docking at the airports and seaports;
(iii) \typo{the} maximum flow between entities;
(iv) product delivery to meet customer demands; 
(v) the return of used products after $LF$ time periods; 
and
(vi) transportation-related decisions such as \typo{the} maximum contracted capacity of airports and seaports and maximum investment in trucks.

To calculate $f_{eco}$, i.e., the NPV, \citet{Mota2018} considered the difference between the profits from selling products and the costs involved in the SSC.  The authors also take into account the depreciation of capital invested in the SSC.
The SSC costs are those related to raw material acquisition, product manufacturing, recovery product acquisition, remanufacturing activities, transportation, handling costs at hub airport and seaport terminals, costs associated with contracts with airline and seaport companies, inventory costs and labor costs.

Function $f_{env}$ is calculated according to ReCiPe 2008 as the normalized sum of the environmental impacts due to production manufacturing and remanufacturing, transportation and entity installation. In this paper, we also consider impacts related to raw material acquisition from suppliers, which are given by 
$\sum_{c \in CC} \eta_c ( \sum_{t \in T, m \in M_{rm}, a\in A, i\in I_p,j \in I_f} es_{ic} X_{maijt})$ 
where $CC$ is the set of indicators in ReCiPe 2008 and $\eta_c$ is a normalization factor for indicator $c \in CC$.

The SB  proposed in \cite{Mota2018} as the social objective function rewards the creation of jobs due to entity installation, manufacturing and remanufacturing technology selection and transportation decisions entities in countries with lower GDP. Such an indicator composes the social objective function $f_{soc}$ of the model studied in this paper to be maximized, that also considers the working conditions in the entities, as presented in Equation \eqref{eq:fo_soc}.

\begin{equation}
\begin{split}
f_{soc} = \sum_{i \in I_f \cup I_w} \mu_i^{gdp} pWork_i (wi_i Y_i + wpsq_i YC_i) + \\ 
 \sum_{i \in I_f} \mu_i^{gdp} pWork_i\sum_{ (m,g) \in H} wg_g Z_{gmi} + \sum_{i \in I} \mu_i^{gdp} pWork_i\sum_{a \in A_{truck}} wa_a  K_{ai} + \\
\sum_{ i,j \in I} \mu_i^{gdp} pWork_i dist_{ij} \sum_{m \in M} pw_m  \sum_{ a \in A_{plane} \cup A_{ship} }\frac{wa_a}{yth} \sum_{t \in T}X_{maijt}
\end{split}
\label{eq:fo_soc}
\end{equation}

\noindent where $yth$ is the number of years in the planning horizon. For more details about the formulation, we refer to \cite{Mota2018}.

The following section discusses the solution methods proposed in this paper.

\section{Proposed method}\label{sec:method}

This section describes the proposed Lagrangian matheuristic -- called \emph{AugMathLagr} ---  to approach  the multi-objective SSCM problem discussed in the previous section. 

\emph{AugMathLagr} heuristically solves the multi-objective SSCM problem by following the same core strategy {as the} Augmented $\epsilon$-Constraint Method (AUGMECON2) \citep{Mavrotas2013}.
\emph{AugMathLagr} introduces the Lagrangian-based heuristic for solving mono-objective problems as an innovation in relation to AUGMECON2.

\citet{Mavrotas2013} proposed AUGMECON2 as an improvement of the $\epsilon$-Constraint method. {Its aim is to identify the Pareto set of a multi-objective problem by systematically solving a sequence of mono-objective problems ($\epsilon$-Constrained problems)}. In line with this, the multi-objective problem \eqref{eq:mo_fo}-\eqref{eq:mo_ctr3} can be  approached by AUGMECON2 {through} the solution of the mono-objective problem  \eqref{eq:monoOrig_fo}-\eqref{eq:monoOrig_ep3}, referred here to as MOP.

 \begin{align}
 \mbox{(MOP):} \ \min \quad & {-f_{eco}(u,v,w',w'')} - eps (  \frac{l_{env}}{r_{env}} + 0.1\frac{l_{soc}}{r_{soc}} ) \label{eq:monoOrig_fo}\\  
s.t.  \quad &  \mathcal{A}^{'} u +  \mathcal{A}^{''} v \leq \beta' w' + \beta'' w'' \label{eq:monoOrig_ctr0}\\ 
  \quad & E^{'}u + E^{''} v \leq d  \\
  & f_{env}(u,v,w',w'') + l_{env} = \epsilon_{env} \label{eq:monoOrig_ep1}\\
  & {-f_{soc}(u,v,w',w'')} + l_{soc} = \epsilon_{soc} \label{eq:monoOrig_ep2}\\
  \quad & u \in \mathds{R}_{\geq 0}, v \in \mathds{Z}_{{\geq 0}}^{n{''}}, w' \in \{0,1\}^{b'}, w'' \in \{0,1\}^{b''},   \label{eq:monoOrig_ctr3} \\
  \quad &  l_{env}, l_{soc} \in \mathds{R}_{\geq 0}\label{eq:monoOrig_ep3}
\end{align}

\noindent where $l_{env}$ and $l_{soc}$ are the slack variables of the  $\epsilon$-Constraints \eqref{eq:monoOrig_ep1} and \eqref{eq:monoOrig_ep2};   $\epsilon_{env}$ and $\epsilon_{soc}$ are scalar values defined as thresholds of $f_{env}$ and {$-f_{soc}$}, respectively; $r_{env}$ and $r_{soc}$ are positive scalars which are the absolute value of the difference between the best and worst possible values of the functions $f_{env}$ and {$-f_{soc}$}, respectively; and  $eps \in {\mathds{R}_+}$ is a small value to promote alternative optimal solutions for {$-f_{eco}$} with the best possible values of $f_{env}$ and {$-f_{soc}$}.  The coefficients $1$ and $0.1$ on the objective function mean that  $f_{env}$ must be prioritized over {$-f_{soc}$}. Their values were defined according to \citet{Mavrotas2013}.

In this paper, we estimate the worst possible value for  functions  $f_{env}$ and $-f_{soc}$ as the 
worst values they assume on the solutions that optimize each objective function.
To define the values of $\epsilon_{env}$ and $\epsilon_{soc}$, AUGMECON2 creates a grid of evenly distributed points in the Cartesian plane limited by the best and worst possible values for $f_{env}$ and {$-f_{soc}$}. The Pareto frontier approximation is composed of the solutions of the MOP considering pre-defined values of $\epsilon_{env}$ and $\epsilon_{soc}$.

Let us consider AUGMECON2 to solve the MOP. {As AUGMECON2 is an iterative strategy, an upper index is used on the functions and variables of the MOP to indicate the solutions of a given iteration.} 
Moreover, to assign values to $\epsilon_{j},  j \in \{env,soc\}$, considering an imposed number of $dg$ points equally distributed on the grid, the step values of $\epsilon_{j}$  are defined and referred to as $step_{j} = \frac{r_j}{dg}, j \in \{env,soc\}$.

In the first iteration, iteration $0$, an initial value for $\epsilon_{soc}$ must be considered, here denoted by $\epsilon_{soc}^{(0)}$, which can be, for example, the nadir point of {$-f_{soc}$}. According to constraints \eqref{eq:monoOrig_fo}, ${-f_{soc}(u,v,w',w'')} + l^{(0)}_{soc}= \epsilon_{soc}^{(0)}$, hence, $l_{soc}^{(0)} = \epsilon^{(0)}_{soc} - {(- f_{soc}(u,v,w',w''))}$. Then, to update $\epsilon_{soc}$, i.e., to define $\epsilon^{(1)}_{soc}$, $\epsilon^{(0)}_{soc}$ is decremented by an scalar $step_{soc}$ and in the next iteration, the values are updated according to: ${-f_{soc}(u,v,w',w'')} + l^{(1)}_{soc} = \epsilon^{(0)}_{soc} - step_{soc}$. First let $l_{soc}^{(0)} - step_{soc}\geq 0$ or, equivalently, $l_{soc}^{(0)} \geq step_{soc}$. {Thereby}, when $l_{soc}^{(1)}=l_{soc}^{(0)} - step_{soc}$, constraint \eqref{eq:monoOrig_ep2} {ensures} ${-f_{soc}({u,v,w',w''})} + l_{soc}^{(0)} - step_{soc} = \epsilon_{soc}^{(0)} - step_{soc}$, i.e., ${-f_{soc}(u,v,w',w'')} + l_{soc} = \epsilon_{soc}$, that is the same problem solved in the current iteration. Therefore, to avoid solving redundant problems, it is {necessary} to choose a step size whose $l_{soc} < step_{soc}$ holds. AUGMECON2 does that by selecting the step size where $\lfloor{ \frac{l_{soc}}{ step_{soc} }} \rfloor$.

\citet{Mavrotas2009} observed that if a problem is infeasible for a given value of $\epsilon_{soc}$, for smaller values of $\epsilon_{soc}$, it will {also be} infeasible. Therefore, in this case, AUGMECON2 halts {at} decrementing $\epsilon_{soc}$ to {avoid solving} unnecessary problems. This mechanism enables AUGMECON2 to investigate a lower number of problems than an enumerative $\epsilon$-Constraint method.

In the next section, we thoroughly explain the \emph{MathLagr} matheuristic, which is the Lagrangian heuristic to solve mono-objective problems in \emph{AugMathLagr}.

 \subsection{MathLagr}

 \emph{MathLagr} finds heuristic solutions through the Lagrangian relaxation of each MOP. 
 The constraints associated with binary variables $\mathcal{A}^{'} u +  \mathcal{A}^{''} v \leq \beta' w' + \beta'' w''$ can be decomposed into equivalent constraints $\mathcal{A}^{'} u +  \mathcal{A}^{''} v \leq \beta' + \beta'' w''$ and
 $\mathcal{A}^{'} u +  \mathcal{A}^{''} v \leq \mathcal{M} w + \beta'' w''$, $\mathcal{M} \in \mathds{R}^{p\times b'}$  being a matrix whose elements are $Big M$ values.
 The method keeps constraints $\mathcal{A}^{'} u +  \mathcal{A}^{''} v \leq \beta' + \beta''w''$ in the problem and relaxes constraints $\mathcal{A}^{'} u +  \mathcal{A}^{''} v \leq \mathcal{M} w' + \beta'' w''$.
 
As the elements of $\beta'$ are a trivial estimation of the  $Big M$ values of $\mathcal{M}$, we shall refer to the relaxed constraints as $\mathcal{A}^{'} u +  \mathcal{A}^{''} v \leq \beta' w' + \beta'' w''$.

 Let $\lambda \in \mathds{R}^p, \lambda \geq 0$, be the $p$-dimensional vector of Lagrange multipliers associated with the relaxed constraints.
 Each Lagrange multiplier, $\lambda_i$, penalizes in the objective function the corresponding 
 violation of constraints $a^{'}_{i,} u + a^{''}_{i,}v \leq \beta'_{i,} w' + \beta''_{i,} w''$, where $a^{'}_{i,}$, $a^{''}_{i,}$, $\beta'_{i,}$ and $\beta''_{i,}$ indicate the $i$-th row of, respectively, $\mathcal{A}^{'}$, $\mathcal{A}^{''}$, $\beta'$ and $\beta''$.

\emph{MathLagr} solves the Lagrangian mono-objective problems PRL defined by Equations \eqref{eq:monoLagr_fo}-\eqref{eq:monoLagr_ep3}. In Equation~\eqref{eq:monoLagr_fo}, $L$ is the Lagrangian function of the problem.

 \begin{align}
 {\mbox{(PRL)}}: \ \min \quad L(u,v,w',w'', \lambda) = & {-f_{eco}(u,v,w',w'')} - eps ( \frac{l_{env}}{r_{env}} + 0.1\frac{l_{soc}}{r_{soc}} )  \notag \\
 + & \lambda^T( \mathcal{A}^{'} u +  \mathcal{A}^{''} v - \beta' w' - \beta'' w''  ) \label{eq:monoLagr_fo}
  \end{align}

 \begin{align}
 s.t. \quad & \mathcal{A}^{'} u +  \mathcal{A}^{''} v \leq \beta' + \beta'' w'' \label{eq:monoLagr_ctr}\\
   \quad & E^{'}u + E^{''} v \leq d \label{eq:monoLagr_ctr2}\\
  & f_{env}(u,v,w',w'') + l_{env} = \epsilon_{env} \label{eq:monoLagr_ep1}\\
  & {-f_{soc}(u,v,w',w'')} + l_{soc} = \epsilon_{soc} \label{eq:monoLagr_ep2}\\
  \quad & u {\in \mathds{R}_{{\geq0}}}, v \in \mathds{Z}_{{\geq 0}}^{n{''}}, w' \in \{0,1\}^{b'}, w'' \in \{0,1\}^{b''}
  \label{eq:monoLagr_ctr3}\\
  \quad &  l_{env}, l_{soc} {\in \mathds{R}_{\geq 0}} \label{eq:monoLagr_ep3}
 \end{align}

The resulting problem is nonlinear. {In order to} approach it {numerically using} mixed\typo{-}integer program methods, \textit{MathLagr} employs the subgradient method \citep{Held1974}. This iterative method assigns, in its first iteration, initial values to the Lagrange multipliers, {which are usually} null values. {Then, at each iteration it relies on the resulting best lower bounds found along the iterations and on estimated upper bounds, possibly from heuristics incorporated in the method, to update the Lagrange multipliers in an attempt to refine them.} In this paper, $\lambda^{(k)}$ refers to the values of the Lagrange multipliers at iteration $k$. Moreover, for ease of notation, we shall refer to the solutions, that is, the values of the decision variables $u,v,w'$ and $w''$, as $x_{RL}^{(k)}$, $x_f^{(k)}$ or $x^{best}$.

We chose the subgradient method because it theoretically converges to the optimal solution. Even though such convergence is not guaranteed in \emph{AugMathLagr} due to the linear relaxations used, the subgradient method allows the developed matheuristic to achieve \typo{high-quality} solutions.

A solution {to} the Lagrangian problem is not necessarily feasible for the original problem. Therefore we apply a feasibility heuristic to the solution of the Lagrangian problem obtained at each iteration of the subgradient method. Since the feasibility heuristic guarantees that the variables are integer, we solve Lagrangian problems with the integer variables relaxed to speed CPLEX up.
Algorithm \ref{alg:subgradient} presents a pseudocode of the proposed Lagrangian heuristic named \emph{Mathlagr} to solve PRL.

 \begin{algorithm}[htb]
 {
 
 \caption{\emph{Mathlagr}.}
 \label{alg:subgradient}
 \SetKwInOut{Input}{Input}
 \SetKwInOut{Output}{Output}
  \Input{Maximum number of iterations, $kMax$, step size, $st$, an initial value for the Lagrange multipliers, $\lambda^{(0)}$}
 \Output{Best feasible solution $x^{best}$}
 
 {$UB:= BigM$}
 
 {$LB:= -BigM$}

 \For{$k=0$ to $k=kMax$ and $UB \neq LB$}{
  $P_{RL}^{(k)}:=$ Relaxed problem \eqref{eq:monoLagr_fo}-\eqref{eq:monoLagr_ep3} considering $\lambda^{(k)}$
  
  $x_{RL}^{(k)} := $ Solve the linear relaxation of $P_{RL}^{(k)}$

  $P_{F}^{(k)} :=$ Modify problem \eqref{eq:monoOrig_fo}-\eqref{eq:monoOrig_ep3} by assigning  0 to real variables $X_{maijt}, m \in M$, $a \in A$, $i,j \in I$, $t \in T$, with value 0 at $x_{RL}^{(k)}$

  $x_{f}^{(k)} := $ Solve $P_{F}^{(k)}$

	$UB:= \min\{{-f_{eco}}(x_{f}^{(k)}),UB\}$ 
	
	\lIf{${-f_{eco}}(x_{f}^{(k)}) < UB$}{ $x^{best} := x_{f}^{(k)}$}
	
	 $LB:= \max\{L(x_{RL}^{(k)}, \lambda^{(k)}),LB\}$ 
  
  {$g(x_f^{(k)}) := \mathcal{A}^{'} u +  \mathcal{A}^{''} v - \beta' w' - \beta'' w''$,  where $u,v,w'$ and $w''$ are the values of the
  corresponding variables in {$x_f^{(k)}$}} 
  
  $\lambda^{k+1} = st \frac{UB - L(x_{RL}^{(k)}, \lambda^{(k)}) }{|| g(x_{f}^{(k)})||}$
  }
  \Return{$x^{best}$}
 }
 \end{algorithm}

Algorithm \ref{alg:subgradient} {uses the following as input: (i) the} maximum number of iterations {of the method}, $kMax${; (ii) } the step size to adjust the Lagrange multipliers at each iteration, $st${; and (iii)} an initial value for the Lagrange multipliers, $\lambda^{(0)}$. In lines 1 and 2, the upper (UB) and lower bounds (LB) for the objectives are initialized.
Then, CPLEX \citep{CPLEX} is employed to solve the Lagrangian problem {PRL} with the integer variables relaxed at iteration $0$, $P_{RL}^{(0)}$. 
In line 6, {problem $P_{F}^{(k)}$ is created. It is} the MOP with  the real variables $X_{maijt}, m \in M, a \in A, i,j \in I, t \in T$, which are null in the relaxed solution of the iteration $k$, i.e $x_{RL}^{(k)}$, fixed at zero.
 In line {7}, CPLEX returns a solution for $P_{F}^{(k)}$ and a feasible solution to the MOP, i.e., $x_f^k$.
 In line {8}, the algorithm updates the upper bound ($UB$) with the best economic objective function value of a solution encountered up to that iteration. In line {9}, $x^{best}$ is updated if the value of its economic objective function is better than $UB$. In the first iteration, $UB$ is exactly the value of ${-f_{eco}}(x_{f}^{(0)})$, since it is the first upper bound obtained by the method. 
The lower bound ($LB$), on the other hand, is the Lagrangian function value of the relaxed solution obtained in that iteration, that is, $L(x_{RL}^{(0)}, \lambda^{(0)})$, defined in line {10}. In the next iterations, $LB$ can be worse than those found in former iterations. Therefore, $LB$ is the {largest} lower bound {up to that iteration. } 
In line 11, $g(x_f^{(k)})$ is defined as the violation of Constraints~\eqref{eq:monoLagr_ctr} according to solution $x_f^{(k)}$. 
 Then, the Lagrange multipliers are updated as indicated in {12}. 
The process {is repeated} until either the upper bound is equal to the lower bound or the maximum number of iterations {has been} reached. 
Algorithm \ref{alg:subgradient} returns the best feasible solution found over the iterations referred to as $x^{best}$ in line {14}.

 \subsection{Multi-objective matheuristic}

Algorithm \ref{alg:aug2} presents a pseudocode of the \emph{AugMathLagr} solution method that creates a sequence of mono-objective problems using the multi-objective method AUGMECON2.  This algorithm has as input the maximum number of iterations of the Lagrangian heuristic, $kMax$, the step size of the subgradient method, $st$, the initial value of the Lagrange multipliers, $\lambda^{(0)}$, and the number of points {of} the grid, $dg$.

 \small
\begin{algorithm}[htb]
{
 \caption{\emph{AugMathLagr}.}
 \label{alg:aug2}
 \SetKwInOut{Input}{Input}
 \SetKwInOut{Output}{Output}
 \Input{maximum number of iterations, $kMax$, the step size, $st$, an initial value for the Lagrange multipliers, $\lambda^{(0)}$, the number of grid points, $dg$}
 \Output{Pareto frontier approximation $\mathcal{P}$}
 Estimate the lower bounds $fL_{i}$ and the upper bounds $fU_{i}$, $i \in \{eco, env, soc\}$ for the objective functions $-f_{eco}, f_{env}$ and $-f_{soc}$

 \small Identify ranges {$r_{j}=fU_{j}-fL_{j}, j \in \{env,soc\}$} of the environmental and social objective functions, respectively

 $step_j = \frac{r_j}{dg}, j \in \{env,soc\}$	
 
 $\mathcal{P} := \emptyset$
 
 ${gr}_{env}=0$
 
 \For{$\epsilon_{env}=fU_{env}$ until ${gr}_{env} < dg$}{
 ${gr}_{soc}=0$
 
  \For{$\epsilon_{soc}=fU_{soc}$ until ${gr}_{soc} < {dg}$}{
	 $MOP:=$ Problem \eqref{eq:monoOrig_fo}-\eqref{eq:monoOrig_ep3}
	 
	 $x_{f}:=$ solve ${MOP}$ by the Algorithm \ref{alg:subgradient} considering $\epsilon_{env}$ and $\epsilon_{soc}$ 	
	 
	  \uIf{$x_{f}$ is feasible}{
	  
		   $\mathcal{P} := \mathcal{P} \cup {\{x_{f}\}}$

        $l_{soc},l_{env}:=$ values for the slack variables in $x_f$
		 
		 $gr_{soc} = gr_{soc} + 1 +\lfloor \frac{l_{soc}}{step_{soc}} \rfloor $
		 
	   $\epsilon_{soc} = \epsilon_{soc} - step_{soc} \left ( 1 +\lfloor \frac{l_{soc}}{step_{soc}} \rfloor \right )$
	  }
	  \lElse{
	    $gr_{soc} = dg$
	  }
	}
   $gr_{env} = gr_{env} + 1$

$\epsilon_{env} = \epsilon_{env} - step_{env}$  
    }
   
   $\mathcal{P}$ := Remove dominated solutions from $\mathcal{P}$
   
   \Return{$\mathcal{P}$}
}
\end{algorithm}

\normalsize

\citet{Mavrotas2009} suggests the use of lexicographic optimization for estimating the lower and upper bounds of each objective function. However, this approach is computationally expensive. Therefore, due to the complexity {of the} studied model, to define these bounds in line 1 of Algorithm \ref{alg:aug2}{,} each objective function {is individually} optimized by applying the Lagrangian heuristic described in Algorithm \ref{alg:subgradient} to the problems $\min \mbox{ } f_{i}(x) + \lambda^T( {g(x)}), s.t.:$ Constraints
$\eqref{eq:monoLagr_ctr}$,$\eqref{eq:monoLagr_ctr2}$ and $\eqref{eq:monoLagr_ctr3}$, for $f_i \in \{-f_{eco}, f_{env}, -f_{soc}$ \}.
The lower bounds $fL_{i}$ and the upper bounds $fU_{i}$ are the lowest and highest values of each objective function $-f_{eco}$, $f_{env}$ and $-f_{soc}$ and are respectively applied to a relaxed solution and an incumbent solution found along the iterations.

In line 2, the algorithm estimates the ranges of the objective functions by calculating the difference between the lower and upper bounds estimated in line 1.
The number of grid points required for each objective function is pre-defined by {$dg$}. The initial values of $\epsilon_{env}$ and $\epsilon_{soc}$ are set as the estimated upper bound values for the respective objective functions. In line 3, the algorithm calculates the step values for $\epsilon_{env}$ and $\epsilon_{soc}$ according to the number of grid points and the ranges of the corresponding objective functions. In line 4, the algorithm initializes the Pareto frontier approximation $\mathcal{P}$ as an empty set and in line 5 a control variable $gr_{env}$ indicating the point in the grid is initialized as zero.

In line 10, \emph{AugMathLagr} employs Algorithm \ref{alg:subgradient} to heuristically solve the Lagrangian mono-objective problem MOP, constructed in line 9.
If solution $x_{f}$ is feasible, in line 12 the approximation of the Pareto frontier, $\mathcal{P}${,} is updated with $x_{f}$.
In lines 13 and 14, the algorithm avoids solving redundant problems by checking the ratio between the slack variable of the social constraint and the {step} size.
Otherwise, if $x_{f}$ is {infeasible}, it means that further decreasing $\epsilon_{soc}$ will only result in {infeasible} problems. Therefore, the algorithm sets ${gr}_{soc}={dg}$ in line 15, preventing the method {from} solving these problems {and making it proceed} to the next value of $\epsilon_{env}$.

In line 20, the dominated solutions are deleted from $\mathcal{P}$ and the updated Pareto frontier approximation $\mathcal{P}$ is returned.

In addition to introducing \emph{AugMathLagr}, we have also adapted the multi-objective matheuristic \emph{AugMathFix}  \citep{Tautenhain2019} to find solutions to the target SSC {management} problem. \typo{The} next section briefly explains this adaptation.

\subsection{Adaptation of {AugMathFix}}
To better evaluate the proposed matheuristic, {\emph{AugMathLagr}}, we compare its performance with an \linebreak AUGMECON2-based matheuristic recently proposed to approach {an} SSC problem, the \emph{AugMathFix}. As the SSC problem is different from the one studied in this paper, we had to adapt such matheuristic, the \emph{AugMathLagr}, to the target SSC management problem. 
The adapted \emph{AugMathFix} selects the same constraints {to be relaxed as} \emph{AugMathLagr}. To find feasible solutions to relaxed problems, in addition to imposing which warehouses and entities will not be opened, \emph{AugMathFix} employs the same strategy as \emph{AugMathLagr} to fix decision variables.

Details about the implementation of \emph{AugMathFix} to the supply chain studied in this paper are presented in Section 2 of the Supplementary Material. The mono-objective heuristic in  \emph{AugMathFix}  is called  \emph{MathFix} and \typo{was} also adapted in this paper.

\typo{The} next section presents the instance generator introduced in this paper.

\section{Instance generator}\label{sec:instance_generator}

The optimization models found in the literature for sustainable supply chains {are usually} problem-specific \citep{Povoa2018} or based on randomly generated parameters. As a consequence, defining a generic instance generator is not an easy task since it should incorporate several traits from this diverse range of characteristics. In this context, this paper introduces a methodology to generate random instances {to} target {the SSC} formulation proposed in \citep{Mota2018}.

In order to define the ranges of the parameters, the data distribution, and so forth, we relied on real data, in particular, from the case study in \citep{Mota2018}. This case study is discussed in Section \ref{sec:experiments_cs}.
Section 3 of the Supplementary Material shows the guidelines for estimating the data parameters of the uniform random variables.

In the next section, we illustrate a small case study created by the instance generator introduced in this paper.

\subsection{Illustrative example}

This section presents a small example to illustrate how the instance generator works when considering a planning horizon with only two periods. In this example, the SSC has one supplier, one factory, two warehouses, two customers, two airports and two seaports. Moreover, it produces only one final product using two types of raw materials. {The example also considers} three  production and remanufacturing {technologies.} In our example, we refer to the production technologies {as} $g_0, g_1$ and $g_2$ and to the remanufacturing technologies {as} $g_3$, $g_4$ and $g_5$.
 The transportation modes are trucks, airplanes and boats{; and} there are two types of trucks to perform land transportation, each of them specified by $k_0$ and $k_1$.

Table \ref{tab:example_payoff} shows the results (time to solutions and upper bounds) achieved by the Lagrangian matheuristic {whose pseudocode is in Algorithm \ref{alg:subgradient}} - \emph{MathLagr}  when optimizing each objective function {of the MOP} individually.
CPLEX \citep{CPLEX} was the tool {that} solved the Lagrangian relaxed problems with a stopping criterion of $1\%$ {of optimality gap.}
The first column of Table \ref{tab:example_payoff} identifies the optimized function. 
The third, fourth and fifth columns report, respectively, the values of the economic, environmental and social functions of the heuristic solution {to} the mono-objective problem whose {objective function} is that indicated in the first column.

 \begin{table}[htb]
 \centering
 \caption{Results for optimizing each objective function individually. Very small and very large values are expressed in E notation.}
 \label{tab:example_payoff}
 \small
 \begin{tabular}{|c|r|r|r|r|}
 \hline
 {\textbf{Function to }} & \multirow{2}{*}{\textbf{Time(s)}} &  \multicolumn{3}{|c|}{\textbf{{Values of the objective functions}}}\\
 \cline{3-5}
 \multicolumn{1}{|c|}{{\textbf{be optimized}}} 		 & 	 \multicolumn{1}{|c|}{}	 & 	 \multicolumn{1}{|c|}{${{f_{eco}}}$} 	 
 & 	 \multicolumn{1}{|c|}{${{f_{env}}}$} 	
 & 	 \multicolumn{1}{|c|}{${{f_{soc}}}$}  \\ \hline
 ${\max {f_{eco}}}$  		 & 	0.306 	 & 	\minor{1.096e+06} 	 & 	\minor{6.185e+05} 	 & 	11.704 \\ \hline
${\min {f_{env}}}$		 & 	0.256 	 & 	\minor{-1.094e-07} 	 & 	\minor{3.885e+05} 	 & 	50.697 \\ \hline
${\max {f_{soc}}}$ 		 & 	0.091	 & 	\minor{-2.503e-08}	 & 	\minor{4.593e+05} 	 & 	499.081 \\ \hline
\end{tabular}
\end{table}

\begin{figure}[htb] 
  \begin{minipage}{0.45\textwidth}  
   \centering
\tikzset{stateDots/.style ={ top color =white , bottom color = white,draw,white, text=black , minimum width =0.2 cm}}
\tikzset{stateLabel/.style ={ top color =white , bottom color = white,draw,white, text=black , minimum width =0.2 cm}}
\tikzset{stateInv/.style ={ top color =white , bottom color = white,draw,white, text=black , minimum width =0.2 cm}}
\begin {tikzpicture}[-latex ,auto ,node distance =3 cm and 3cm ,on grid ,
semithick ,
state/.style ={ rectangle,rounded corners ,top color =white , bottom color = gray!20 ,
draw,gray , text=black , minimum width = 0.2 cm,
}]
\scriptsize
\node[state] (P) {$S$}; 
\node[state] (F) [right=2.5cm of P]{$F$}; 
\node[state] (W1) [right=2.2cm of F]{$W_1$}; 
\node[state] (C1) [right=2.2cm of W1]{$C_1$}; 
\node[state] (C2) [below=2cm of C1]{$C_2$}; 
\node[state] (W2) [left=2cm of C2]{$W_2$}; 
\node[state] (Air2) [left =2cm of W2]{$Air_2$}; 
\node[state] (Air1) [left =2cm of Air2]{$Air_1$}; 
\node[stateInv] (Sea2) [below =2cm of Air2]{}; 
\node[stateInv] (Sea1) [left =2cm of Sea2]{}; 
\scriptsize
\node[stateLabel] (Fprod) [above = 0.6cm of F]{\tiny $P_{2010} = 60588$}; 
\node[stateLabel] (W1S) [above= 0.6cm of W1]{\tiny $S_{220} = 334$};
\node[stateLabel] (W2S) [above= 0.5cm of W2]{\tiny $S_{230} = 15969$};
%
\path (P) edge [left =25] node[above =0 cm] {\tiny $X{00010} = \minor{11236}$} (F); 
\path (P) edge [bend right =20] node[below =0.02 cm] {\tiny $X{10010} =\minor{2012}$} (F); 
\path (W1) edge [right =25] node[right =0 cm] {\tiny $X{20260} = 46189$} (Air1);
\path (W1) edge [left = 20] node[above =0.12 cm] {\tiny $X{20240} = 14065$} (C1);
\path (W2) edge [right =25] node[below =0.2 cm] {\tiny $X{20350} = 30220$} (C2);
\path (F) edge [right =25] node[above =0.12 cm] {\tiny $X{20120} = 60588$} (W1);
\path (Air2) edge [right =25] node[below =0.2 cm] {\tiny $X{20730} = 46189$} (W2);
\path (Air1) edge [left =20] node[below =0.2 cm] {\tiny $X{22670} = 46189$} (Air2);


\end{tikzpicture}
   {{\footnotesize \centering First period of the planning horizon.}}
  \end{minipage}
   \begin{minipage}{0.50\textwidth}  
    \centering
\tikzset{stateDots/.style ={ top color =white , bottom color = white,draw,white, text=black , minimum width =0.2 cm}}
\tikzset{stateLabel/.style ={ top color =white , bottom color = white,draw,white, text=black , minimum width =0.2 cm}}
\tikzset{stateInv/.style ={ top color =white , bottom color = white,draw,white, text=black , minimum width =0.2 cm}}
\begin {tikzpicture}[-latex ,auto ,node distance =3 cm and 3cm ,on grid ,
semithick ,
state/.style ={ rectangle,rounded corners ,top color =white , bottom color = gray!20 ,
draw,gray , text=black , minimum width = 0.2 cm,
}]
\scriptsize
\node[state] (P) {$S$}; 
\node[state] (F) [right=2.5cm of P]{$F$}; 
\node[state] (W1) [right=3.5cm of F]{$W_1$}; 
\node[state] (C1) [right=3.5cm of W1]{$C_1$}; 
\node[state] (C2) [below=2cm of C1]{$C_2$}; 
\node[state] (W2) [left=2.7cm of C2]{$W_2$}; 
\node[stateInv] (Air2) [left =1cm of W2]{}; 
\node[stateInv] (Air1) [left =1cm of Air2]{}; 
\node[state] (Sea2) [below =1.5cm of Air2]{$Sea_2$}; 
\node[state] (Sea1) [left =4.5cm of Sea2]{$Sea_1$}; 
\scriptsize
\node[stateLabel] (Fprod) [above = 0.6cm of F]{\tiny $P_{2011} = \minor{21585}$};  
\node[stateLabel] (Freman) [above= 1.0cm of F]{\tiny $R_{2511}= \minor{1657}$}; 
\node[stateLabel] (W1S) [above= 0.6cm of W1]{\tiny $S_{221} = 334$};  
\node[stateLabel] (W2S) [above= 0.6cm of W2]{\tiny $S_{231} = 384$}; 
%
\path (P) edge [left =25] node[above =0 cm] {\tiny $X{00011} = \minor{4003}$} (F); 
\path (P) edge [bend right =20] node[below =0.02 cm] {\tiny $X{10011} = \minor{717}$} (F); 
\path (F) edge [bend left =10] node[above =0.01 cm] {\tiny $X{20121} = 23242$} (W1);
\path (W1) edge [bend left =10] node[below =0.01 cm] {\tiny $X{30211} = \minor{6630}$} (F); 
\path (W1) edge [bend right =10] node[left =0 cm] {\tiny $X{20281} = 7476$} (Sea1);
\path (W2) edge [bend left =10] node[right =0.1 cm] {\tiny $X{30391} = \minor{4524}$} (Sea2); 
\path (W1) edge [bend left = 10] node[above =0.05 cm] {\tiny $X{20241} = 15766$} (C1);
\path (W2) edge [bend left =10] node[above =0.05 cm] {\tiny $X{20351} = 23061$} (C2);
\path (Sea1) edge [bend right =10] node[right =0.1 cm] {\tiny $X{30821} = \minor{4524}$} (W1); 
\path (Sea2) edge [bend left =10] node[left=0.1cm] {\tiny $X{20931} = 7476$} (W2);
\path (C1) edge [bend left = 10] node[below =0.1 cm] {\tiny $X{30421} = \minor{2106}$} (W1); 
\path (C2) edge [bend left=10] node[below =0.01 cm] {\tiny $X{30531} = \minor{4524}$} (W2); 
\path (Sea1) edge [bend left =10] node[above =0.1 cm] {\tiny $X{22671} = 7476$} (Sea2);
\path (Sea2) edge [bend left =10] node[below =0.1 cm] {\tiny $X{33981} = \minor{4524}$} (Sea1); 


\end{tikzpicture}
 \
    {{\footnotesize \centering Second period of the planning horizon.}}
    \end{minipage}
  \caption{{Illustrative example of the solution obtained by the optimization of the social criterion of {an} SSC generated by the proposed instance generator.}} 
  \label{fig:example_soc}
 \end{figure}
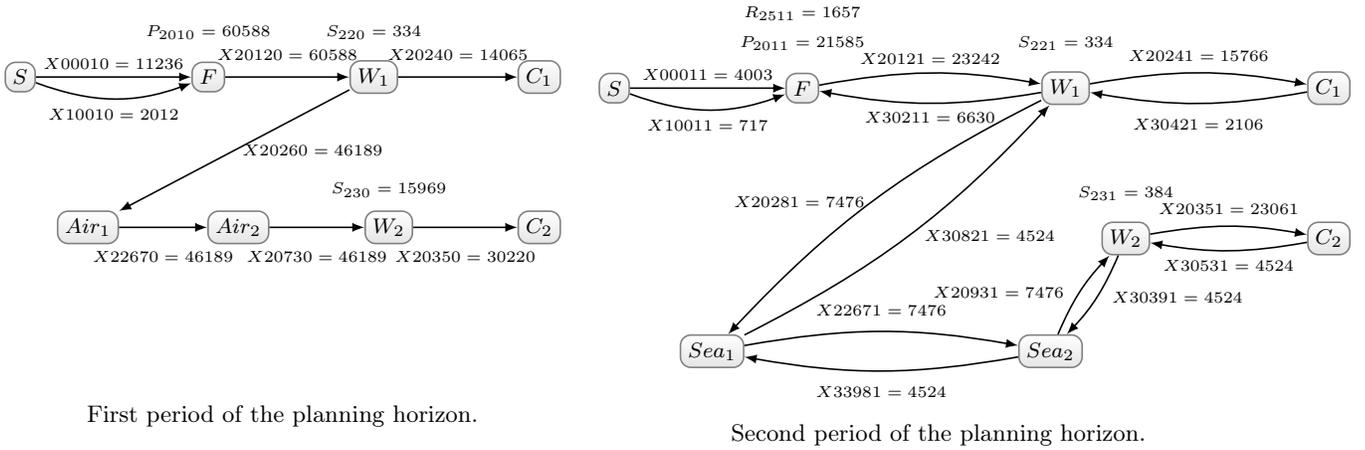

Figures \ref{fig:example_eco} to \ref{fig:example_soc} display the decision variable values regarding each time period of the planning horizon obtained by \emph{MathLagr} when optimizing the economic, environmental and social functions, respectively\footnote{\minor{The values were rounded for readability purposes.}}. These figures illustrate the flow of items {indicated by}  the values of the variables of production, remanufacturing and storage {on the arcs of the network}. Labels $S$, $F$, $W_1$, $W_2$, $C_1$, $C_2$, $Air_1$, $Air_2$, $Sea_1$ and $Sea_2$ represent, respectively, the supplier, the factory, the first warehouse, the second warehouse, the first customer, the second customer, the first airport, the second airport, the first seaport and the second seaport.

Moreover, to identify these entities in the indices of the variables, we assign integer numbers to them. Entity $S$ {corresponds} to index $0$, $F$ to $1$, and so on until $Sea_2${, which} is identified by index $9$ in the variables.
For the same reason, technologies $g_0$, $g_1$, $g_2$, $g_3$, $g_4$ and $g_5$ are assigned to sequential indices $0$ to {$5$}  and trucks $k_0$ and $k_1$ {to} indices $0$ and $1$. 
The labels on the arcs correspond to the decision variables, according to Table \ref{tab:decision_vars}. 
In this example, warehouse $W_2$, customer $C_2$, airport $Air_2$ and $Port_2$ are located in a continent different from the remaining entities. Thereby to transport items between entities, in this case, air or sea transportation is mandatory.

 \begin{figure}[H] 
 \centering
   \begin{minipage}{0.49\textwidth}  
    \centering
  { 
\tikzset{stateDots/.style ={ top color =white , bottom color = white,draw,white, text=black , minimum width =0.2 cm}}
\tikzset{stateLabel/.style ={ top color =white , bottom color = white,draw,white, text=black , minimum width =0.2 cm}}
\tikzset{stateInv/.style ={ top color =white , bottom color = white,draw,white, text=black , minimum width =0.2 cm}}
\begin {tikzpicture}[-latex ,auto ,node distance =3 cm and 3cm ,on grid ,
semithick ,
state/.style ={ rectangle,rounded corners ,top color =white , bottom color = gray!20 ,
draw,gray , text=black , minimum width = 0.2 cm,
}]
\scriptsize
\node[state] (P) {$S$}; 
\node[state] (F) [right=2.5cm of P]{$F$}; 
\node[stateInv] (W1) [right=2.5cm of F]{}; 
\node[stateInv] (W2) [below=2cm of W1]{}; 
\node[state] (C1) [right=3cm of W1]{$C_1$}; 
\node[state] (C2) [below=2cm of C1]{$C_2$}; 
\node[stateInv] (Air2) [left =0.5cm of W2]{}; 
\node[stateInv] (Air1) [left =2cm of Air2]{}; 
\node[state] (Sea1) [left =2cm of Air2]{$Sea_1$}; 
\node[state] (Sea2) [left =0.1cm of W2]{$Sea_2$}; 
\scriptsize
\node[stateLabel] (Fprod) [above = 0.6cm of F]{\tiny $P_{2110} =50827$}; 
\node[stateLabel] (Freman) [above= 1.0cm of F]{\tiny $S_{210}=\minor{6542}$}; 
\path (F) edge [bend right =25] node[left =0 cm] {\tiny $X_{20180} = 30220$} (Sea1);
\path (F) edge [ right =10] node[above =0 cm] {\tiny $X_{20140} = 14065$} (C1);
\path (Sea2) edge [right =25] node[above =0.1 cm] {\tiny $X_{20950} = 30220$} (C2);
\path (P) edge [left =25] node[above =0 cm] {\tiny $X_{01010} =  \minor{14383}$} (F); 
\path (P) edge [bend right =20] node[below =0 cm] {\tiny $X_{11010} = \minor{9386}$} (F); 
\path (Sea1) edge [right =25] node[above =0.1 cm] {\tiny $X_{23890} = 30220$} (Sea2);



\end{tikzpicture}
}
  {\footnotesize \centering First period of the planning horizon.}
  \end{minipage}
  \begin{minipage}{0.49\textwidth}  
   \centering
 {
\tikzset{stateDots/.style ={ top color =white , bottom color = white,draw,white, text=black , minimum width =0.2 cm}}
\tikzset{stateLabel/.style ={ top color =white , bottom color = white,draw,white, text=black , minimum width =0.2 cm}}
\tikzset{stateInv/.style ={ top color =white , bottom color = white,draw,white, text=black , minimum width =0.2 cm}}
\begin {tikzpicture}[-latex ,auto ,node distance =3 cm and 3cm ,on grid ,
semithick ,
state/.style ={ rectangle,rounded corners ,top color =white , bottom color = gray!20 ,
draw,gray , text=black , minimum width = 0.2 cm,
}]
\scriptsize
\node[state] (P) {$S$}; 
\node[state] (F) [right=2.5cm of P]{$F$}; 
\node[stateInv] (W1) [right=2.5cm of F]{}; 
\node[stateInv] (W2) [below=2cm of W1]{}; 
\node[state] (C1) [right=3cm of W1]{$C_1$}; 
\node[state] (C2) [below=2cm of C1]{$C_2$}; 
\node[stateInv] (Air2) [left =0.5cm of W2]{}; 
\node[stateInv] (Air1) [left =2cm of Air2]{}; 
\node[state] (Sea1) [left =2cm of Air2]{$Sea_1$}; 
\node[state] (Sea2) [left =0.1cm of W2]{$Sea_2$}; 
\scriptsize
\node[stateLabel] (Fprod) [above = 0.6cm of F]{\tiny $P_{2111} = \minor{30628}$};  
\node[stateLabel] (Freman) [above= 1.0cm of F]{\tiny $R_{2411} = \minor{1657}$}; 
\node[stateLabel] (Freman) [above= 1.4cm of F]{\tiny $ $};
%
\path (P) edge [left =25] node[above =0 cm] {\tiny $X_{01011} = \minor{8667}$} (F); 
\path (P) edge [bend right =20] node[below =0 cm] {\tiny $X_{11011} = \minor{5656}$} (F); 
\path (F) edge [bend right =25] node[left =0 cm] {\tiny $X{20181} = 23061$} (Sea1);
\path (F) edge [left =10] node[above =0 cm] {\tiny $X{20141} = 15766$} (C1);
\path (Sea1) edge [left =20] node[right =0.1 cm] {\tiny $X{30811} = \minor{4524}$} (F); 
\path (Sea2) edge [right =25] node[above =0.1 cm] {\tiny $X{20951} = 23061$} (C2);
\path (C1) edge [bend left =10] node[below =0.1 cm] {\tiny $X{30411} = \minor{2106}$} (F); 
\path (C2) edge [bend left =20] node[below =0.1 cm] {\tiny $X{30591} = \minor{4524}$} (Sea2); 
\path (Sea2) edge [bend left =30] node[below =0.1 cm] {\tiny $X{33981} = \minor{4524}$} (Sea1); 
\path (Sea1) edge [right =25] node[above =0.1 cm] {\tiny $X{23891} = 23061$} (Sea2);

\end{tikzpicture}
}
 {{\footnotesize \centering Second period of the planning horizon.}}
  \end{minipage}
\caption{{Illustrative example of the solution obtained by the optimization of the economic criterion of {an} SSC generated by the proposed instance generator.}}
\label{fig:example_eco}
 \end{figure}
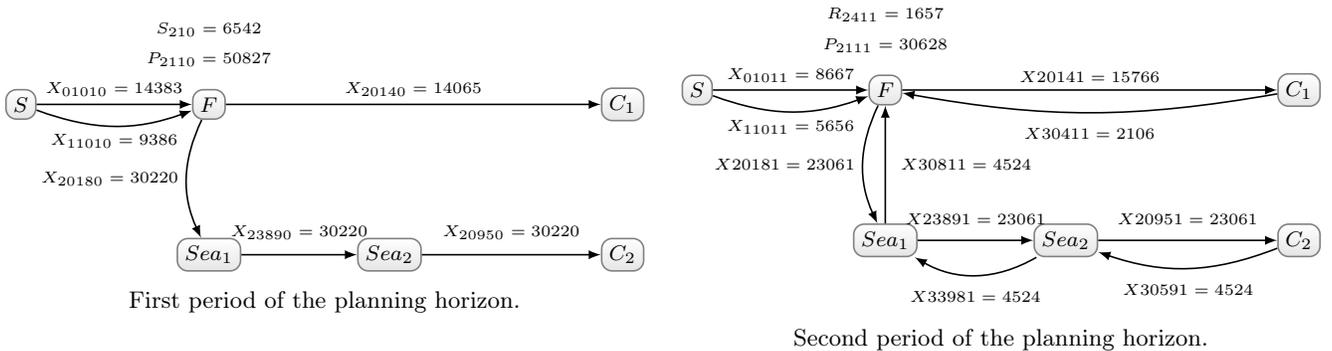

 \begin{figure}[htb] 
 \centering
   \begin{minipage}{0.49\textwidth}  
    \centering
\tikzset{stateDots/.style ={ top color =white , bottom color = white,draw,white, text=black , minimum width =0.2 cm}}
\tikzset{stateLabel/.style ={ top color =white , bottom color = white,draw,white, text=black , minimum width =0.2 cm}}
\tikzset{stateInv/.style ={ top color =white , bottom color = white,draw,white, text=black , minimum width =0.2 cm}}
\begin {tikzpicture}[-latex ,auto ,node distance =3 cm and 3cm ,on grid ,
semithick ,
state/.style ={ rectangle,rounded corners ,top color =white , bottom color = gray!20 ,
draw,gray , text=black , minimum width = 0.2 cm,
}]
\scriptsize
\node[state] (P) {$S$}; 
\node[state] (F) [right=2.5cm of P]{$F$}; 
\node[stateInv] (W1) [right=2.5cm of F]{}; 
\node[stateInv] (W2) [below=2cm of W1]{}; 
\node[state] (C1) [right=3cm of W1]{$C_1$}; 
\node[state] (C2) [below=2cm of C1]{$C_2$}; 
\node[stateInv] (Air2) [left =0.5cm of W2]{}; 
\node[stateInv] (Air1) [left =2cm of Air2]{}; 
\node[state] (Sea1) [left =2cm of Air2]{$Sea_1$}; 
\node[state] (Sea2) [left =0.1cm of W2]{$Sea_2$}; 
\scriptsize
\node[stateLabel] (Fprod) [above = 0.6cm of F]{\tiny $P_{2010} = \minor{71047}$};  
\node[stateLabel] (Freman) [above= 1.0cm of F]{\tiny $S_{210}= \minor{26762}$}; 
%
\path (P) edge [left =25] node[above =0 cm] {\tiny $X_{01010} = \minor{13175}$} (F); 
\path (P) edge [bend right =20] node[below =0 cm] {\tiny $X_{11010} = \minor{2360}$} (F); 

\path (F) edge [bend right =25] node[left =0 cm] {\tiny $X_{21180} = 30220$} (Sea1);
\path (F) edge [ right =10] node[above =0 cm] {\tiny $X_{21140} = 14065$} (C1);
\path (Sea2) edge [right =25] node[above =0.1 cm] {\tiny $X_{21950} = 30220$} (C2);
\path (Sea1) edge [right =25] node[above =0.1 cm] {\tiny $X_{23890} = 30220$} (Sea2);

\end{tikzpicture}
   {{\footnotesize \centering First period of the planning horizon.}}
  \end{minipage}
   \begin{minipage}{0.49\textwidth}  
    \centering
\tikzset{stateDots/.style ={ top color =white , bottom color = white,draw,white, text=black , minimum width =0.2 cm}}
\tikzset{stateLabel/.style ={ top color =white , bottom color = white,draw,white, text=black , minimum width =0.2 cm}}
\tikzset{stateInv/.style ={ top color =white , bottom color = white,draw,white, text=black , minimum width =0.2 cm}}
\begin {tikzpicture}[-latex ,auto ,node distance =3 cm and 3cm ,on grid ,
semithick ,
state/.style ={ rectangle,rounded corners ,top color =white , bottom color = gray!20 ,
draw,gray , text=black , minimum width = 0.2 cm,
}]
\scriptsize
\node[state] (P) {$S$}; 
\node[state] (F) [right=2.5cm of P]{$F$}; 
\node[stateInv] (W1) [right=2.5cm of F]{}; 
\node[stateInv] (W2) [below=2cm of W1]{}; 
\node[state] (C1) [right=3cm of W1]{$C_1$}; 
\node[state] (C2) [below=2cm of C1]{$C_2$}; 
\node[stateInv] (Air2) [left =0.5cm of W2]{}; 
\node[stateInv] (Air1) [left =2cm of Air2]{}; 
\node[state] (Sea1) [left =2cm of Air2]{$Sea_1$}; 
\node[state] (Sea2) [left =0.1cm of W2]{$Sea_2$}; 
\scriptsize
\node[stateLabel] (Fprod) [above = 0.6cm of F]{\tiny $P_{2011} = \minor{994}$};  
\node[stateLabel] (Freman) [above= 1.0cm of F]{\tiny $R_{2411} = \minor{11071}$}; 
\node[stateLabel] (Freman) [above= 1.4cm of F]{\tiny $ $};
%
\path (P) edge [left =25] node[above =0 cm] {\tiny $X{01011} = \minor{184}$} (F); 
\path (P) edge [bend right =20] node[below =0 cm] {\tiny $X{11011} = 33$} (F);
\path (F) edge [bend right =25] node[left =0 cm] {\tiny $X{21181} = 23061$} (Sea1);
\path (F) edge [left =10] node[above =0 cm] {\tiny $X{21141} = 15766$} (C1);
\path (Sea1) edge [left =20] node[right =0.1 cm] {\tiny $X{31811} = 30220$} (F);
\path (Sea2) edge [right =25] node[above =0.1 cm] {\tiny $X{21951} = 23061$} (C2);
\path (C1) edge [bend left =10] node[below =0.1 cm] {\tiny $X{31411} = 14065$} (F);
\path (C2) edge [bend left =20] node[below =0.1 cm] {\tiny $X{31591} = 30220$} (Sea2);
\path (Sea2) edge [bend left =30] node[below =0.1 cm] {\tiny $X{33981} = 30220$} (Sea1);
\path (Sea1) edge [right =25] node[above =0.1 cm] {\tiny $X{23891} = 23061$} (Sea2);

\end{tikzpicture}
    {{\footnotesize \centering Second period of the planning horizon.}}
    \end{minipage}
  \caption{{Illustrative example of the solution obtained by the optimization of the environmental criterion of {an} SSC generated by the proposed instance generator.}}
  \label{fig:example_env}
 \end{figure}
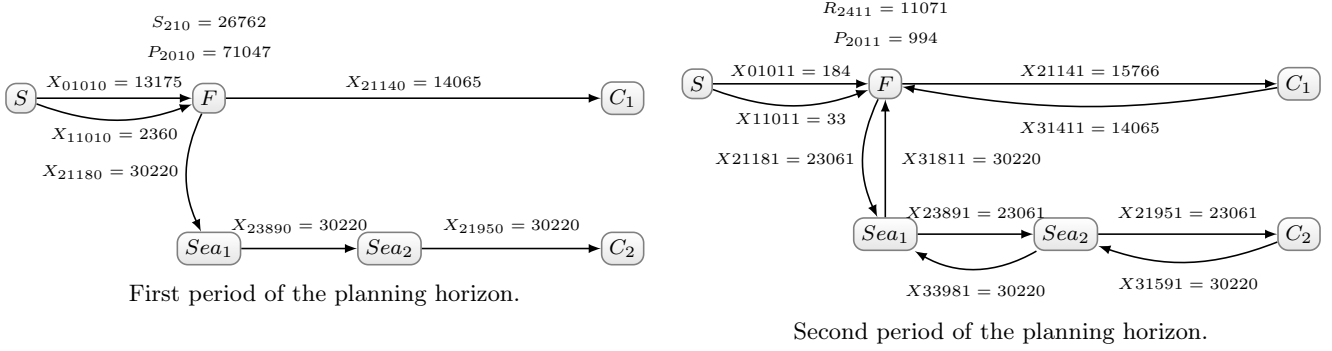

No warehouse is opened in the solutions that \typo{optimize} the economic function, i.e., maximizes the NPV, and minimizes the environmental impact using the ReCiPe 2008 methodology. The solution that minimizes the environmental function indicates that avoiding to install warehouses contributes to reducing the environmental impacts.
Since the decisions to open warehouses contribute positively to the social objective function, the solution to the maximization of this indicator indicates the opening of the two warehouses at maximum capacity.

Truck $k_1$ is selected for land transportation in all the entities in the environmental solution, whereas truck $k_0$ is preferable in the economic solution. In particular, overall, truck $k_1$ has a higher price and a lower environmental impact than truck $k_0$. 
In the environmental optimization solution, apart from land transportation, only  seaports are selected due to the lower environmental impact {they cause}.
In the economic, environmental and in the second period of the social optimization solutions, seaports are responsible for delivering products to customer $C_2$ and returning the used products to the factory.
In the first period of the solution that optimizes the social function, in contrast, airports are responsible for delivering products to customer $C_2$.
The stock in warehouse $W_2$ at the end of the first time period in the social optimization solution is used to meet the demand by customer $C_2$ \typo{in} the second time period.

The environmental and social solutions employ the manufacturing \typo{technology} $g_0$, whereas the economic solution uses the cheaper technology $g_1$.
The remanufacturing levels at factories are higher in the environmental solution, suggesting \typo{that} reverse logistics can reduce the environmental impacts of the SSC. One of the advantages of recycling is to reduce the negative environmental impacts due to raw material acquisition from suppliers.

\typo{The} next section presents the computational experiments carried out with a real case study and artificially generated instances.

\section{Computational experiments} \label{sec:experiments}

In this section, we shall discuss two experiments carried out using a case study and artificial instances drawn from the instance generator presented in this paper.
The results obtained by the proposed matheuristic, \emph{AugMathLagr}, are contrasted with two multi-objective methods found in the literature: AUGMECON2 and an adaptation of \emph{AugMathFix}. Note that we refer to the heuristics employed in \emph{AugMathLagr} and in \emph{AugMathFix} to solve mono-objective problems by, respectively, \emph{MathLagr} and \emph{MathFix}. 
The first experiment concerns the case study presented by \citet{Mota2018}. The second presents an analysis of the computational results achieved by the methods for a set of benchmark instances.

The number of grid points \typo{was} set at 10 for all the methods. The \emph{AugMathLagr} {parameters were}: $kMax=10$; $\lambda^{(0)}=0$; 
and $st=1 \times 10^{-10}$, $st=1 \times 10^{-10}$ and $st=1 \times 10^{-8}$ when optimizing the economic, environmental and social functions, respectively.
All experiments were performed on a computer with an Intel Xeon E5-2680v2 2.8 GHz processor and 128 GB of main memory. CPLEX was limited to use only 15GB of memory in all the experiments. Moreover, CPLEX is limited to use 8  threads in both  experiments. A maximum time limit of $1h$ was imposed for CPLEX to solve each mono-objective problem, MOP. Exceptionally for AUGMECON in the experiment with the case study, which \typo{was} computationally more challenging, we imposed a time limit of 3h for CPLEX. Very large and small absolute values are expressed in E notation.

Before detailing the results of the experiments, we shall describe the measures used to evaluate the multi-objective optimization methods.

\subsection{Measures of assessment}

$GAP$, a measure of proximity, was used to evaluate
the quality of {solutions} with respect to {their} optimal value.
$GAP$ quantifies how close a solution $S^{'}$ is from another $S^*$, considering an objective function 
$f_i$ to be minimized, i.e. $f_i \in \{-f_{eco}, f_{env}, -f_{soc}$\} indicated in Equation \eqref{eq:gap}.

\begin{equation}
{GAP(f_i(S^{'}),f_i(S^*))} = \frac{|{f_i}(S^{'}) - {f_i}(S^*)|}{\max({f_i}(S^*), {f_i}(S^{'}))} 
 \label{eq:gap}
\end{equation}

The stopping criterion for CPLEX was when it reached a solution whose $GAP$ between the upper and lower bounds was less than or equal to $1\%$. 
The remaining parameters were set as the default of the solver.

To better assess the results achieved by \textit{MathLagr}, similar to \citet{Tautenhain2019}, we used a set of multi-objective metrics to evaluate the quality of the Pareto frontiers obtained by the methods. They were: the R2 Indicator and two variations of the Mean Ideal Distance (MID) and Spread of Non-Dominated Solutions (SNS).

The R2 Indicator \citep{Hansen1998} checks the quality of a Pareto frontier approximation, $P^A$, by comparing it with a representative set of Pareto frontier, $P^Z$, as shown by Equation \eqref{eq:r2}.
In the experiments reported in this paper, $P^Z$ {was} estimated by AUGMECON2.

\begin{equation}\label{eq:r2}
 I_{R2}(P^A,U) = \frac{1}{|U|}\left( \sum_{{\mu} \in U} \max_{S \in P^Z} \{ {\mu}(S) \} - \sum_{{\mu} \in U} \max_{a \in P^A} \{ {\mu}(a) \} \right)
\end{equation}

\noindent where $U$ is the set of utility functions. A utility function ${\mu} \in U$, ${\mu}: \mathds{R}^{3} \rightarrow \mathds{R}$, maps a solution of the multi-objective problem to a scalar value. Lower values of $I_{R2}(P^A,U)$ indicate better approximations. Negative values of the R2 Indicator express that solutions from $P^A$ are closer to the Ideal point than solutions from $P^Z$. Let the Ideal point be defined
by {$f^I=\{f^*_{eco},f^*_{env},f^*_{soc}\}$}, such that {$f^*_i$} is the optimal value of the problem that optimizes objective function {$f_i \in \{-f_{eco}, f_{env}, -f_{soc}$\}}.

\citet{Hansen1998} suggest {using} the Weighted Sum and weighted \textit{Tchebycheff} utility functions. The Weighted Sum function only takes into account points inside the convex hull of a feasible region and therefore is unsuitable for cases where the solution space is not convex. Thereby, \citet{Brockhoff2012} suggest using the weighted \textit{Tchebycheff} function, formulated as presented in Equation \eqref{eq:tcheb}.

\begin{equation}\label{eq:tcheb}
{u_{{\gamma}}(S) = \max_{{i \in \{eco, env,soc\}}} {\gamma}_{{i}} |{{f_i^*}} - {f_{i}(S)}|}
\end{equation}

\noindent where $\gamma_{eco}, \gamma_{env}$ and $\gamma_{soc}$ are weights associated \typo{with}  $f_{eco}, f_{env}$ and $f_{soc}$, respectively.

For ease of notation, we refer to $I_{R2}(P^A)$ as $I_{R2}(P^A,U)$, with $U$ representing the \textit{Tchebycheff} utility function.

The Mean Ideal Distance (MID) and Spread of Non-Dominated Solutions (SNS) \citep{Behnamian2009} calculate, respectively, the mean and the standard deviation of the distance between the solutions of an approximation of the Pareto frontier, $P^A$, and the Ideal point.  Due to the different scales of the objective function values, {\citet{Tautenhain2019}} employed the adaptation of these measures to calculate the multidimensional $GAP$ {($GAPM$)} instead of the distance between the solutions. Equation \eqref{eq:gapMO} presents how to calculate 
$GAPM$, where $GAP$ is defined according to Equation \eqref{eq:gap}.

 \begin{equation}
 GAPM(f(S^{'}),f^I) = \sqrt{ \sum_{{i \in \{eco, env,soc\} }} \left( GAP({f_{i}}(S^{'}),{f^*_i}) \right )^2} 
 \label{eq:gapMO}
\end{equation}

\noindent where $f(S^{'})=\{-f_{eco}(S^{'}), f_{env}(S^{'}), -f_{soc}(S^{'})\}$.


Equations \eqref{eq:mid} and \eqref{eq:sns} present the adapted $MID$ and $SNS$ metrics\minor{, here called $aMID$ and $aSNS$,} for $P^A$.

\begin{equation}
aMID(P^A) = \sum_{S^{'} \in P^A}\frac{\| GAPM(f(S^{'}), f^I) \| }{ |P^A| }
 \label{eq:mid}
\end{equation}

\begin{equation}
aSNS(P^A) = \sqrt{\sum_{S^{'} \in P^A} \frac{ (aMID(P^A) - \| GAPM(f(S^{'}), f^I) \|)^2 }{ |P^A|-1 } }
 \label{eq:sns}
\end{equation}
Lower values of $aMID(P^A)$ and $aSNS(P^A)$ indicate better solutions in $P^A$.

 \subsection{Experiment I: case study}\label{sec:experiments_cs}

 In the first experiment, a case study of an electronics company is investigated (see the supplementary report on Mendeley Data \cite{Mendeley2020}  or \citep{Mota2018} for additional details {on} this case study). 
  The experiment was carried out with the original problem of \cite{Mota2018}, which \minor{did not consider} suppliers' environmental \minor{impacts} and working conditions. 
 It is worth pointing out that \citet{Mota2018} selected the set of locations that could be used to install the facilities beforehand, which guaranteed adequate working conditions for the employees.

 Table \ref{tab:facility} {only shows details regarding} entities {in} this case study -- location, country GDP per capita where  an entity is located (GDPPC) \citep{Mota2018}, \minor{measured in Purchasing Power Parity (PPP)\footnote{\minor{PPP is an exchange rate that translates the purchasing power of different monetary units.}},}
  and costs in euros related to construction and labor -- in order to contrast such information with the results obtained by the methods. 
The first, second, third and fourth columns \typo{respectively} indicate the entities, their labels, geographical locations and regions. The fifth, sixth and seventh columns respectively report the values of {GDPPC} and construction and labor costs of entities.

\begin{table}[!htb]
\centering
\caption{{Details about the entities of the supply chain.}} 
\small
\label{tab:facility}
\begin{tabular}{|c|c|l|c|c|c|c|}
\hline
\multirow{2}{*}{\textbf{Entity}}   &\multirow{2}{*}{\textbf{Label}}	& \multicolumn{1}{c|}{\multirow{2}{*}{\textbf{Location}}} & \multirow{2}{*}{\textbf{Region}} & \multirow{2}{*}{{GDPPC} \minor{(PPP)}} & \textbf{Construction} & \textbf{Labor} \\
			   &				&			 &	  &	   &	\textbf{cost} \minor{(\euro)} 	& \textbf{cost} \minor{(\euro)} \\ \hline \hline
\multirow{3}{*}{Supplier}  & $S_1$			& Verona, Italy		 & Europe  & 0.98   &	-	& 28.1	\\ \cline{2-7} 
                            & $S_2$           		& Hannover, Germany     & Europe   & 1.24   &	-	& 30.4  \\ \cline{2-7} 
                            & $S_3$                 	& Leeds, United Kingdom	 & Europe   & 1.06   & 	 - 	& 15.3  \\  \hline 
\multirow{3}{*}{Factories}  & $F_1$		 	& Verona, Italy        	 & Europe   & 0.98   & 595  	& 28.1  \\ \cline{2-7} 
                            & $F_2$                  	& Hannover, Germany    	 & Europe   & 1.24   & 661  	& 30.4  \\ \cline{2-7} 
                            & $F_3$                  	& Leeds, United Kingdom	 & Europe   & 1.06   & 601  	& 15.3  \\ \hline 
\multirow{9}{*}{Warehouses} & $W_1$		 	& Verona, Italy	         & Europe   & 0.98   & 595  	& 28.1  \\ \cline{2-7} 
			    & $W_2$                 	& Hannover, Germany     & Europe   & \minor{1.24}   & 661  	& 30.4  \\ \cline{2-7} 
                            & $W_3$ 	            	& Leeds, United Kingdom	 & Europe   & 1.06   & 601  	& 15.3  \\ \cline{2-7} 
                            & $W_4$             	& Zaragoza, Spain       & Europe   & 0.95   & 373  	& 21  \\ \cline{2-7} 
                            & $W_5$	            	& Lisbon, Portugal	 & Europe   & 0.75   & 318 	& 12.2  \\ \cline{2-7} 
                            & $W_6$                  	& São Paulo, Brazil   	 & Brazil   & 0.355  & 538  	& 8.98  \\ \cline{2-7} 
                            & $W_7$                  	& Recife, Brazil	 & Brazil   & 0.355  & 538  	& 8.98  \\ \cline{2-7} 
                            & $W_8$                  	& Budapest, Hungary	 & Europe   & 0.67   & 282  	& 7.5  \\ \cline{2-7} 
                            & $W_9$                  	& Sofia, Bulgaria	 & Europe   & 0.47   & 270  	& 3.7  \\ \hline 
\multirow{7}{*}{Customers}  & $C_1$		 	&  Italy       	 & Europe   & 0.98   & 	-  	& 28.1  \\ \cline{2-7} 
                            & $C_2$                  	& Germany   	 & Europe   & 1.24   & 	-  	& 30.4  \\ \cline{2-7} 
                            & $C_3$                  	& United Kingdom	 & Europe   & 1.06   & 	-   	& 15.3 \\ \cline{2-7}
                            & $C_4$                  	& Spain      	 & Europe   & 0.95   & 	-  	& 21  \\ \cline{2-7} 
                            & $C_5$                  	& Portugal     	 & Europe   & 0.75   & 	-  	& 12.2   \\ \cline{2-7} 
                            & $C_6$                  	& São Paulo, Brazil    	 & Brazil   & 0.355  & 	-  	& 8.98  \\ \cline{2-7} 
                            & $C_7$                  	& Recife, Brazil      	 & Brazil   & 0.355  & 	-  	& 8.98   \\ \hline 
\multirow{4}{*}{Airports}   & $Air_1$			& Zaragoza	 & Europe   & 1.19   & 	-  	& 21  \\ \cline{2-7}                             
                            & $Air_2$                  	& Paris-Charles de Gaulle & Europe & 1.08   & 	-  	& 32.4  \\ \cline{2-7} 
                            & $Air_3$                  	& Kortrijk-Wevelgem      & Europe  & 0.95   & 	-  	& 37.2  \\ \cline{2-7} 
                            & $Air_4$                  	& São Paulo              & Brazil  & 0.355  & 	-  	& 8.98   \\ \hline 
\multirow{2}{*}{Seaports}   & $Sea_1$		 	& Hamburg{, Germany}		  & Europe  & 1.24   & 	-  	& 30.4  \\ \cline{2-7} 
                            & $Sea_2$            	& Santos{, Brazil}		  & Brazil  & 0.355  & 	-  	& 8.98 \\ \hline
\end{tabular}
\end{table}

The supply chain entities are located in several countries in Europe and Brazil. The company {has its} factory $F_1$ and warehouse $W_1$ installed in Verona. Since the suppliers and factories {are located} in Europe, the company must use {either air} or sea transportation to {carry} the final products to Brazil.

The company produces only two types of final products from four different raw materials. The returned final products are modeled as recovered products. 
There are also four technologies available for factories to manufacture products from raw materials and two technologies for recovering used products into final products. 

  Table \ref{tab:exp1_payoff} depicts the results of the exact {method} employed in AUGMECON2, of  \emph{MathLagr} and of \emph{MathFix} in the optimization of each objective function of the problem.

  The third, fourth and fifth columns indicate the results {obtained by} optimizing the economic, environmental {and} social functions. For each method, the results reported are the computational running time to find the solution and the values of economic, {$f_{eco}$}, environmental, {$f_{env}$}, and social, {$f_{soc}$}, functions. 
  To calculate the $GAPs$ for \emph{MathFix} and \emph{MathLagr} according to Equation \eqref{eq:gap}, we assume the solutions achieved by the exact method to be the optimal ones.
  Moreover, the rows identified by ``Speed up'' report the ratio of the exact method execution time and the matheuristics execution time.
    All the methods estimate the ideal and nadir points {using}, respectively, the best and worst values of each objective solution.

\begin{table}[!htb]
\begin{center}
\caption{Optimization of the individual objective functions by an exact strategy, \emph{MathLagr} and \emph{MathFix}.  } 
\label{tab:exp1_payoff}
\small
\begin{tabular}{cc|c|r|r|r|}
\cline{4-6}
\multicolumn{1}{l}{\textbf{}}             &     \multicolumn{1}{l}{\textbf{}}                     &                        & \multicolumn{3}{c|}{{\textbf{Function to be optimized}}}                                                                      \\ [0.25cm]  \cline{4-6} 
\multicolumn{1}{l}{\textbf{}}         &  \multicolumn{1}{l}{\textbf{}}   &                   & \multicolumn{1}{c|}{{$\max {f_{eco}}$}} & \multicolumn{1}{c|}{{${\min} {f_{env}}$}}  &  \multicolumn{1}{c|}{{$\max {f_{soc}}$}} \\ [0.25cm]  \hline 
\multicolumn{1}{|c|}{\multirow{4}{*}{{Exact method}}}                       & \multicolumn{2}{c|}{Time (s)}              & 

154.218             & 15.676              & 3742.222           \\ \cline{2-6} 
\multicolumn{1}{|c|}{}   &       Values  of the           & ${f_{eco}} $ & 
\minor{1.405e+09}         & \minor{1.052e+09}         & \minor{-5.713e+09}       \\ \cline{3-6} 
\multicolumn{1}{|c|}{}   &         objective                 & ${f_{env}} $  &
\minor{9.766e+08}         & \minor{9.108e+08}         & \minor{9.249e+08}        \\ \cline{3-6} 
\multicolumn{1}{|c|}{}   &          functions                & ${f_{soc}}$ &
895.804          & 2176.588         & 40800.432       \\  \hline
\multicolumn{1}{|c|}{\multirow{6}{*}{\emph{MathLagr}}}   & \multicolumn{2}{c|}{Time (s)}              &
43.035              & 22.469              & 150.603            \\ \cline{2-6}  
\multicolumn{1}{|c|}{}  &      Values  of the                    & ${f_{eco}}$ &
 \minor{1.404e+09}         & \minor{8.906e+08}         & \minor{-5.646e+09}       \\ \cline{3-6}  
\multicolumn{1}{|c|}{}   &    objective                    & ${f_{env}}$  & 
\minor{9.774e+08}         & \minor{9.108e+08}         &  \minor{9.721e+08}        \\ \cline{3-6} 
\multicolumn{1}{|c|}{}       &   functions                     & ${f_{soc}}$ &
943.260          & 2573.369        & 40432.385      \\  \cline{3-6}   
\cline{2-6} 
\multicolumn{1}{|c|}{}                                                             & \multicolumn{2}{c|}{$GAP$ (\%)}               & 
 0.041\%            &  \minor{2.273e-3\%}            & 0.902\%           \\ \cline{2-6}  
\multicolumn{1}{|c|}{}                                                             & \multicolumn{2}{c|}{Speed up}               &
3.584               & 0.698               & 24.848             \\ \hline  
\multicolumn{1}{|c|}{\multirow{6}{*}{\emph{MathFix}}} & \multicolumn{2}{c|}{Time (s)}               & 
51.515              & 34.231              & 152.861            \\ \cline{2-6} 
\multicolumn{1}{|c|}{}                                     &     Values  of the                   & ${f_{eco}}$ &
\minor{1.401e+09}         & \minor{1.040e+09}         & \minor{-5.669e+09}       \\ \cline{3-6} 
\multicolumn{1}{|c|}{}                                     &       objective                 & ${f_{env}}$  & 
\minor{9.820e+08}         & \minor{9.108e+08}         & \minor{9.288e+08}        \\ \cline{3-6} 
\multicolumn{1}{|c|}{}                                   &      functions                    & ${f_{soc}}$ & 
828.169          & 2038.037         & 40532.255       \\ \cline{3-6}
\cline{2-6} 
\multicolumn{1}{|c|}{}                                                             & \multicolumn{2}{c|}{$GAP$ (\%)}       &
0.230            & \minor{6.038e-3}             &  0.657         \\\cline{2-6} 
\multicolumn{1}{|c|}{}                                                             & \multicolumn{2}{c|}{Speed up}               &
2.994               & 0.458               & 24.481         \\ \hline         
\end{tabular}
\end{center}
\end{table}

\emph{MathLagr} is approximately \minor{$3.5$ and $25$} times faster than  the exact method  for optimizing  the economic and social objective functions, respectively. Despite being slower than \emph{MathLagr}, \emph{MathFix} was still \minor{$3$ and $24$}
times faster than  the exact method for optimizing the economic and social objective functions, respectively. AUGMECON was faster than \emph{MathLagr} and \emph{MathFix} to optimize the environmental objective function.

\emph{MathLagr} and \emph{MathFix} obtained solutions for the economic and environmental objective functions whose $GAP$s to the corresponding solutions obtained by  the exact method were lower than $0.041\%$ and $0.230\%$, respectively.
On the one hand, \emph{MathLagr} was faster and obtained a smaller $GAP$ than \emph{MathFix}  regarding the economic and environmental function optimization solution. On the other, \emph{MathFix}  obtained a solution for the social objective function whose $GAP$ to the solution obtained by the exact method was only $0.657\%$.

 
  Table \ref{tab:exp1_metrics} presents the numbers of non-dominated solutions found by the three multi-objective methods AUGMECON2, \emph{AugMathLagr} and \emph{AugMathFix} as well as the $aMID$, $aSNS$, {and R2} Indicator values. It also reports the total computational running times in seconds that the methods took to approximate the Pareto frontier for the case study.

\begin{center}
\begin{table}[!htb]
\caption{Results of the multi-objective metrics for the Pareto frontier approximations obtained by AUGMECON2, \emph{AugMathLagr} and \emph{AugMathFix}. 
} 
\label{tab:exp1_metrics}
\small
\begin{center}
\begin{tabular}{c|c|r|r|r|r|}
\cline{2-6}
\multicolumn{1}{l|}{}            & \multicolumn{1}{c|}{\textbf{Number}} & \multicolumn{1}{c|}{$\mathbf{aMID}$}& \multicolumn{1}{c|}{$\mathbf{aSNS}$} & \multicolumn{1}{c|}{$\mathbf{R2}$} & \multicolumn{1}{c|}{\textbf{{Time (h)}}} \\ \hline
\multicolumn{1}{|c|}{AUGMECON2}         & {31}       & {12.177}  &  {3.010}       & 0.000                   & {153.359}                 \\ \hline
\multicolumn{1}{|c|}{\emph{AugMathLagr}}   & 31      & 13.252          & 2.943      &  \minor{5.558e-04}                 &      2.656               \\ \hline
\multicolumn{1}{|c|}{\emph{AugMathFix}} & 26       & 11.076          & 3.021       & \minor{-6.096e-04}                & 4.913                       \\ \hline
\end{tabular}
\end{center}
\end{table}
\end{center}

Table \ref{tab:exp1_metrics} {shows that}  \emph{AugMathLagr} and \emph{AugMathFix} obtained solutions with lower $aMID$ and $aSNS$ values than {AUGMECON2}. 
In addition, {the R2} indicator of the solutions found by \emph{AugMathLagr} and \emph{AugMathFix} were low. These results attest \typo{to} the good quality of the Pareto frontier approximations found by the methods when compared to AUGMECON2.

On the one hand, \emph{AugMathLagr}  is \minor{approximately $57$ and $1.8$} times faster than AUGMECON2 and \emph{AugMathFix}, respectively, to estimate the Pareto frontier. On the other, both AUGMECON2 and \emph{AugMathLagr} find more non-dominated solutions than \emph{AugMathFix}.

\begin{figure}[!htb]
\subfigure[fig:exp1paretoaug][{Projection over $f_{eco}$ and $f_{env}$.}]{\label{fig:exp1_pareto_proj_eco_env} \includegraphics[width=0.32\textwidth]{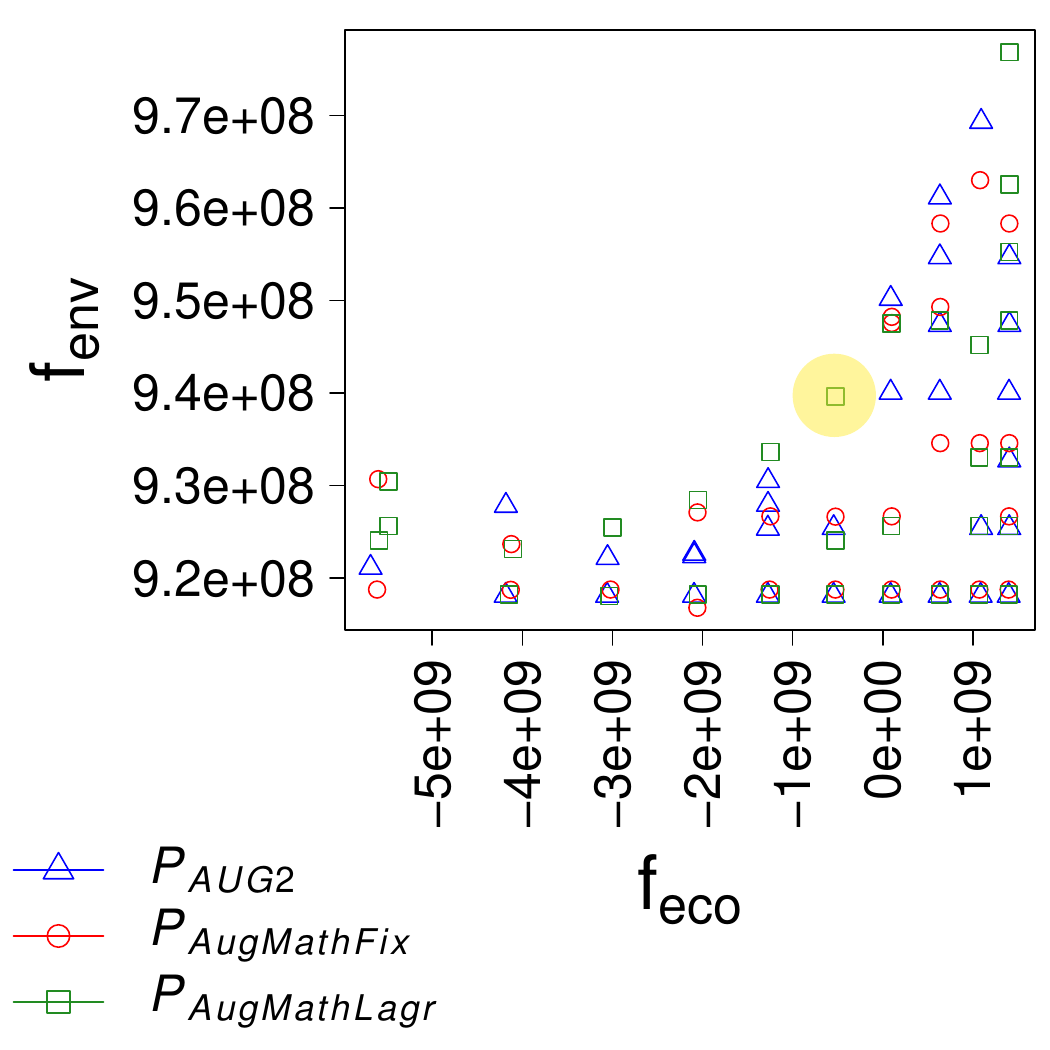}}
\subfigure[fig:exp1paretolagr][{Projection over $f_{eco}$ and $f_{soc}$.}]{\label{fig:exp1_pareto_proj_eco_soc}
\includegraphics[width=0.32\textwidth]{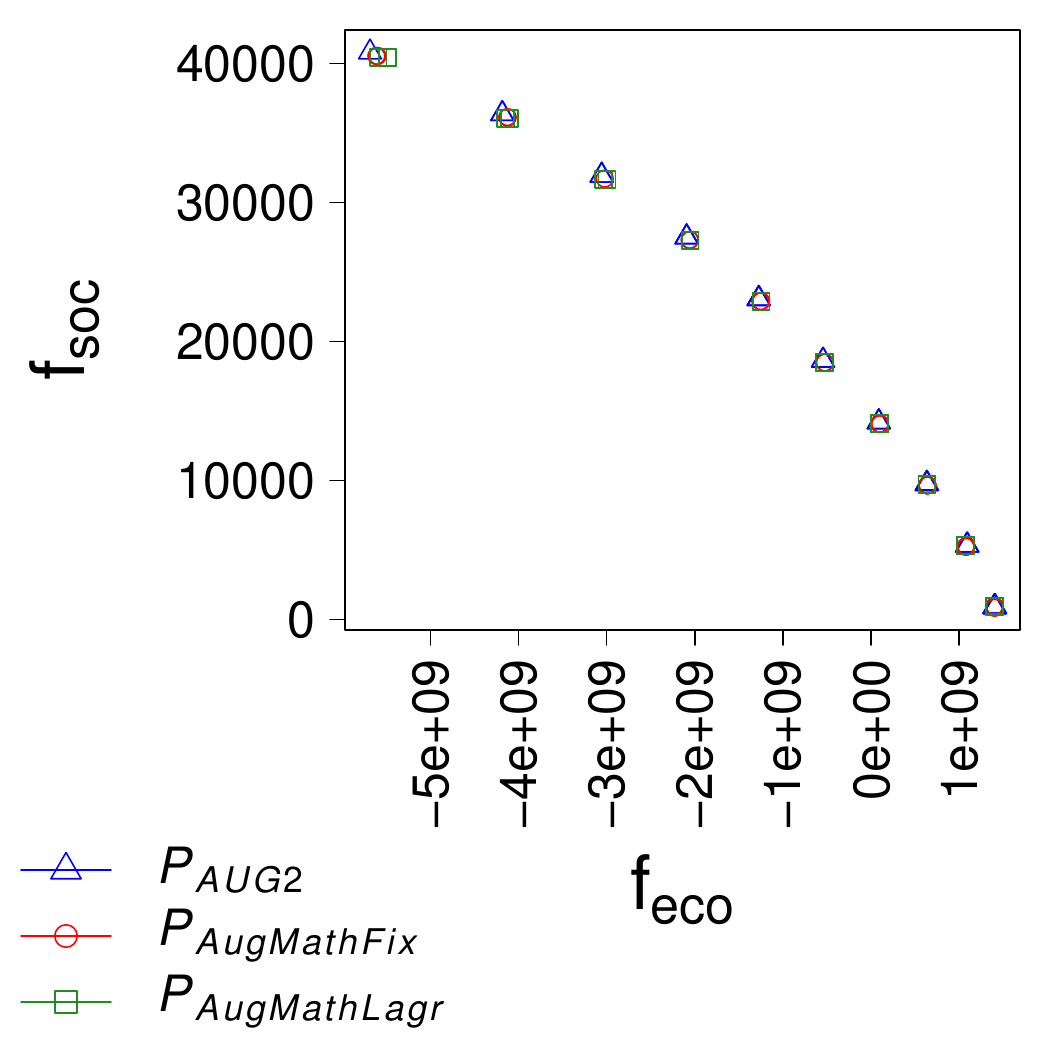}}
\subfigure[fig:exp1paretofix][{Projection over $f_{env}$ and $f_{soc}$.}]{\label{fig:exp1_pareto_proj_env_soc}
\includegraphics[width=0.32\textwidth]{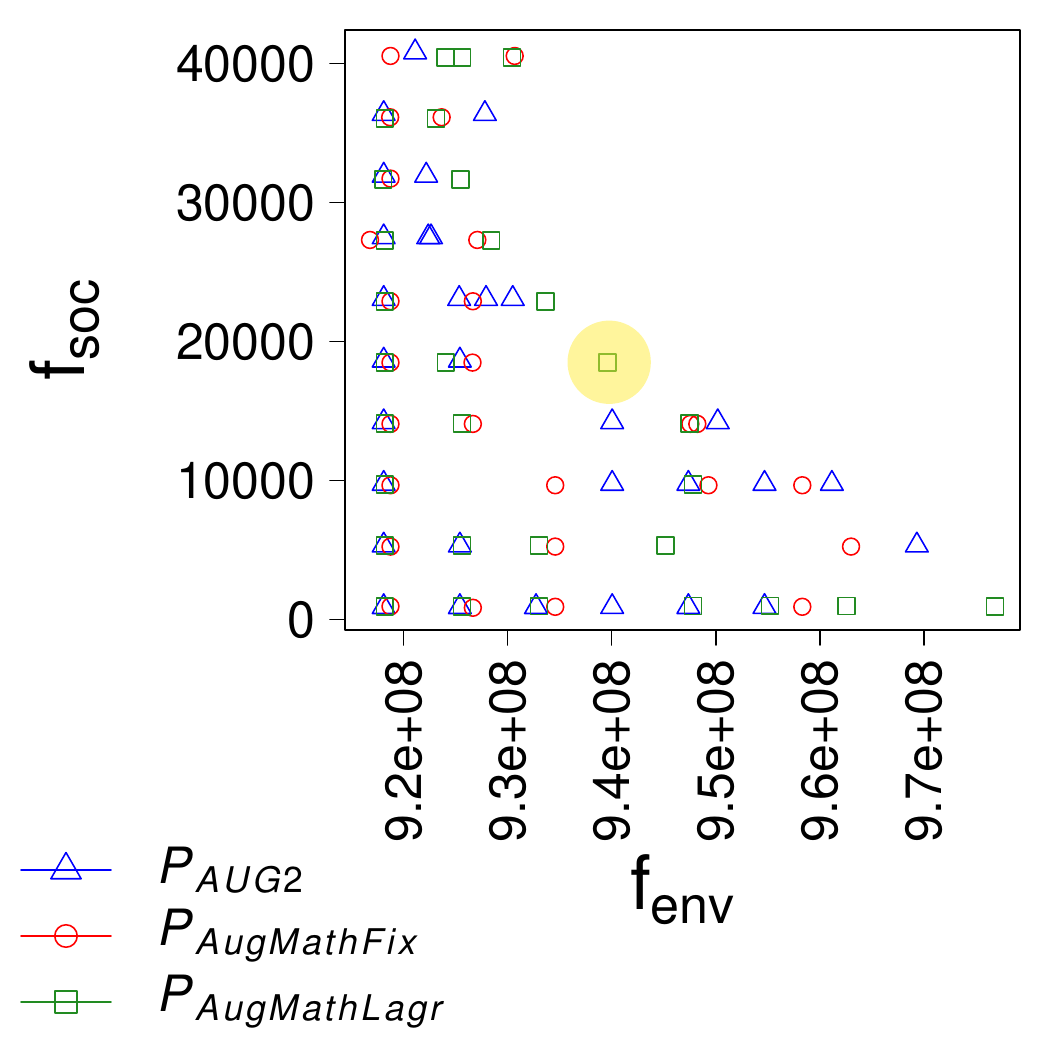}}
\caption{\label{fig:example_pareto}{Projections of the Pareto frontier approximations obtained by the methods for the case study.  }}
\end{figure}

Figure \ref{fig:example_pareto} exhibits projections of the objective \minor{function} values of the solutions in the Pareto frontier approximations  achieved by  AUGMECON2, \emph{AugMathLagr} and \emph{AugMathFix} into two dimensional spaces.
 
 Despite \emph{AugMathLagr} and \emph{AugMathFix} finding very close objective function  values, Figures \ref{fig:exp1_pareto_proj_eco_env} and \ref{fig:exp1_pareto_proj_env_soc} show that \emph{AugMathFix} did not obtain some solutions achieved by \emph{AugMathLagr}.
In particular, \emph{AugMathLagr} found a non-dominated solution  whose $f_{env}$ value is  9.4e+08, highlighted in yellow circles in Figures \ref{fig:exp1_pareto_proj_eco_env} and \ref{fig:exp1_pareto_proj_env_soc}. No non-dominated solution close to such a solution was obtained by \emph{AugMathFix}.

\minor{In Figures \ref{fig:exp1_pareto_proj_eco_env} and \ref{fig:exp1_pareto_proj_env_soc}, it is not possible to notice a clear trade-off between the environmental  and the other objectives. Figure \ref{fig:exp1_pareto_proj_eco_soc}, however, more clearly illustrates the trade-off between the economic ($f_{eco}$) and social ($f_{soc}$) objectives, since the higher the economic function, the lower the social value. Thereby, for a given environmental objective, the higher the $f_{eco}$, the lower the $f_{soc}$. Besides, in Figures \ref{fig:exp1_pareto_proj_eco_env} and \ref{fig:exp1_pareto_proj_env_soc}, there are some solutions with the same $f_{soc}$ whose $f_{eco}$ values, despite being very close, ensure that the solutions are non-dominant.}

\minor{In this context, Figure \ref{fig:example_pareto} {graphically demonstrates} the good quality  of the Pareto frontier approximations  found by the methods, since it is possible to visualize that the non-dominated solutions are well distributed in the hyperplane.
Multi-objective optimization can add value to the supply chain by allowing decision makers to select the most suitable solutions, according to their goals. When presenting solutions to decision makers, the improvement and worsening of each objective should be taken into account. For example, a slight decrease in the economic objective might represent a significant improvement in the environmental objective.}


Table \ref{tab:exp1_cs_inst} {shows} which  factories and warehouses are opened in the solutions found by \emph{MathLagr}, \emph{MathFix} and the exact method in AUGMECON2 to optimize the economic, environmental and social criteria. Additionally, {in this table we} point out the entities installed at maximum and minimum {capacities}. The marks $A,L$ and $M$ indicate, respectively, the exact method, \emph{MathLagr},  and \emph{MathFix}.

\begin{table}[htb]
\centering
\caption{{Case study: location and installation capacities of factories and warehouses in the solutions {obtained by} the exact method, \emph{MathLagr} and \emph{MathFix}.}}
\label{tab:exp1_cs_inst}
\small
\tabcolsep=0.05cm
\begin{tabular}{cc|c|c|c|c|c|c|c|c|c|}
\cline{3-11}
			&		&  \multicolumn{9}{c|}{\textbf{{Function to be optimized}}} \\ [0.25cm]   \cline{3-11}
			&   		& \multicolumn{3}{c|}{$\max {f_{eco}}$}                          & \multicolumn{3}{c|}{${\min {f_{env}}}$}                        & \multicolumn{3}{c|}{${\max {f_{soc}}}$}                            \\ \cline{3-11} 
			&   		& Location   	& Maximum & Minimum 	 & Location   	& Maximum 	& Minimum    & Location   	& Maximum 	& Minimum  \\                                
			& 		&		            & capacity	    & capacity	    &		    & capacity  	& capacity	&		            & capacity	    & capacity  \\ \hline                                 
\multicolumn{1}{|c|}{\multirow{3}{*}{Factories}}
			           &  $F_1$ 	& $A,L,M$	& - 	  & -       & $A,L,M$  & -	      &	-    & $A,L,M $   	& -     & -	\\ \cline{2-11}	 
\multicolumn{1}{|c|}{}	& $F_2$ 	& $A,L,M$ &	          &		    & $A,L;M$  & $A,L,M$  &		 & $A,L,M $   	& $A,L,M$ 	&	\\ \cline{2-11}	    
\multicolumn{1}{|c|}{}	& $F_3$ 	& $A,L,M$ & $A,L,M$   &		    & $A,L,M$  & $A,L,M$  &  	 & $A,L,M $   	& $A,L,M$ 	&	\\ \hline 
\multicolumn{1}{|c|}{\multirow{9}{*}{Warehouses}} 
			            & $W_1$		& $A,L,M$	&  -      &	-       &  $A,L,M $&  -        &   -  & $A,L,M $ & -     & - \\\cline{2-11}	 	
\multicolumn{1}{|c|}{}	& $W_2$     &           &	      &         &  $L,M $  &  $L,M$  &       &$A$    & $A$	  &     	\\ \cline{2-11}	   
\multicolumn{1}{|c|}{}	& $W_3$ 	&           &	      &   	    &  $A,L,M $& $A,L,M $  &     & $A,L,M$   & $A,L,M$	  &       	\\ \cline{2-11}	 	
\multicolumn{1}{|c|}{}  & $W_4$ 	&        	&	      &   	    &  $A$    &	 $A$      &      &$A,L,M $ & $A,L,M $  &       	\\ \cline{2-11}	    
\multicolumn{1}{|c|}{}	& $W_5$ 	& $M$      	&	      &         &  $A$    &	 $A$      &      &$A,L,M $ & $A,L,M $  &       	\\ \cline{2-11}	 	
\multicolumn{1}{|c|}{}	& $W_6$ 	&        	&	      &     	&  $A,L $ &	 $A,L $   &      &$A,L,M $ & $A,L,M $ &       	\\ \cline{2-11}	   
\multicolumn{1}{|c|}{}	& $W_7$ 	&     	    &	      &         &  $L,M$  &	 $L,M$    &      &$A,L,M $ & $A,L,M $  &        \\ \cline{2-11}	  
\multicolumn{1}{|c|}{}	& $W_8$ 	&       	&	      &   	    &         &	          &      &$A,M $ & $A,M $ &       \\ \cline{2-11}	 
\multicolumn{1}{|c|}{}	& $W_9$ 	& $A,L$     &         &   	    &         &	          &      &$A,L,M $ & $A,L,M $  &        	\\ \hline      
 
\end{tabular}
\end{table}

It is possible to observe in Table \ref{tab:exp1_cs_inst} that there is no difference in the selection of factories to be opened at maximum or minimum capacity in the solutions found by any of the methods.
%
%
Both the exact method and \emph{MathLagr} opened warehouse $W_9$ in addition to warehouse $W_1$, which is opened at a fixed capacity in every solution, in the solutions that optimize the economic function. 
\emph{MathFix}, on the other hand, selected warehouse $W_5$ besides $W_1$.

The solutions obtained by the three methods for the environmental minimization problem open more warehouses  than  the economic objective maximization solutions.

When minimizing the environmental objective function, the most relevant decisions for the environmental impact are those related to manufacturing and remanufacturing. Such decisions correspond to approximately $99.94\%$ of the environmental impacts on the solution that minimizes the environmental objective function. Transportation and entity location decisions  represent only \minor{$0.4\%$ and $0.2\%$} of the overall environmental impacts on such a solution, respectively.

 By comparing the environmental impacts due to the installation of additional warehouses on the solution  of the mono-objective problem that minimizes the environmental objective function, they were 43.43\% greater than such impacts on the solution that maximizes the economic objective. 
Nonetheless, the environmental impacts due to manufacturing and remanufacturing on the environmental minimization solution are \minor{approximately $7\%$}
lower than such impacts on the economic optimization solution. Note that, as mentioned earlier, manufacturing and remanufacturing represent the overwhelming majority of the environmental impacts. As a consequence, this difference in the overall solution is \minor{approximately $7\%$.}

\minor{From these results, we can see that the model and the methods can properly identify the manufacturing and remanufacturing activities as the most representative activities for the environmental impact. Thereby,  it is recommended for the company to perform a careful analysis and estimations on the parameters that affect both the manufacturing and remanufacturing technology costs and environmental impacts. We also recommend performing a sensitivity analysis on these parameters when dealing with different industrial segments.}

Except for warehouse {$W_1$}, required by the case study to be opened at a fixed capacity in every solution, all the remaining warehouses are opened at maximum capacity in the solutions obtained by the exact method through the optimization of  the social function. \emph{MathLagr} and \emph{MathFix} opened, respectively, 7 and 8 out of the 9 available warehouses when optimizing the social function.

\minor{The social solution benefits from more entities and transportation. We suggest, thereby, to analyze in detail choices associated with solutions that benefit the social criterion. Moreover, companies might have strategic goals or preferences in selecting certain entities. Therefore, it is advisable to take these issues into account when selecting the best trade-off solution involving social optimization. }

All the solutions obtained by the methods to individually optimize the objective functions use the two seaports.
In every social optimization solution, all  airports are operational, because they positively contribute to creating jobs. 
Because air transportation is more expensive than road transportation, the airports are not employed in the solution that optimizes the economic objective function.
In this case study, air transportation causes lower environmental impacts than one of the trucks used in road transportation \citep{Mota2018}. Thereby, it is preferred when optimizing the environmental function. Airports $Air_1$ and $Air_3$ are operational in all the solutions that optimize this function. In the solution found by \emph{MathLagr}, airport $Air_2$ operates as well.

\minor{Even though the manufacturing and remanufacturing represent the overwhelming majority of the environmental impacts, all the methods still prefer air transportation due to its lower environmental impacts.}

\subsection{Experiment II: artificial instances}
{In} this experiment, we created artificial instances through the introduced instance generator. For this purpose, we considered instances with a predefined number of entities, items and technologies. We generated four types of instances, differentiated by the prefix in their labels ``STD'', ``TECHC'', ``RAWC'', ``SUP'' and ``CAP''. The instances prefixed by ``STD''  follow the same patterns of the model introduced in \cite{Mota2018}.
In the instances prefixed by ``TECHC'', the technology acquisition costs {were} the same for all technologies.  The manufacturing costs {were} the same for all the suppliers in the instances prefixed by ``RAWC''. In the instances prefixed by ``SUP'', no minimum order {was} imposed for the suppliers.  There {was} no minimum use of technologies in the instances prefixed by ``CAP''. The corresponding remaining parameters of the instances prefixed by  ``STD'', ``TECHC'', ``RAWC'', ``SUP'' and ``CAP'' {were} generated as described in Section 3 of the Supplementary Material.

The generated set {consisted} of {$60$} artificial instances, whose primary characteristics are summarized in Table  \ref{tab:exp2_sets_horizon}. \minor{The first and second columns of this table indicate the prefixes and suffixes for the names of the instances, respectively. The name of an instance is given by its prefix followed by the suffix.}
Columns ``$|I|$'', ``$|M_{rm}|$'', ``$|M_{fp}|$'', ``$|M_{rp}|$'' and ``$T$'' express, respectively, the number of entities, raw materials, final products and periods  of the instances.

Instances with 17 entities have 3 suppliers, 3 factories, 3 warehouses, 4 customers, 2 airports and 2 seaports. Instances with 25 entities have 3 suppliers, 3 factories, 9 warehouses, 4 customers, 4 airports and 2 seaports.  The number of production and remanufacturing technologies {was} 3, totaling 6 technologies.  For ease of identification, we included in the name of the instances information regarding the number of entities, items and periods in the planning horizon and specific patterns for the parameters, as described in the following paragraph. 
The supplementary report in \citep{Mendeley2020} presents the values chosen for the parameters to generate such instances discussed in Section 3 of the Supplementary Material.

\begin{table}[htb]
\centering
\caption{{Number of entities, items and periods in each set of instances.}}
\label{tab:exp2_sets_horizon}
\small
\tabcolsep=0.1cm
\begin{tabular}{|c|c|c|c|c|c|c|}
\hline
\multicolumn{2}{|c|}{Instance} &  \multirow{2}{*}{{$|I|$}} &  \multirow{2}{*}{{$|M_{rm}|$}} &  \multirow{2}{*}{{$|M_{fp}|$}} &  \multirow{2}{*}{{$|M_{rp}|$}}  &  \multirow{2}{*}{{$|T|$}}              \\  
    \cline{1-2} 
{{Prefixes}}	& {{Suffix}} &  &  &  &   &             \\ 
\hline

 STD, TECHC, RAWC, SUP, CAP & \_I17\_M4\_T3 	& 17 & 2 & 1  & 1& 3                    \\ \hline
 STD, TECHC, RAWC, SUP, CAP, & \_I17\_M4\_T6 	& 17 & 2 & 1  & 1 & 6                \\ \hline
  STD, TECHC, RAWC, SUP, CAP & \_I17\_M4\_T12 & 17 & 2 & 1  & 1 & 12               \\ \hline 
 STD, TECHC, RAWC, SUP, CAP & \_I17\_M8\_T3 	& 17 & 4 & 2  & 2 & 3                    \\ \hline
 STD, TECHC, RAWC, SUP, CAP & \_I17\_M8\_T6 	& 17 & 4 & 2 & 2 & 6                \\ \hline
 STD, TECHC, RAWC, SUP, CAP & \_I17\_M8\_T12 	& 17 & 4 & 2 &  2 & 12               \\ \hline
 STD, TECHC, RAWC, SUP, CAP & \_I25\_M4\_T3 	& 25 & 2 & 1  & 1 & 3                    \\ \hline
 STD, TECHC, RAWC, SUP, CAP & \_I25\_M4\_T6 	& 25 & 2 & 1 & 1 & 6                \\ \hline
 STD, TECHC, RAWC, SUP, CAP & \_I25\_M4\_T12 	& 25 & 2 & 1  & 1 & 12               \\ \hline 
 STD, TECHC, RAWC, SUP, CAP & \_I25\_M8\_T3 	& 25 & 4 & 2  & 2 & 3                    \\ \hline
 STD, TECHC, RAWC, SUP, CAP & \_I25\_M8\_T6 	& 25 & 4 & 2  & 2 & 6                \\ \hline
 STD, TECHC, RAWC, SUP, CAP & \_I25\_M8\_T12 	& 25 & 4 & 2 & 2 & 12               \\ \hline
\end{tabular}
\end{table}

\begin{figure}[!htb]
\subfigure[fig:midI17M4][{$|I|=17$ and $|M|=4$.}
]{ \includegraphics[width=0.234\textwidth]{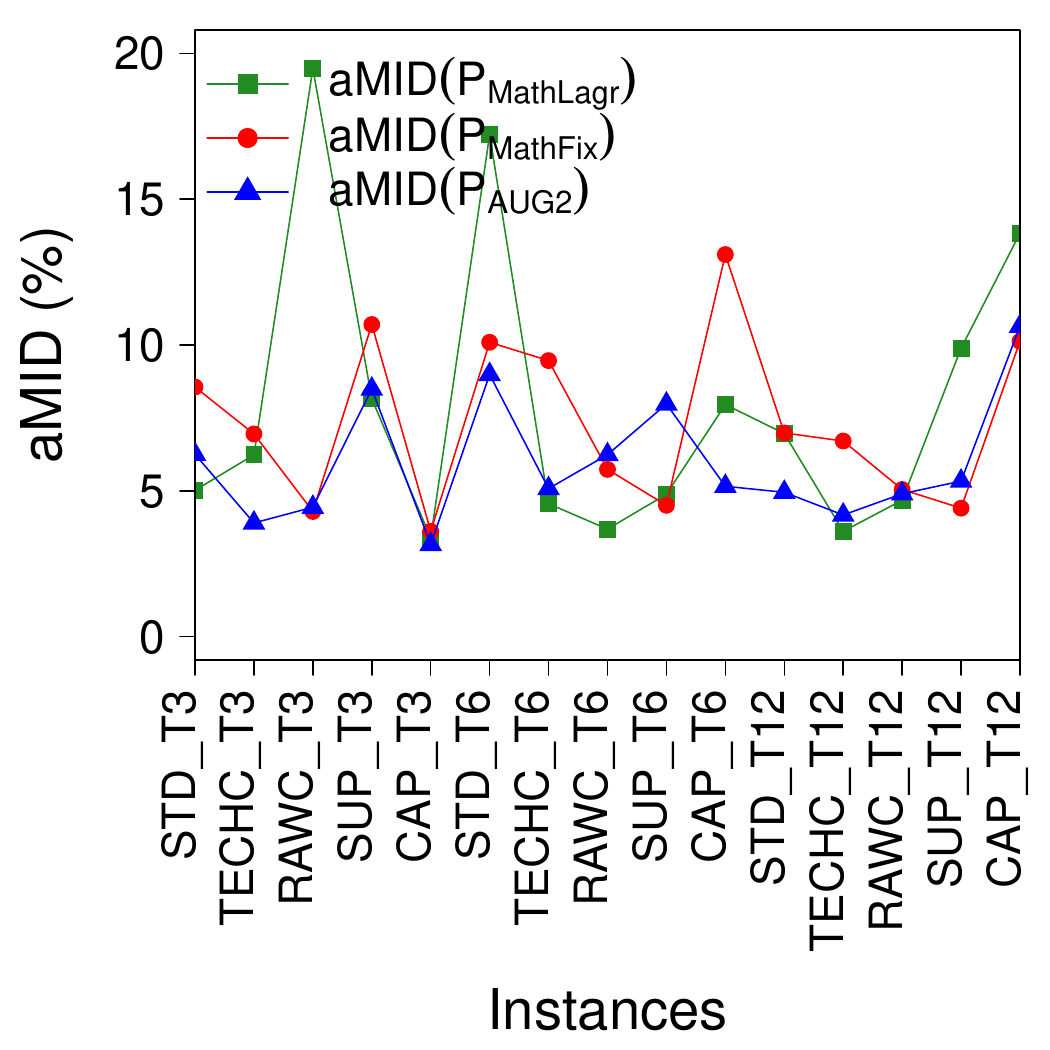}}
\subfigure[fig:midI17M8][{$|I|=17$ and $|M|=8$.}
]{ \includegraphics[width=0.234\textwidth]{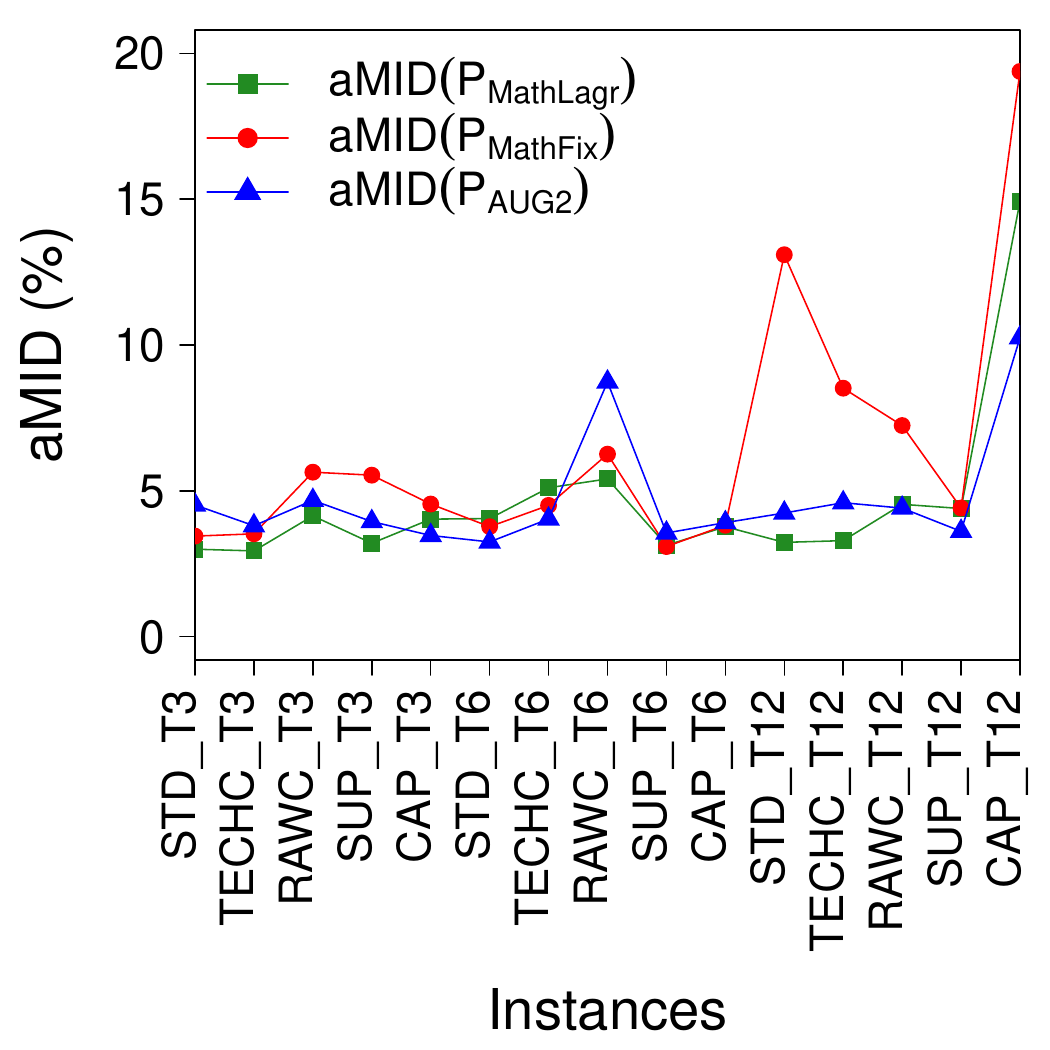}}
\subfigure[fig:midI25M4][{$|I|=25$ and $|M|=4$.}
]{ \includegraphics[width=0.234\textwidth]{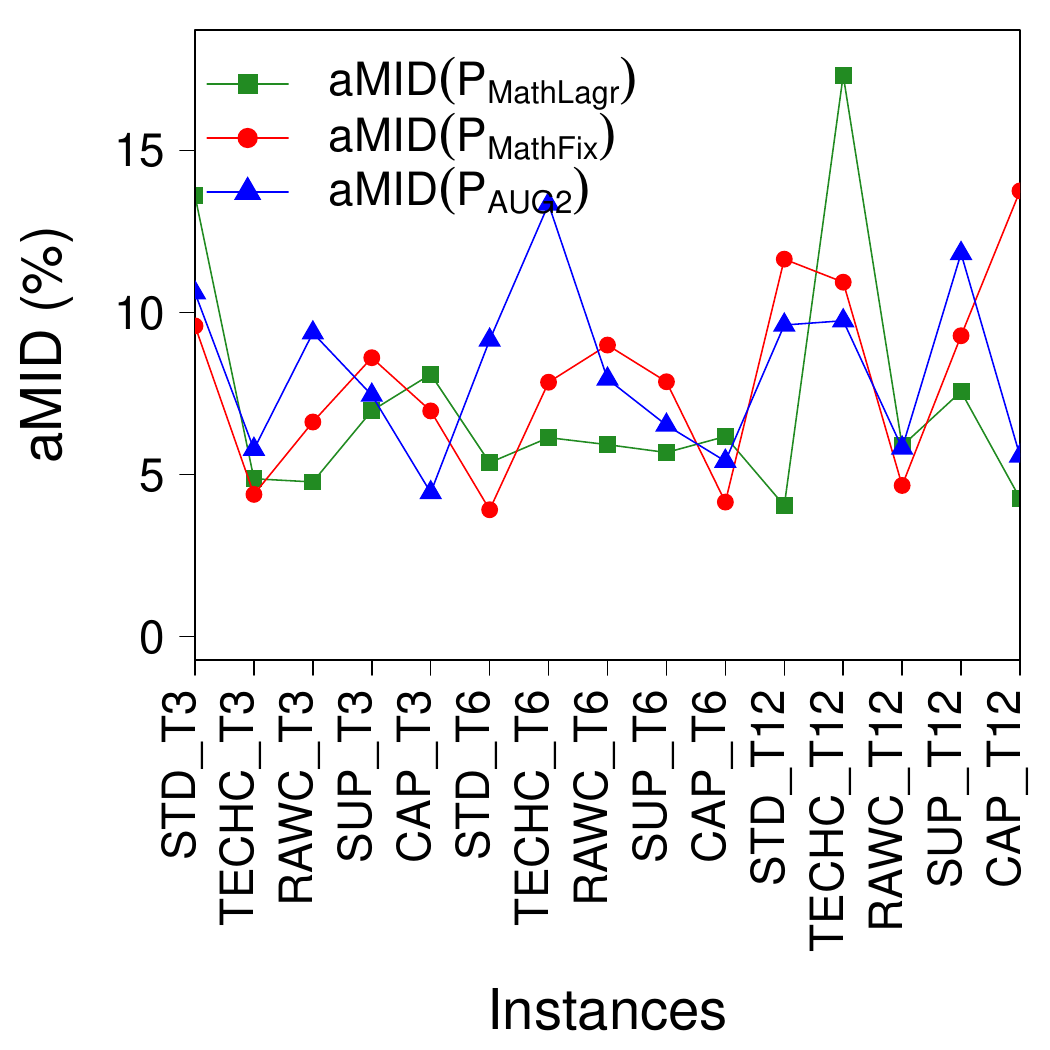}}
\subfigure[fig:midI25M8][{$|I|=25$ and $|M|=8$.}
]{ \includegraphics[width=0.234\textwidth]{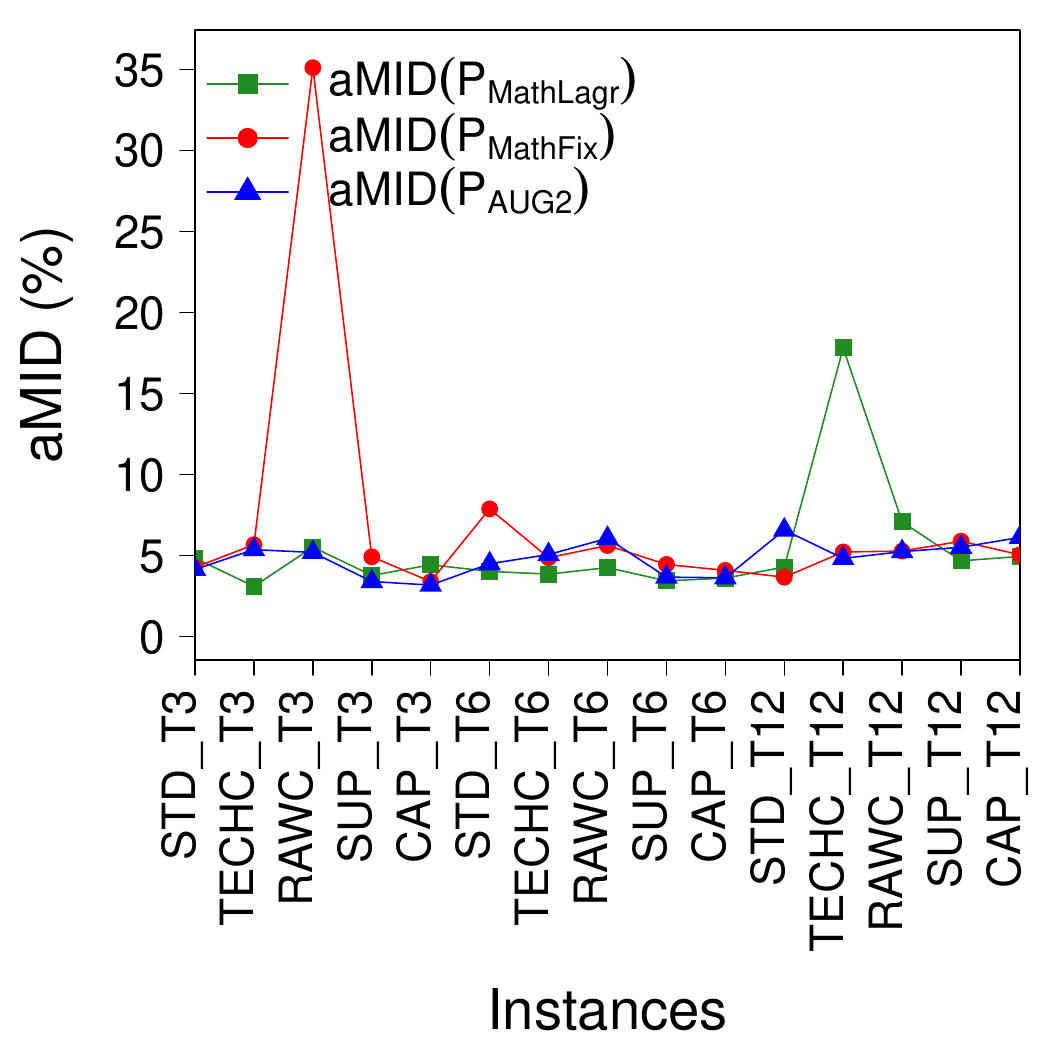}}
       \caption{\typo{Measure $aMID$} for $P_{MathLagr}$, $P_{MathFix}$ and $P_{AUG2}$. 
       }
      \label{fig:metricsM_aMID}      
\end{figure}
\begin{figure}[!htb]
 \subfigure[fig:snsI17M4][{$|I|=17$ and $|M|=4$.}]{ 
 \includegraphics[width=0.234\textwidth]{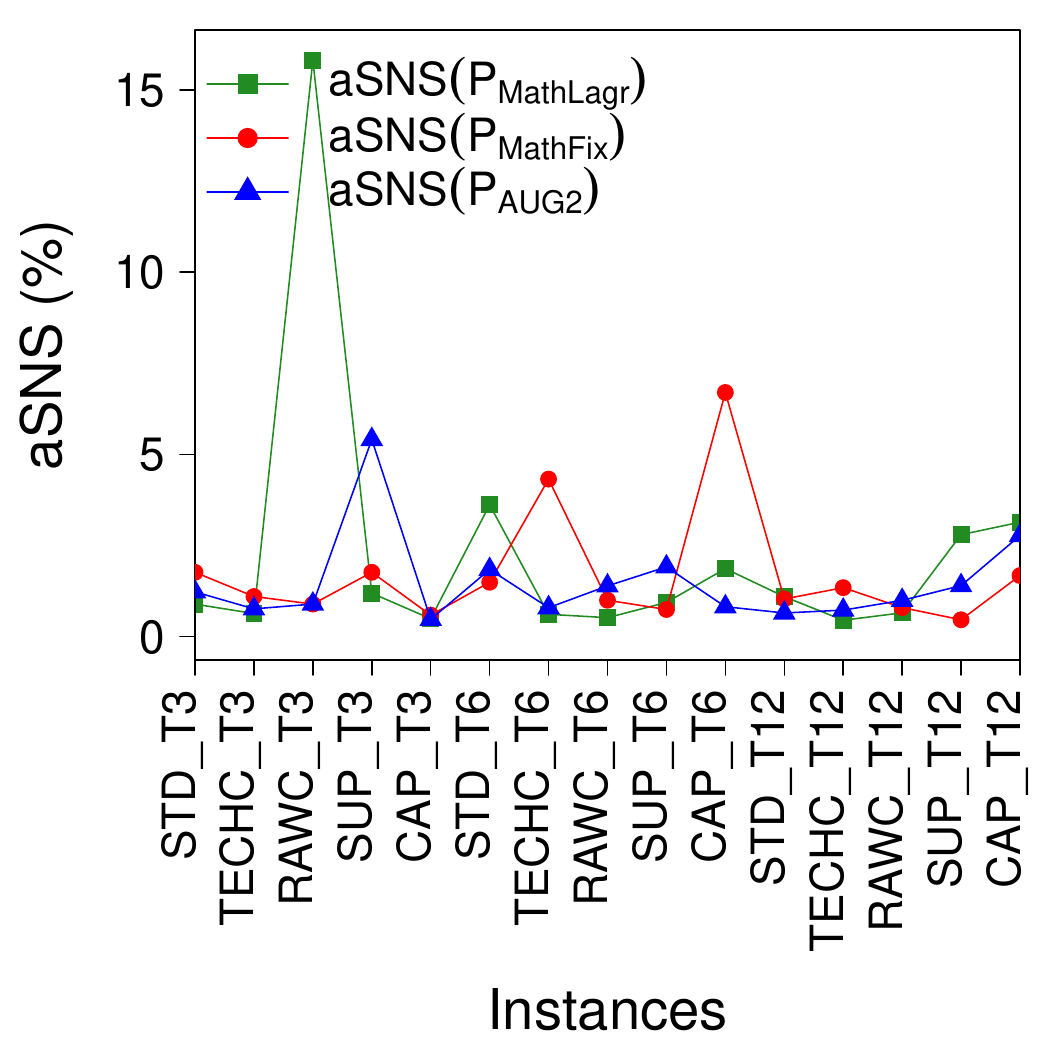}}
\subfigure[fig:snsI17M4][{$|I|=17$ and $|M|=8$.}]{ \includegraphics[width=0.234\textwidth]{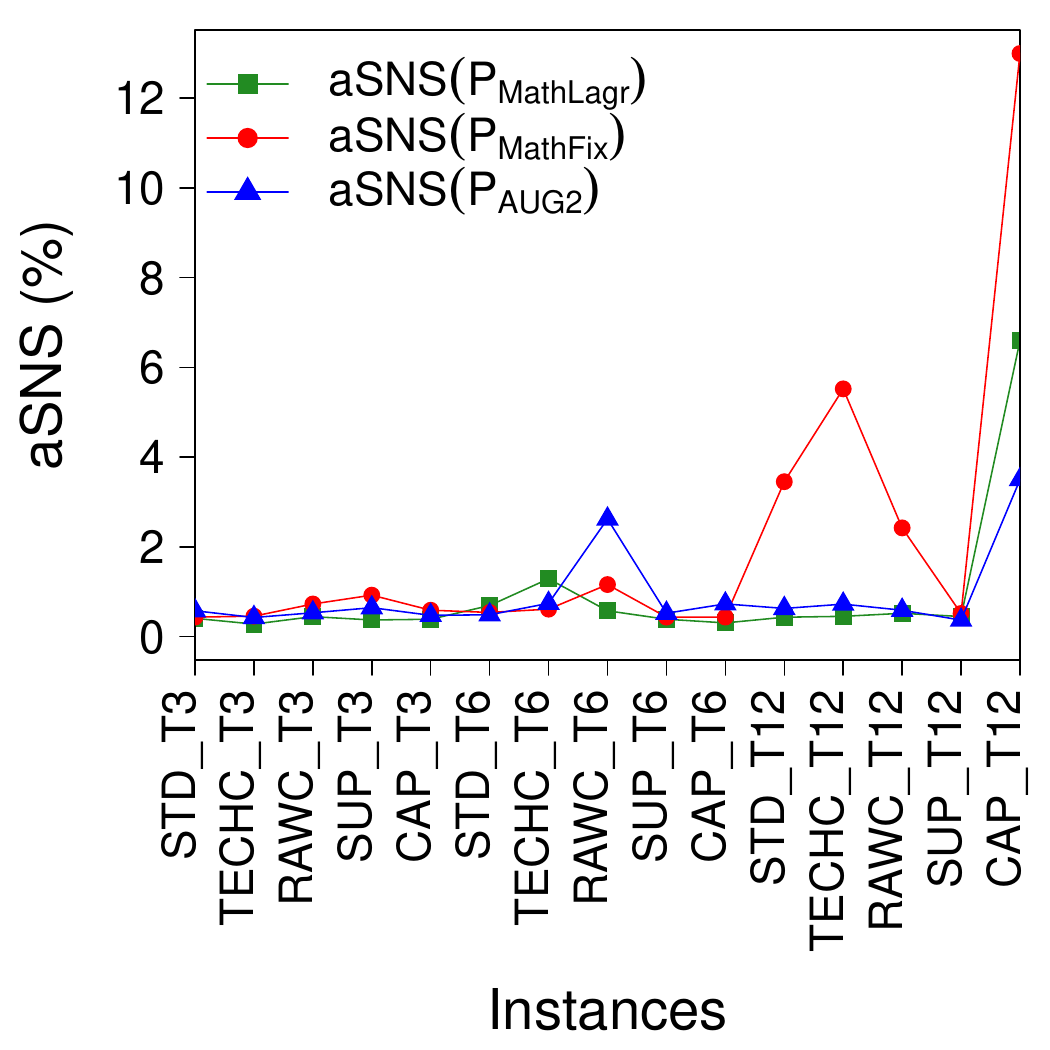}}
\subfigure[fig:snsI17M4][{$|I|=25$ and $|M|=4$.}]{ \includegraphics[width=0.234\textwidth]{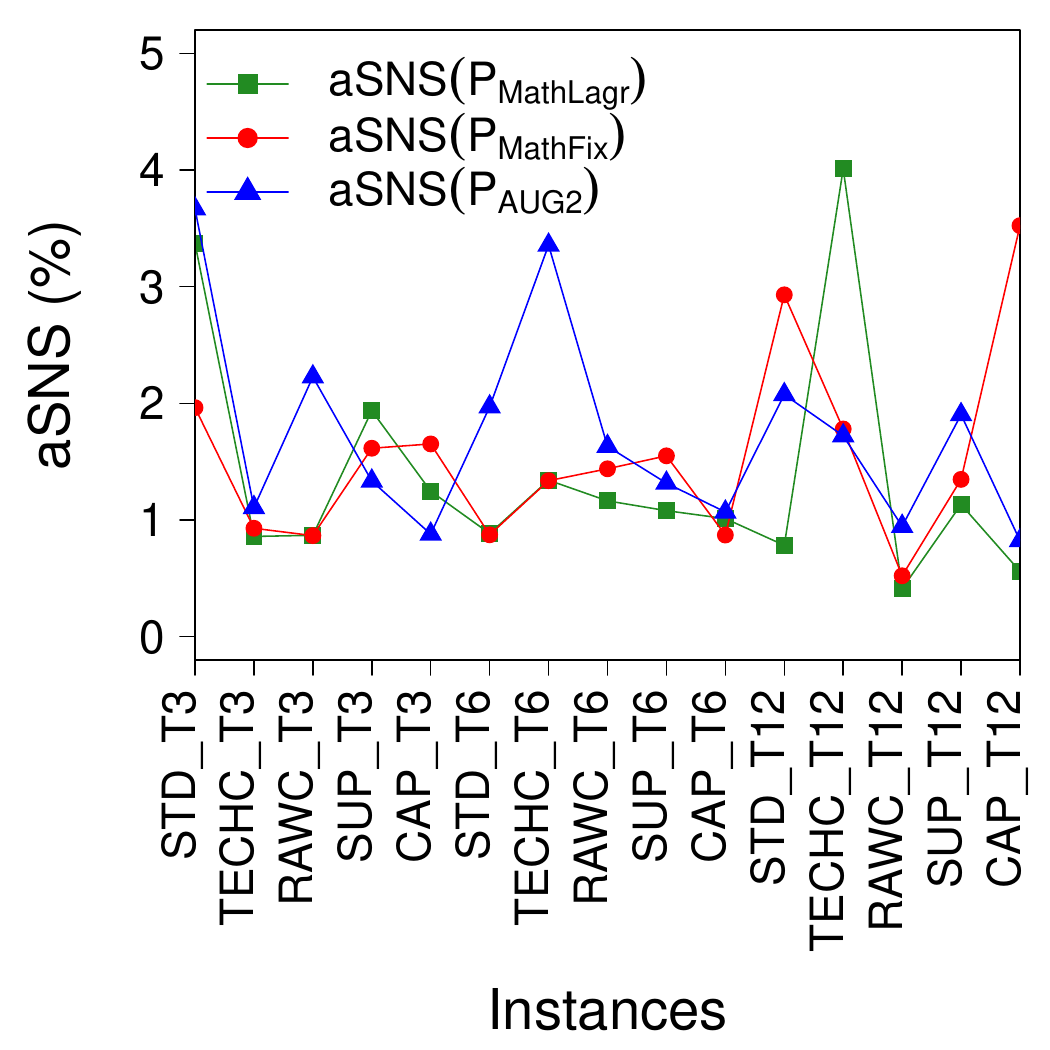}}
\subfigure[fig:snsI25M8][{$|I|=25$ and $|M|=8$.}]{ \includegraphics[width=0.234\textwidth]{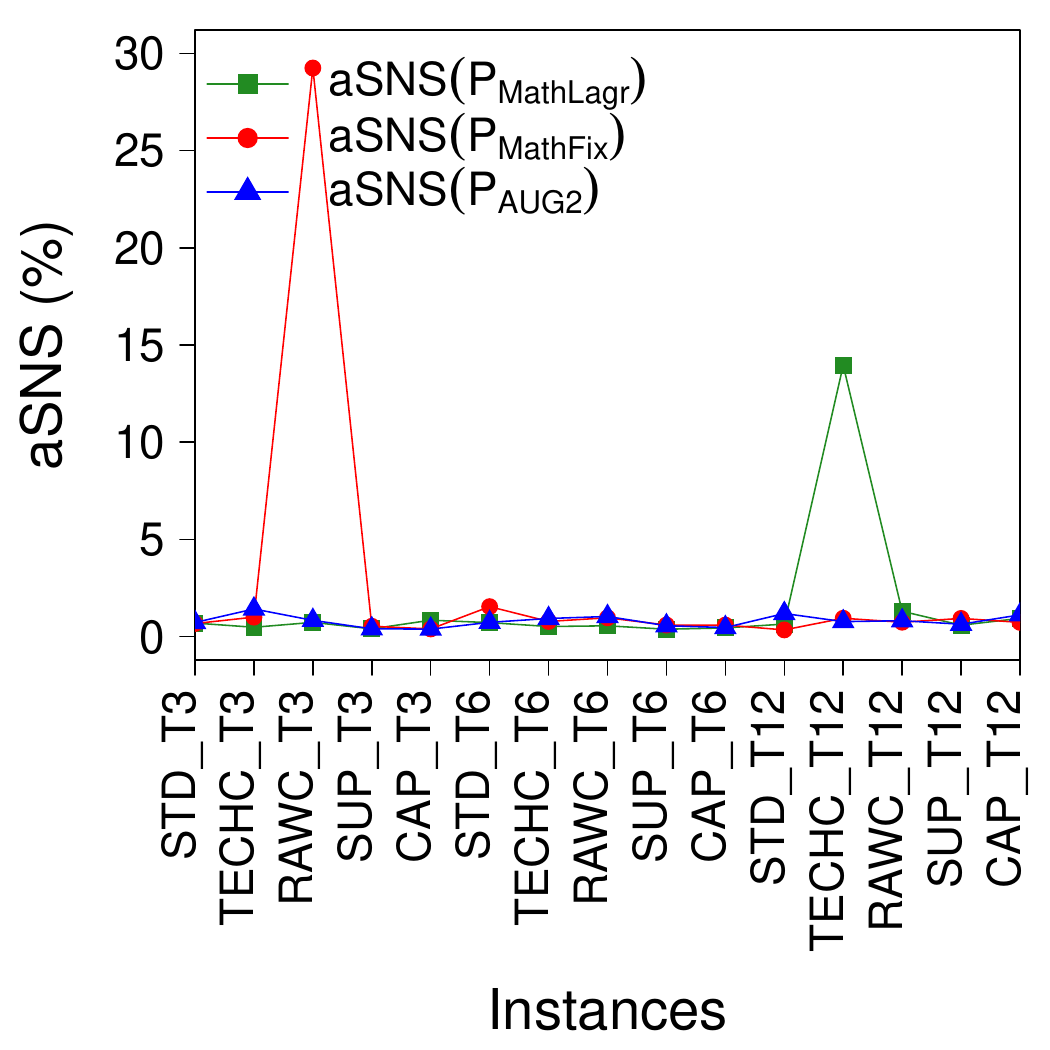}}
	\caption{\typo{Measure $aSNS$} for $P_{MathLagr}$, $P_{MathFix}$ and $P_{AUG2}$.}
      \label{fig:metricsM_aSNS}     
    \end{figure}
\begin{figure}[!htb]
 \subfigure[fig:R2I17M4][{$|I|=17$ and $|M|=4$.}]{ \includegraphics[width=0.234\textwidth]{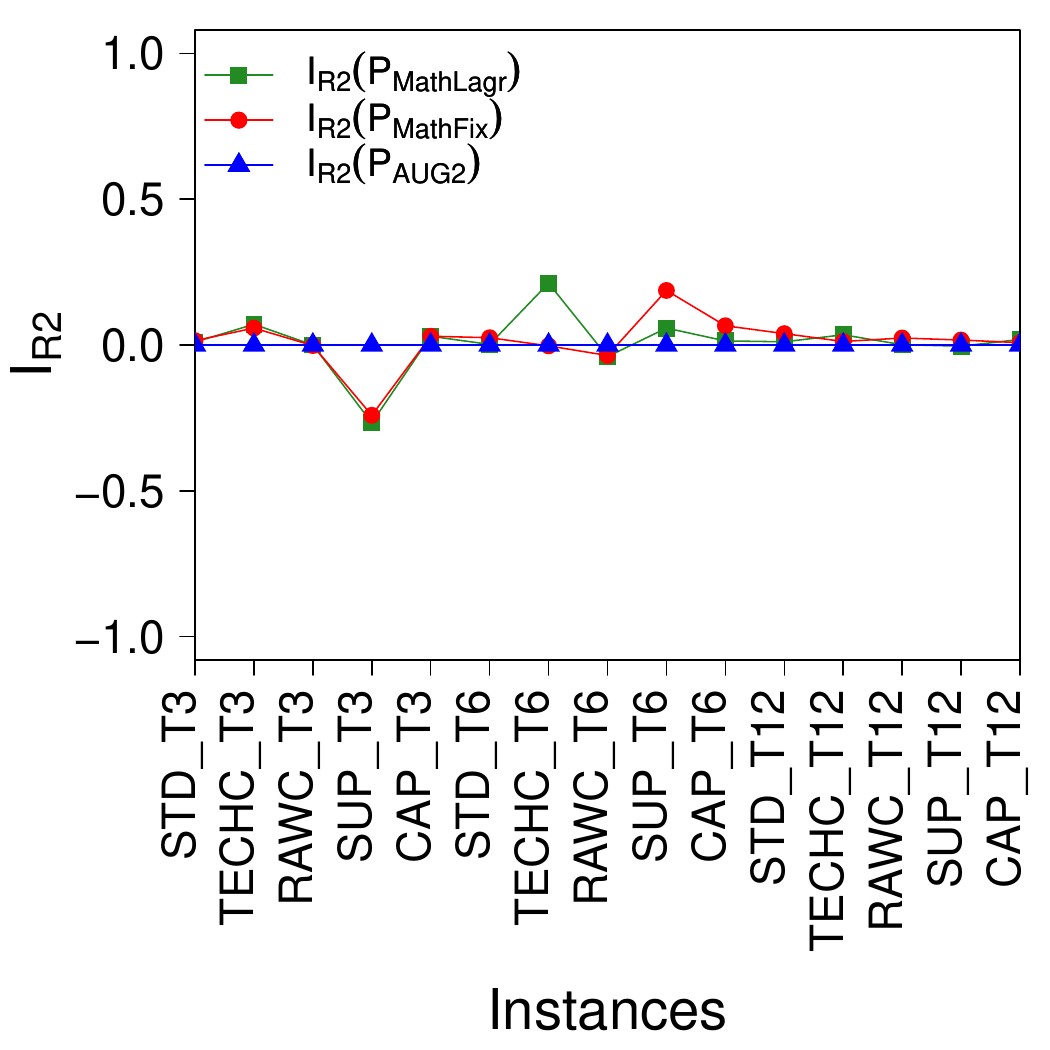}}
\subfigure[fig:R2I17M4][{$|I|=17$ and $|M|=4$.}]{ \includegraphics[width=0.234\textwidth]{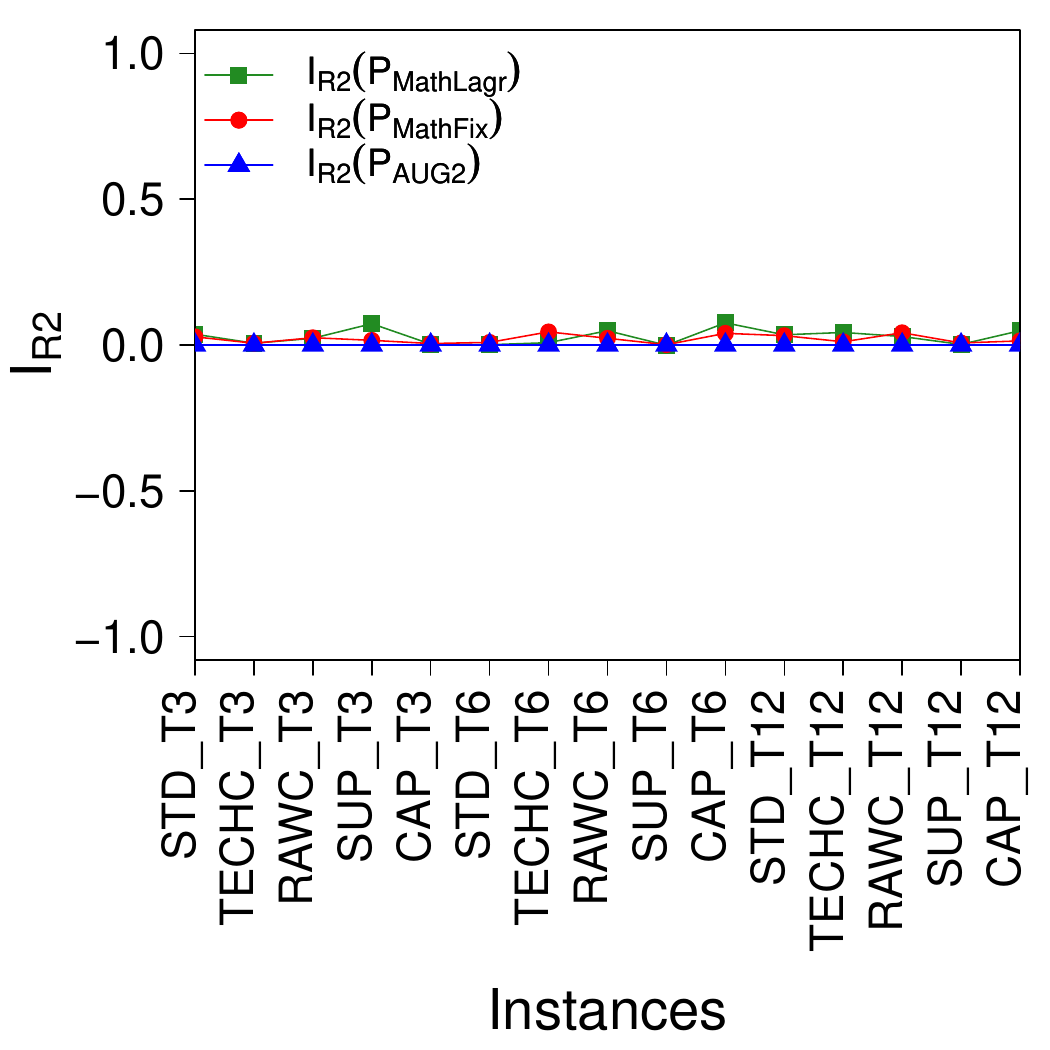}}
\subfigure[fig:R2I17M4][{$|I|=25$ and $|M|=4$.}]{ \includegraphics[width=0.234\textwidth]{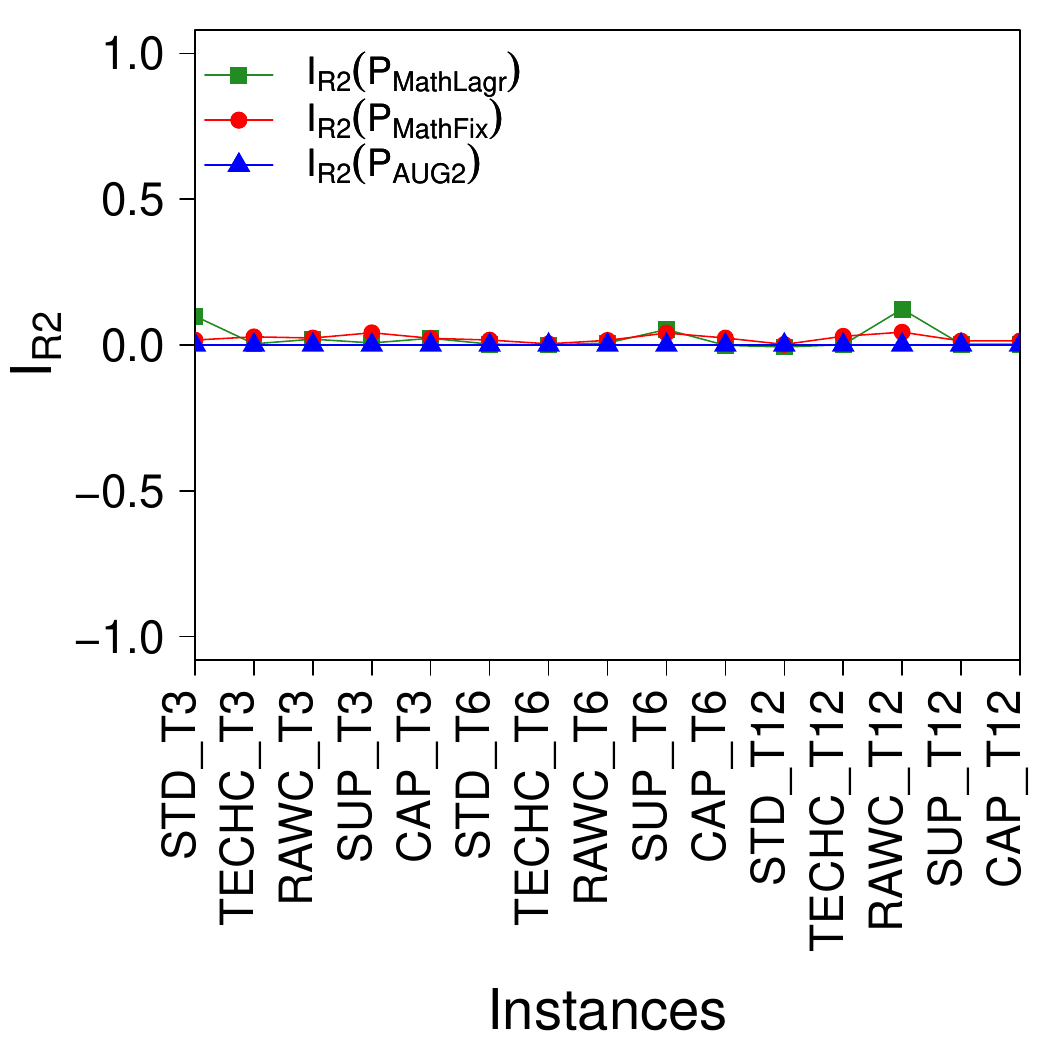}}
\subfigure[fig:R2I25M8][{$|I|=25$ and $|M|=8$.}]{ \includegraphics[width=0.234\textwidth]{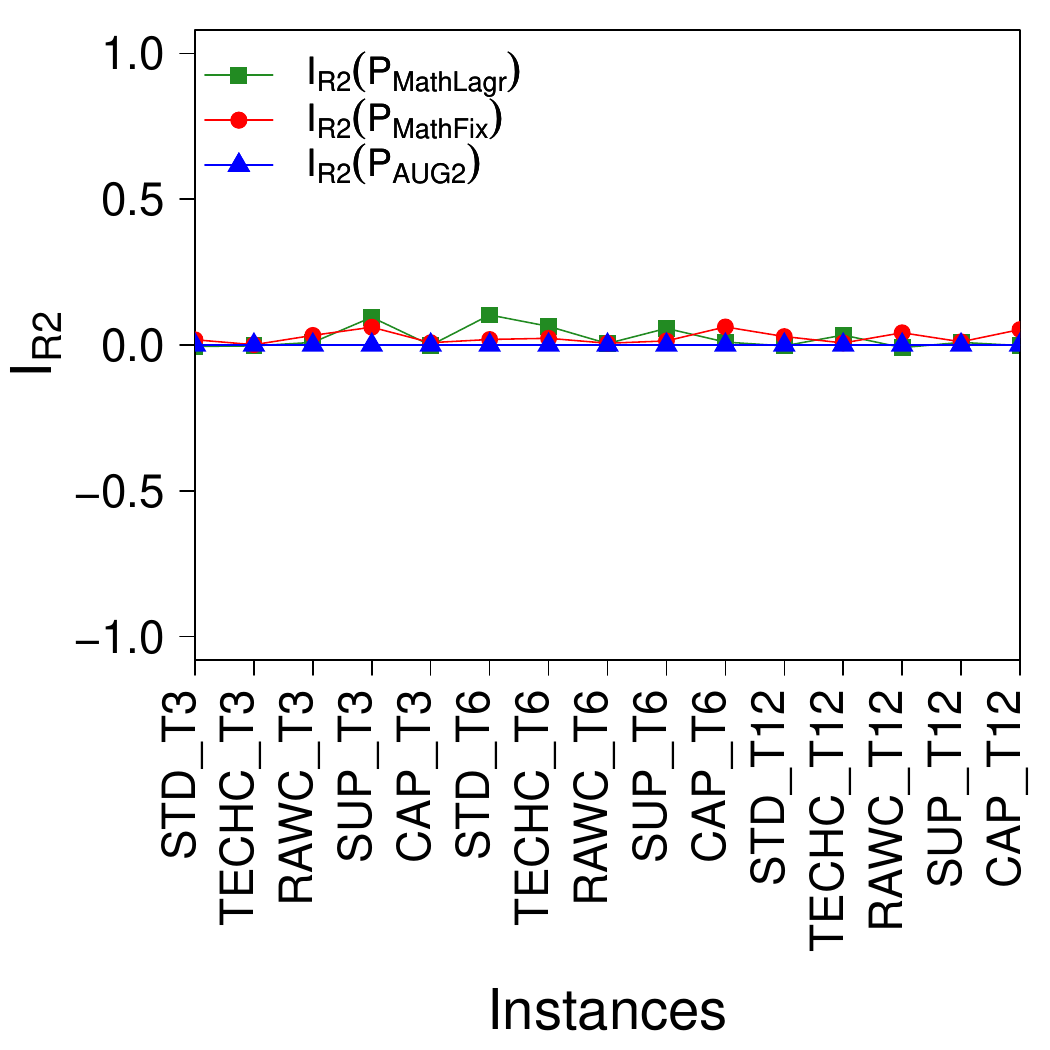}}
	\caption{{R2 Indicator \typo{for} $P_{MathLagr}$, $P_{MathFix}$ and $P_{AUG2}$.}}
      \label{fig:metricsM_R2}
\end{figure}
\begin{figure}[!htb]
\subfigure[fig:timeI17M4][{$|I|=17$ and $|M|=4$.}]{ \includegraphics[width=0.234\textwidth]{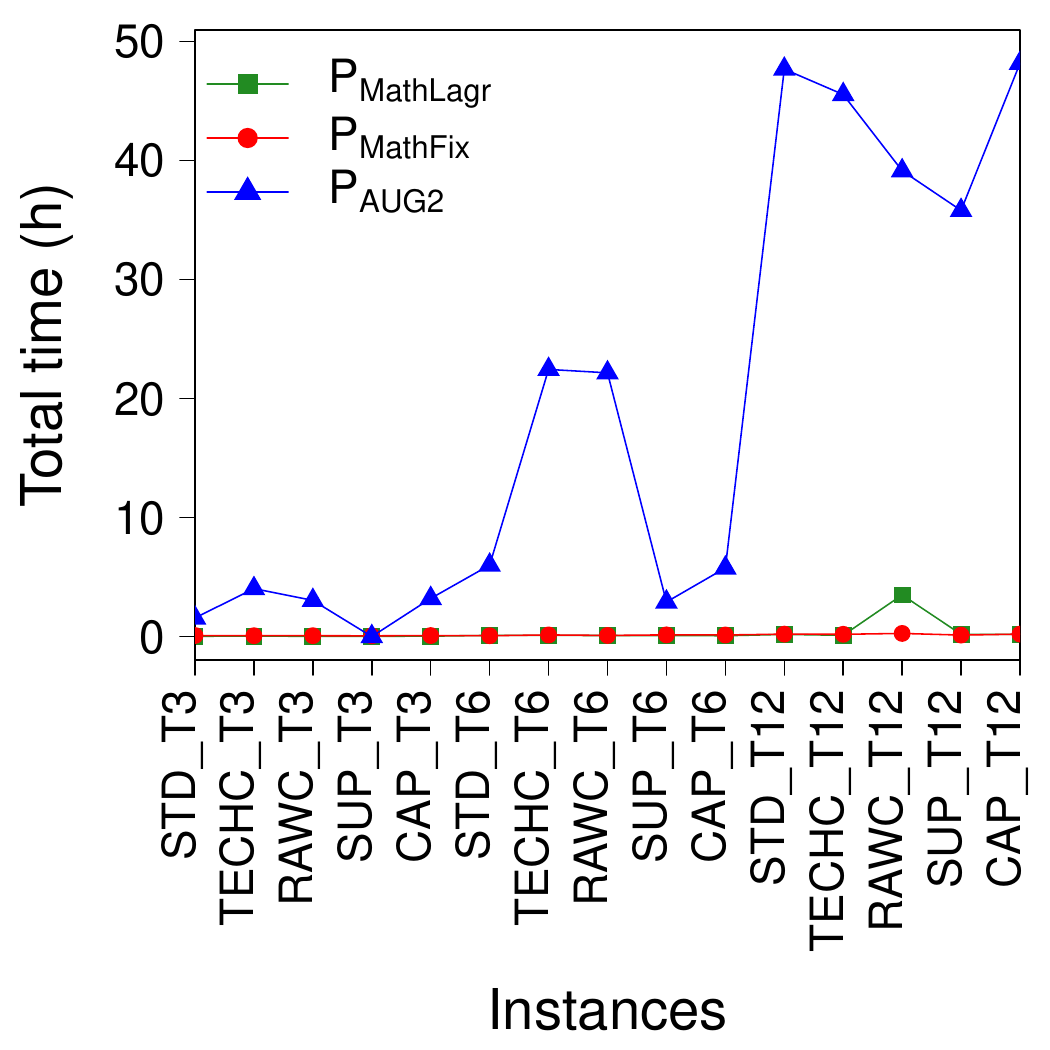}}
\subfigure[fig:timeI17M4][{$|I|=17$ and $|M|=8$.}]{ \includegraphics[width=0.234\textwidth]{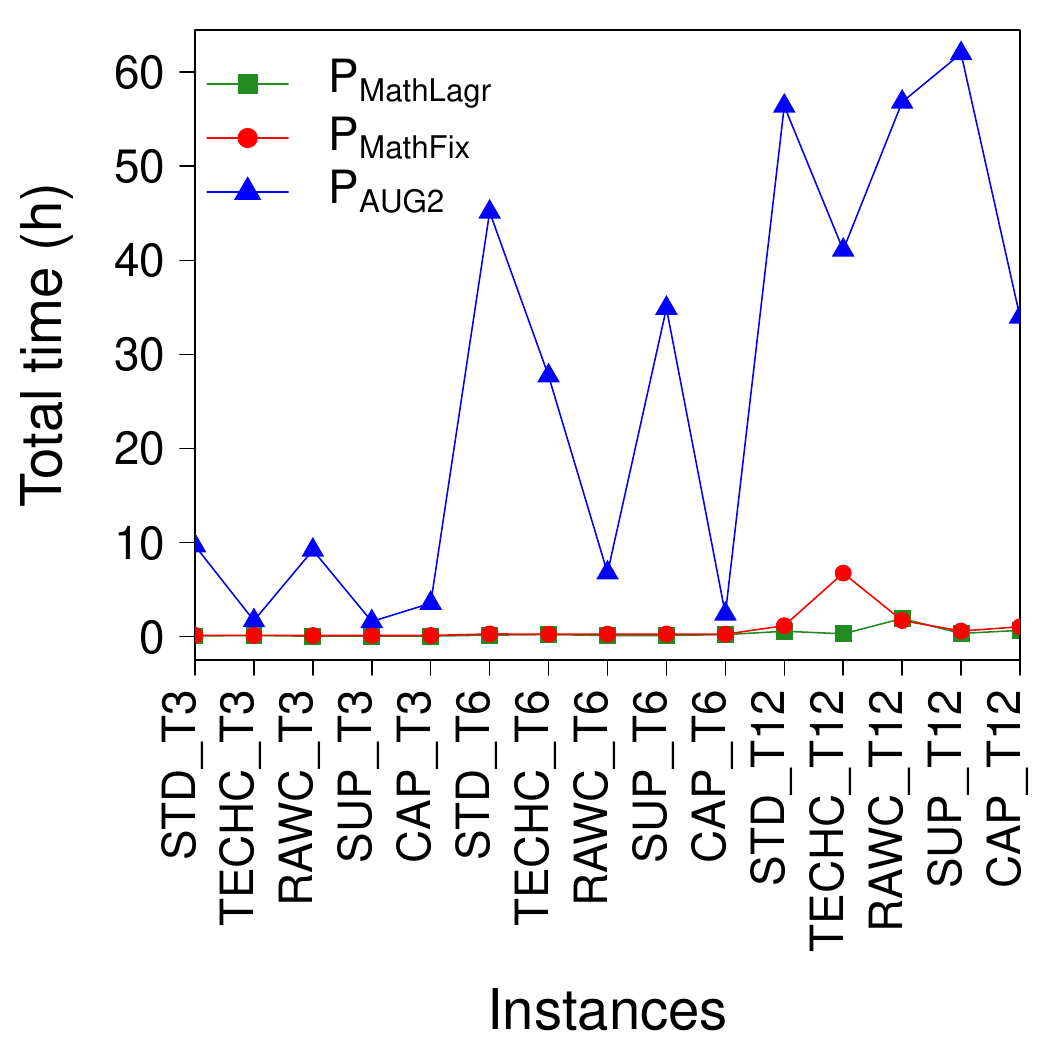}}
\subfigure[fig:timeI17M4][{$|I|=25$ and $|M|=4$.}]{ \includegraphics[width=0.234\textwidth]{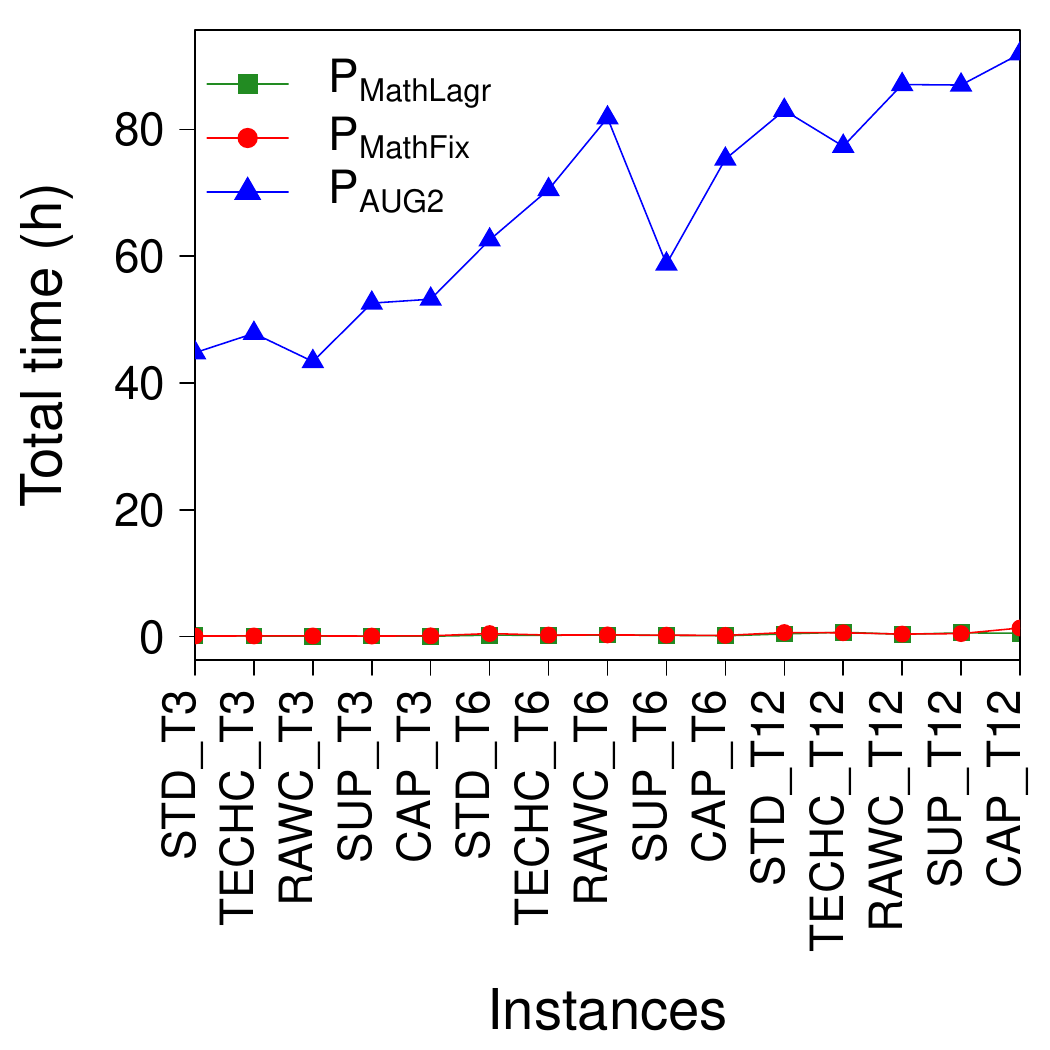}}
\subfigure[fig:timeI25M8][{$|I|=25$ and $|M|=8$.}]{ \includegraphics[width=0.234\textwidth]{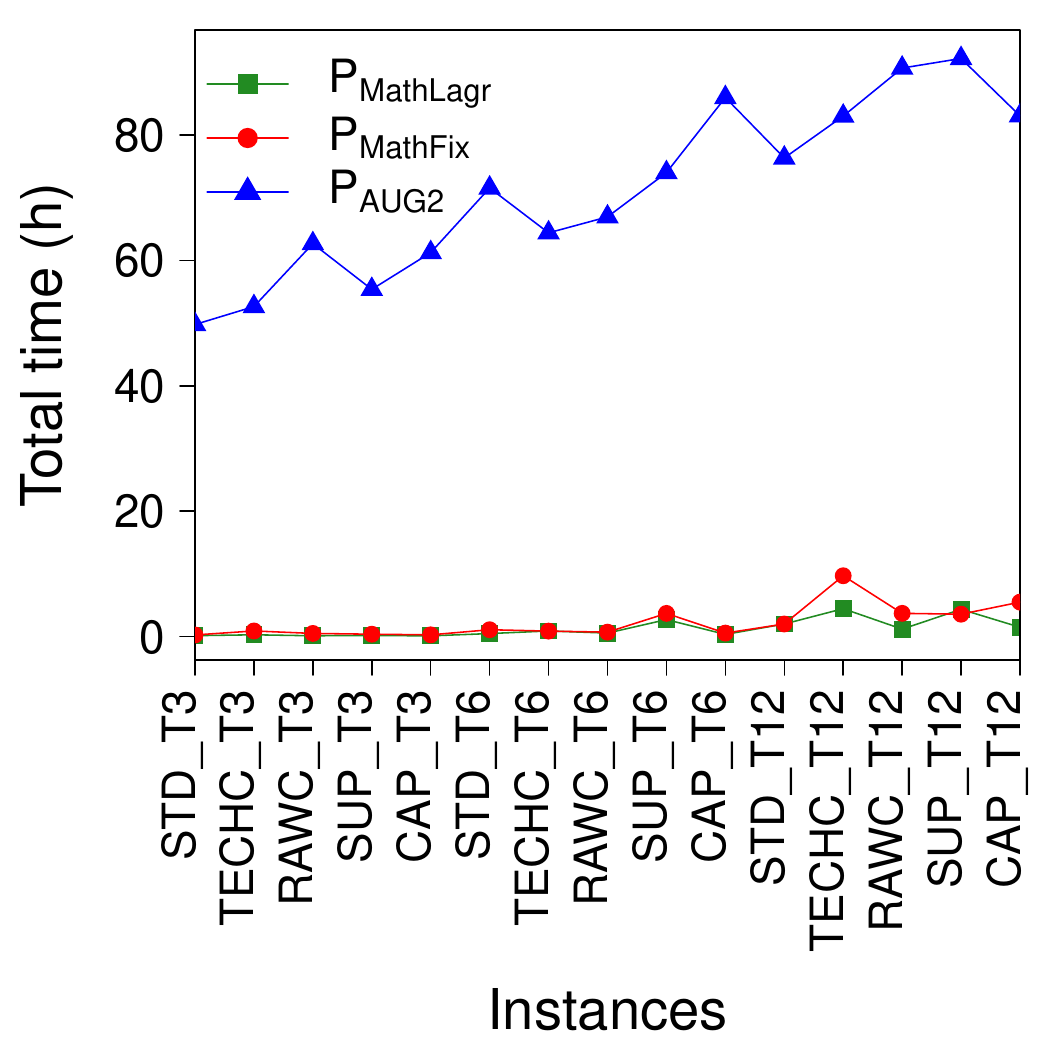}}
       \caption{{Total running time for obtaining  $P_{MathLagr}$, $P_{MathFix}$ and $P_{AUG2}$.}}
       \label{fig:metricsM_timeTot}
 \end{figure}

Let $P_{AUG2}$,  $P_{MathLagr}$ and $P_{MathFix}$ be, respectively, the Pareto frontier approximations obtained by AUGMECON2, \emph{AugMathLagr} and \emph{AugMathFix}.
Figures \ref{fig:metricsM_aMID}, \ref{fig:metricsM_aSNS} and \ref{fig:metricsM_R2} present, respectively, the $aMID$ measure, the $aSNS$ measure and the R2 Indicator for $P_{AUG2}$,  $P_{MathLagr}$ and $P_{MathFix}$. Figure \ref{fig:metricsM_timeTot} presents the total computational running times for obtaining $P_{AUG2}$, $P_{MathLagr}$ and $P_{MathFix}$.

As can be noted in Figure \ref{fig:metricsM_aMID}, the $aMID$ measures for the solutions of \emph{AugMathLagr} and \emph{AugMathFix} are lower than or approximately the same as the solutions achieved by AUGMECON2 for $42$ and $32$, respectively, out of the {$60$} instances. Note that these values correspond to \minor{$70\%$ and approximately $53\%$}
respectively of the instances considered. Although the $aMID$ values of $P_{MathLagr}$ are significantly high in 4 other instances, their overall results indicate that the solutions provided by \emph{AugMathLagr}  present low {$GAPM$}s in comparison to the Ideal point for most of the instances.

The standard deviations of $aMID$, measured by the $aSNS$ measure, for solutions obtained by
\emph{AugMathLagr} and \emph{AugMathFix} were considerably larger than those achieved by AUGMECON2, which were significantly high for 6 and 10 instances, respectively. AUGMECON2, in particular, presented unexpectedly large $aSNS$ values for 4 other instances.

The R2 indicator value for the solutions obtained by AUGMECON2 is null, since we {considered} $P_{AUG2}$ as the set of representative solutions for the measure.
{The maximum distance from a solution obtained by \emph{AugMathLagr} and \emph{AugMathFix} to the Ideal point is close to the maximum distance from a solution found by AUGMECON2 to the Ideal point.}

The total computational running times required {by} AUGMECON2 were at least \minor{$11$ and $6$}  times higher than {those} required by \emph{AugMathLagr} and \emph{AugMathFix} in instances with at least 6 periods in the planning horizon.
 Concerning the running times, the {advantages} of \emph{AugMathLagr} and \emph{AugMathFix} over AUGMECON2 are even more evident when considering the largest instances with 25 entities, i.e. $|I|=25$. \emph{AugMathLagr}, in particular, was from \minor{$18$ to $550$}  times faster than AUGMECON2 and from \minor{$0.8$  to $3.7$}
  times faster than \emph{AugMathFix} in instances with $25$ entities and at least $6$ periods in  the planning horizon.

 The results of this experiment attest \typo{to} the good quality of the Pareto frontier approximations obtained by \emph{AugMathLagr}. Moreover, the adaptation of \emph{AugMathFix} presented significantly inferior running times in comparison to AUGMECON2. \emph{AugMathLagr} was even faster than \emph{AugMathFix} in 81.667\% of the tested instances.
 This enables us to infer that the proposed Lagrangian matheuristic is suitable for solving large instances.

 \minor{In the supplementary report in \citep{Mendeley2020}, we present experiments with instances whose demands have increasing, decreasing and alternating trends along the time periods. The results of the experiments attest that the proposed method is robust in finding good quality solutions for these scenarios. Therefore, we can infer that the proposed method is suitable for dealing with supply chains whose demands vary in the planning horizon.}

\section{Final remarks} \label{sec:conclusion}

SSC management has recently drawn the attention of industrial and academic sectors to new research developments.
Most studies approach SSC management by multi-objective optimization, focusing on the construction of optimization models based on case studies.
This approach, however, is commonly very time-consuming and, even though it belongs to the strategic and tactical levels of planning, to refine, validate and thoroughly study the supply chain, it is necessary to employ efficient solution methods \citep{Eskandarpour2015}. 
Perhaps due to the difficulty in solving these problems and because the studies are usually problem-oriented, little effort has been devoted to keeping a benchmark data repository. 
 This paper introduces an artificial instance generator and an efficient Lagrangian matheuristic, here called \emph{AugMathLagr}, for SSC management multi-objective problems. In particular, \emph{AugMathLagr} is implemented to the optimization model introduced by \citet{Mota2018}.
In addition to introducing \emph{AugMathLagr}, this paper also efficiently adapts the matheuristic \emph{AugMathFix} proposed by \citet{Tautenhain2019} to the target problem.

Experiments conducted with a case study found in the literature and with a test bed of artificial instances \typo{indicated} that in comparison with a classical exact method for multi-objective problems, known as AUGMECON2, \emph{AugMathLagr} and \emph{AugMathFix} \minor{were} significantly faster. In particular, \emph{AugMathLagr} \typo{was} faster than \emph{AugMathFix} and \typo{performed} even better when approximating the Pareto frontier in cases where larger instances \typo{were} considered. We assessed the Pareto frontier approximations {using} three multi-objective measures to evaluate how close the solutions from the Pareto frontier approximation were from the Ideal point. These results indicated that the solutions of the Pareto frontier approximation achieved by \emph{AugMathLagr} and \emph{AugMathFix} \typo{were} close to the Ideal point, in comparison with the solutions from AUGMECON2. In particular, the values of the  R2 indicator expressed the high quality of the solutions found by \emph{AugMathLagr} and \emph{AugMathFix}.

\minor{In the case we studied in this paper, the solution methods identified the most representative decisions for environmental and social objectives. On the one hand, the social criterion benefits from more entities being opened. On the other, manufacturing and remanufacturing decisions represent most of the environmental impacts. We acknowledge, however, that this result might be biased towards the fact that our case study was based on an electronics component producer, for which manufacturing and remanufacturing activities are higher than for other industries.}

\minor{Nonetheless, the results related to the environmental impacts are strongly dependent on the scope of the investigation and on the case under study. Since data associated with environmental and social objectives are particularly challenging to be accessed, companies should carefully study the parameters associated with the most representative decision in order to better meet the reality. In particular, data available for similar supply chains should be employed as a starting point. Additionally, sensitivity analysis can also be conducted on these parameters to ensure that the obtained supply chain results and conclusions remain unaltered.}

This paper also contributes with a framework for generating artificial instances. Artificial instances are particularly useful for assessing the performance of optimization models since they allow statistical inferences on a wide variety of numerical scenarios. Due to the specificity of each optimization model, the implementation of the instance generator to other SSC management optimization models is a suggestion for future work.

Another direction for future work regards the proposed matheuristic. The primary goal of such study would be to further improve the quality of the solutions in order to obtain even closer values to the optimum ones achieved by AUGMECON2 using local search procedures.

\section*{Acknowledgments} \label{Acknowledgements}

Authors Camila P. S. Tautenhain and Mariá C. V. Nascimento would like to acknowledge {the funding granted by} São Paulo Research Foundation (FAPESP), grant numbers: 14/27334-9, 15/21660-4 and 16/02203-4; Conselho Nacional de Desenvolvimento Científico e Tecnológico (CNPq), grant numbers: 306036/2018-5; and Coordenação de Aperfeiçoamento de Pessoal de Nível Superior – Brazil (CAPES) – Finance Code 001. The authors Ana Paula Barbosa-Povoa and Bruna Mota would like to acknowledge the financial support from FCT and Portugal 2020 FCT under the project PTDC/EGEOGE/28071/2017, Lisboa -01.0145-Feder-28071. The author Mariá C.V. Nascimento is also grateful to Leonardo V. Rosset for giving her a hand.
Research carried out using the computational resources of the Center for Mathematical Sciences Applied to Industry (CeMEAI) funded by FAPESP (grant 2013/07375-0).


 \bibliographystyle{model5-names} 
 \small
	\bibliography{referencias}
 
 

\section*{Supplementary material}

# Supplementary material to the paper "An efficient Lagrangian-based heuristic to solve a multi-objective sustainable supply chain problem"

## 1 Model generalization for different networks

In this section, we describe how to generalize the studied model to address different SSC networks. Moreover, we provide some examples to illustrate how to deal with additional planning situations as, e.g., two levels of warehouses, assembled or process industries and sorting centers for the reverse logistics.

The allowed flow of items between entities in the studied SSC is defined by tuples $\{(m, i, j)$ such that $m \in M, i \in I, j \in J\}$. Table 1 presents the sets of tuples that indicate the in-flow and out-flow of items in each entity. In addition to defining the SC network, these tuples also specify the items that each entity can receive and send.

Let us take set $F_{OUTFFP}$ to exemplify how to add new features to the SSC network regarding product flows. This set represents the out-flow in factories which establishes the delivery of final products from factories: (i) directly to customers or (ii) first to warehouses, which can also act as distribution centers. The in-flow in customers, $F_{INCFP}$, must also be adjusted accordingly to ensure the proper flow through the network. Analogously, set $F_{OUTCRP}$ controls where to send each product recovered from customer: (i) directly to factories for recovery or recycling or (ii) to warehouses, acting as collection or sorting centers, to be later shipped to factories. The in- and out-flow of recovered products in factories and in warehouses must also be defined accordingly.

The formulation can be generalized in such a form to address more than one level of warehouses. For example, consider a network with two levels of warehouses. The first warehouse, $w_1 \in I_w$, stores a given final product, $fp_1 \in M_{fp}$, to be later shipped to a second warehouse, $w_2 \in I_w$, which acts as a distribution center. In this example, sets $F_{OUTW}$ and $F_{INW}$ must both include tuple $(fp_1, w_1, w_2)$ to represent the out-flow and in-flow of product $fp_1$, respectively.

Assembled goods that require manufacturing intermediate products can be modeled by including the intermediate products in the bill of products required to manufacture the final product. For this, let $BOM_{mng}^{prod}$ be the bill of raw materials or final products $m \in M_{rm} \cup M_{fp}$ to produce another final product $n \in M_{fp}$ using technology $g \in G_{prod}$. Additionally, assembled or process industries can be defined by setting the final products that each factory can receive and send to other factories in the corresponding in- and out-flow sets, $F_{INFFP}$ and $F_{OUTFFP}$, respectively.

In chemical industries, SCs transform raw materials into a set of products using chemical processes. These industries often require different sets of raw materials and the processing times are long. But although having special characteristics they can be handled by the current approach.

Table 1: Description of the in- and out-flow of items in entities in the studied SC.

| Definition | Set | Description |
|---|---|---|
| In-flow of final products ($fp$) in factories | $F_{INFFP}$ | Factories can receive $fp \in M_{fp}$ from:<br>– other factories;<br>– warehouses (in case of assembled goods, for example);<br>– airports and airports through cross-docking. |
| Out-flow of final products ($fp \in M_{fp}$) in factories | $F_{OUTFFP}$ | Factories can send $fp$ to:<br>– factories;<br>– warehouses;<br>– directly to customers;<br>– seaports and airports through cross-docking. |
| In-flow of recovered products ($rp$) in factories | $F_{INFRP}$ | Factories can receive $rp$ from:<br>– warehouses;<br>– directly from customers;<br>– seaports or airports through cross-docking. |
| In-flow in warehouses | $F_{INW}$ | - Warehouses can receive:<br>– $fp$ from factories or other warehouses;<br>– $rp$ from customers or other warehouses;<br>– $fp$ or $rp$ to airports or seaports through cross-docking. |
| Out-flow in warehouses | $F_{OUTW}$ | - Warehouses can send:<br>– $fp$ to factories, other warehouses or customers;<br>– $rp$ to factories or other warehouses;<br>– $fp$ or $rp$ to airports or seaports through cross-docking. |
| In-flow in airports and seaports | $F_{INAIR}$, $F_{INPORT}$ | Airports and seaports can receive:<br>– $rm$ from suppliers;<br>– $fp$ from factories or warehouses;<br>– $rp$ from warehouses or customers;<br>– $rm$, $fp$ or $rp$ from other airports or seaports through cross-docking. |
| Out-flow in airports and seaports | $F_{OUTAIR}$, $F_{OUTPORT}$ | Airports and seaports can send:<br>– $rm$ to factories;<br>– $fp$ to factories, warehouses or customers;<br>– $rp$ to factories or warehouses;<br>– $rm$, $fp$ or $rp$ to other airports or seaports through cross-docking. |
| In-flow of final products ($fp$) in customers | $F_{INCFP}$ | Customers can receive $fp$ from:<br>–factories or warehouses. |
| Out-flow of recovered products ($rp$) in customers | $F_{OUTCRP}$ | Customers can send $rp$ to:<br>– factories;<br>– warehouses (warehouses can act as sorting centers);<br>– airports and airports through cross-docking. |

## 2 Adapted *AugMathFix*

This appendix shows the matheuristic *AugMathFix* originally introduced in [2], and which we adapted in this paper to solve the model introduced in [1]. The pseudocode of the adapted method is presented in Algorithm 1.

First, let us present the problem with $\epsilon$-Constraints to be solved by the adapted *AugMathFix* – referred to as MOP in the paper.

$$(\text{MOP}): \min \quad -f_{eco}(u, v, w', w'') - eps\left(\frac{l_{env}}{r_{env}} + 0.1 \frac{l_{soc}}{r_{soc}}\right) \tag{1}$$

$$s.t. \quad \mathcal{A}' u + \mathcal{A}'' v \leq \beta' + \beta'' w'' \tag{2}$$

$$E' u + E'' v \leq d \tag{3}$$

$$f_{env}(u, v, w', w'') + l_{env} = \epsilon_{env} \tag{4}$$

$$-f_{soc}(u, v, w', w'') + l_{soc} = \epsilon_{soc} \tag{5}$$

$$u \in \mathbb{R}_{\geq 0}, v \in \mathbb{Z}_{\geq 0}^{n''}, w' \in \{0,1\}^{b'}, w'' \in \{0,1\}^{b''} \tag{6}$$

$$l_{env}, l_{soc} \in \mathbb{R}_{\geq 0} \tag{7}$$

The adapted *AugMathFix* is an iterative AUGMECON2-based matheuristic that solves a sequence of

mono-objective problems MOP by a heuristic method known as *MathFix*. These problems are generated by varying the values of $\epsilon$ as in *AugMathLagr*. Therefore, the pseudocode of *AugMathFix* – Algorithm 1 – is the same as *AugMathLagr* except for line 10, which corresponds to the solution method to solve the MOP. Moreover, for ease of notation, we shall refer to the solutions, that is, the values of the decision variables $u, v, w'$ and $w''$, as $x_{FS}^{(it)}, x_{RL}^{(it)}, x^{best}$ or $x_f$.

---

**Algorithm 1:** Adapted *AugMathFix*.

**Input** : maximum number of iterations, $kMax$, number of grid points $dg$
**Output:** Pareto frontier approximation $\mathcal{P}$

1. Estimate the lower bounds $fL_i$ and the upper bounds $fU_i$, $i \in \{eco, env, soc\}$ for the objective functions $-f_{eco}, f_{env}$ and $-f_{soc}$
2. Identify ranges $r_j = fU_j - fL_j, j \in \{env, soc\}$ of the environmental and social objective functions, respectively
3. $step_j = \frac{r_j}{dg}, j \in \{env, soc\}$
4. $\mathcal{P} := \emptyset$
5. $gr_{env} = 0$
6. **for** $\epsilon_{env} = fU_{env}$ until $gr_{env} < dg$ **do**
7. $\quad gr_{soc} = 0$
8. $\quad$ **for** $\epsilon_{soc} = fU_{soc}$ until $gr_{soc} < dg$ **do**
9. $\quad\quad MOP :=$ Problem (1)-(7)
10. $\quad\quad x_f :=$ solve MOP using Algorithm 2 – *MathFix* – considering $\epsilon_{env}$ and $\epsilon_{soc}$
11. $\quad\quad$ **if** $x_f$ *is feasible* **then**
12. $\quad\quad\quad \mathcal{P} := \mathcal{P} \cup \{x_f\}$
13. $\quad\quad\quad l_{soc}, l_{env} :=$ values for the slack variables in $x_f$
14. $\quad\quad\quad gr_{soc} = gr_{soc} + 1 + \lfloor \frac{l_{soc}}{step_{soc}} \rfloor$
15. $\quad\quad\quad \epsilon_{soc} = \epsilon_{soc} - step_{soc} \left(1 + \lfloor \frac{l_{soc}}{step_{soc}} \rfloor\right)$
16. $\quad\quad$ **else** $gr_{soc} = dg$ ;
17. $\quad$ **end**
18. $\quad gr_{env} = gr_{env} + 1$
19. $\quad \epsilon_{env} = \epsilon_{env} - step_{env}$
20. **end**
21. $\mathcal{P} :=$ Remove dominated solutions from $\mathcal{P}$
22. **return** $\mathcal{P}$

---

*MathFix* is a heuristic to solve mono-objective problems based on the fixing and re-solving strategy. Algorithm 2 shows a pseudocode of this heuristic.

The inputs of Algorithm 2 are the original MOP, referred to as $P^{(0)}$, the maximum number of iterations of the algorithm, $IT$, and positive real scalars $\alpha$, $\beta$, $\theta$, $dec\alpha$, $dec\beta$ and $dec\theta$ such that $\alpha + \beta + \theta = 1$ and $dec\alpha + dec\beta + dec\theta = 0$. More details about these scalars can be found in the paper that first proposed the method. The values employed are the ones suggested in [2]. Therefore, the following parameter values were employed in the experiments described in this paper: $IT = 5, \alpha = 0.5, \beta = 0.5, \theta = 0, dec\alpha = -0.25, dec\beta = 0$ and $dec\theta = -0.25$.

In the first line of Algorithm 2, a relax-and-fix heuristic in *MathFix*, called *Heuristic*, is used to solve the initial problem $P^{(0)}$. The pseudocode of the *Heuristic* is presented in Algorithm 3 and shall be discussed later. Alike the original *Mathfix*, a control variable $ch$ dictates the choice for defining the closing of entities. In this adaptation, we consider two alternatives: if the choice is 0, i.e., if value 0 is assigned to $ch$ the option is to close a warehouse; if $ch$ is 1, the option is to close a factory. Moreover, in line 2, auxiliary variables of the algorithm related to this choice are initialized: $num_o = 0$ (number of times that choice $o$ was applied in the heuristic) and $max_o = IT$ (maximum number of times that choice $o$ can be applied in the heuristic before restarting the problem), for each option $o \in \{0, 1\}$.

**Algorithm 2:** *MathFix*.

**Input** : Problem $P^{(0)}$, maximum number of iterations, $IT$ and positive scalars $\alpha$, $\beta$, $\theta$, $dec\alpha$, $dec\beta$ and $dec\theta$ such that $\alpha + \beta + \theta = 1$ and $dec\alpha + dec\beta + dec\theta = 0$

**Output:** The best feasible solution $x^{best}$

1 $x_{FS}^{(0)} := Heuristic(P^{(0)})$
2 $num_o = 0, max_o = IT, o \in \{0, 1\}$
3 $x^{best} := x_{FS}^{(0)}$
4 $it = 1, ch = 0$
5 **while** $it < IT$ **do**
6     $vuse_i = \sum_{t \in T} \sum_{m \in M_{fp} \cup M_{rp}} S_{mit}, \forall i \in I_f \cup I_w$
7     $veco_i = wpsq_i YC_i + lc_i\ wwh\ wpt\ wpsq_i\ YC_i + lc_i\ wi_i\ wwh\ wpt, \forall i \in I_f \cup I_w$
8     $venv_i = \arg\min_{i \in I_w} \sum_{c \in C}(\eta_c ei_c YC_i), \forall i \in I_f \cup I_w$
9     **if** $ch = 0$ **then** $\mathcal{I} := \arg\min_{i \in I_w}(\alpha\ vuse_i + \beta\ veco_i + \theta\ venv_i)$ ;
10     **else if** $ch = 1$ **then** $\mathcal{I} := \arg\min_{i \in I_f}(\alpha\ vuse_i + \beta\ veco_i + \theta\ venv_i)$ ;
11     $P^{(it)} :=$ Fix non-binary variables index by $\mathcal{I}$ at problem $P^{(it-1)}$
12     $x_{FS}^{(it)} := Heuristic(P^{(it)})$
13     **if** $x_{FS}^{(it)}$ *is feasible* **then**
14         $num_{ch} = num_{ch} + 1$
15         **if** $-f_{eco}(x_{FS}^{(it)}) > -f_{eco}(x^{best})$ **then** $x^{best} := x_{FS}^{(it)}$ ;
16         **if** $ch = 0$ **then** $ch = 1$;
17         **else if** $ch = 1$ **then** $ch = 0$;
18         **if** $num_{ch} \geq max_{ch}$ **then**
19             $P^{(it)} := P^{(0)}, x_{FS}^{(it)} := x_{FS}^{(0)}$
20             $max_{ch} = num_{ch}, num_o = 0, o \in \{0, 1\}$
21         **end**
22     **else**
23         $\alpha = \alpha + dec\alpha, \beta = \beta + dec\beta, \theta = \theta + dec\theta$
24     **end**
25     $it = it + 1$
26 **end**
27 **return** $x^{best}$

At each iteration of *MathFix*, $veco_i$, $venv_i$ and $vuse_i$ are defined in lines 6, 7 and 8 as, respectively, the economic cost, the environmental impact and the storage or production related to warehouse or factory $i$. In line 7, $wwh$, $wpt$ and $wpsq_i$ are, respectively, the number of weekly working hours, the number of weeks per time period and the minimum number of workers per square meter in entity $i \in I_f \cup I_w$. In line 9, if $ch = 1$, $\mathcal{I}$ is defined as the warehouse $i$ whose sum $\alpha veco_i + \beta venv_i + \theta vuse_i$ is the minimum amongst all $i \in I_w$. If $ch = 1$, in line 10, $\mathcal{I}$ is defined as the factory $i$ whose value of $\alpha veco_i + \beta venv_i + \theta vuse_i$ is the minimum amongst all $i \in I_f$. In line 11, $P^{(it)}$ is obtained by fixing $Y_i = 0, YC_i = 0, KT_{ait} = 0, K_{ai} = 0$ and $S_{mit} = 0$ and removing the constraints $S_{mit} \geq ic_{mi}^{min}Y_i, m \in M_{fp}, a \in A, t \in t$ and $YC_i \geq ea_i^{min}Y_i$ at $P^{(it-1)}$. Moreover, if $ch = 1$, i.e., $i$ is a factory, $P_{mgit} = 0, (m,g) \in H_{prod}$, and $R_{mgit} = 0, (m,i) \in H_{rem}$ and the constraints $P_{mgit} \geq pc_g^{min}Z_{gmi}, (m,g) \in H_{prod}$ and $R_{mgit} \geq pc_g^{min}Z_{gmi}, (m,g) \in H_{rem}, t \in T$ are also removed. In line 12, *Heuristic* is applied to $P^{(it)}$ and its solution is called $x_{FS}^{(it)}$. If solution $x_{FS}^{(it)}$ is feasible, $num_{ch}$ is incremented in line 14. The best solution achieved up to the current iteration $x^{best}$ is checked for update in line 15. In lines 16 and 17, $ch$ is updated to either 1 or 0, depending on its current value. In case $num_{ch} \geq max_{ch}$, the problem and feasible solution of the current iteration are reset to the respective problem and feasible solution of the first iteration. Moreover, $max_{ch}$ is limited by $num_{ch}$. In case $x_{FS}^{(it)}$ is infeasible, $\alpha$, $\beta$ and $\theta$ are updated in line 23. At the end, Algorithm 2 returns the best solution found over the iterations, $x^{best}$.

Algorithm 3 presents the heuristic to solve an input problem $P^{(it)}$.

4

**Algorithm 3:** *Heuristic*

    **Input**  : Problem $P^{(it)}$
    **Output:** Feasible solution $x_{FS}^{(it)}$
1   $P_{RL}^{(it)} :=$ Linear relaxation of problem $P^{(it)}$
2   $x_{RL}^{(it)} :=$ Solve $PR_{RL}^{(it)}$
3   $P_{FS}^{(it)} :=$ Modify problem $P^{(it)}$ by assigning 0 to real variables $X_{maijt}, m \in M, a \in A, i, j \in I, t \in T$, with value 0 at $x_{RL}^{(it)}$
4   $x_{FS}^{(it)} :=$ Solve $P_{FS}^{(it)}$
5   **return** $x_{FS}^{(it)}$

In line 1 of Algorithm 3, the linear relaxation of problem $P^{(it)}$ is assigned to $PR^{(it)}$. In line 2, the problem is solved using the CPLEX solver. Then, the problem given as input is assigned to $P_{FS}^{(it)}$ and in the next lines of the algorithm some variables of $P_{FS}^{(it)}$ are fixed according to the values of the variables of the solution of the linear relaxation. In line 3, problem $P_{FS}^{(k)}$ is created from $P_{RL}^{(it)}$ by fixing at zero the real variables $X_{maijt}, m \in M, a \in A, i, j \in I, t \in T$, which are null in the relaxed solution of the iteration $it$, i.e $x_{RL}^{(it)}$. In lines 4 and 5, CPLEX returns $x_{FS}^{(it)}$ as a solution for $P_{FS}^{(it)}$, which is feasible to the original problem.

## 3   Instance generator

This section describes at length the introduced instance generator. The generator estimates the maximum capacities according to the demand of the SSC customers. The minimum capacities are always estimated as a fraction of the maximum capacities.

In this section, the notation $Xv \sim U[lb, ub]$ indicates that a random variable labeled $Xv$ is randomly drawn from a continuous uniform interval ranging from $lb$ to $ub$, $lb, ub \in \mathbb{R}_{\geq 0}$.

The following sections present the guidelines for estimating the data parameters of the uniform random variables.

### 3.1   Demands and specification of items

Table 2 reports the parameters related to demand, bills and items. The demands for products are defined as random variables drawn from uniform distributions in the interval $[lbdc, ubdc]$. Suggestions to define all the input values required by the instance generator are available in [3].

Let $lbBOM^{prod}, ubBOM^{prod}, lbBOM^{rem}$ and $ubBOM^{rem}$ be scalars such that $lbBOM^{prod} < ubBOM^{prod}$ and $lbBOM^{rem} < ubBOM^{rem}$. The bill of raw materials $m \in M_{rm}$ required to manufacture product $n \in M_{fp}$ using technology $g \in G_{prod}$ is chosen with constant probability in the interval $[lbBOM^{prod}, ubBOM^{prod}]$. The bill of recovery products $m \in M_{rp}$ to remanufacture a final product $n \in M_{fp}$ follows a uniform distribution which falls within the interval $[lbBOM^{rem}, ubBOM^{rem}]$. The weight ($pw_m$) and necessary area per unit ($apu_m$) of raw material are drawn from uniform distributions in the intervals $[lbpw, ubpw]$ and $[lbapu, ubapu]$, respectively. The values of $pw_m$ and $apu_m$ of a final product $m$ are the average values of, respectively, $pw_n$ and $apu_n$ of raw materials $n$ required to manufacture $m$. The remanufacturing bill must require several recovery products to refurbish a new final product. Therefore, the values of $pw_m$ and $apu_m$ for a recovered product $m$ are the sums of, respectively, $pw_n$ and $apu_n$ of final products $n$ in the bill to remanufacture $m$.

The return rate depends on the type of industry and is a value randomly chosen between 0 and an input upper bound *ubret*. In the case study of the electrical industry, for example, the return rate is approximately 15%. The number of workers required to manage a technology $g$ is defined as the operating cost of $g$ ($opc_g$) multiplied by a pre-defined scalar $fracwg$.

Table 2: Instance generator description: customer demands, bills and item specification.

| Parameter | Value | Variables |
|---|---|---|
| $dmd_{mit}$ | $\sim U[lbdc, ubdc]$ | $lbdc, ubdc \in \mathbb{R}_{\geq 0}$ |
| $BOM_{mng}^{prod}$ | $\sim U[lbBOM^{prod}, ubBOM^{prod}]$ | $lbBOM^{prod}, ubBOM^{prod} \in \mathbb{R}_{\geq 0}$ |
| $BOM_{mn}^{rem}$ | $\sim U[lbBOM^{rem}, ubBOM^{rem}]$ | $lbBOM^{rem}, ubBOM^{rem} \in \mathbb{R}_{\geq 0}$ |
| $pw_m$ | $\begin{cases} \sim U[lbpw, ubpw], & \text{if } m \in M_{rm} \\ \frac{\sum_{(n,g) \in H_{prod}} BOM_{nmg}^{prod} pw_n}{|G_{prod}|}, & \text{if } m \in M_{fp} \\ \sum_{n \in M_{fp}} BOM_{nm}^{rem} pw_n, & \text{if } m \in M_{rp} \end{cases}$ | $lbpw, ubpw \in \mathbb{R}_{\geq 0}$ |
| $apu_m$ | $\begin{cases} \sim U[lbapu, ubapu], & \text{if } m \in M_{rm} \\ \frac{\sum_{(n,g) \in H_{prod}} BOM_{nmg}^{prod} apu_n}{|G_{prod}|}, & \text{if } m \in M_{fp} \\ \sum_{n \in M_{fp}} BOM_{nm}^{rem} apu_n, & \text{if } m \in M_{rp} \end{cases}$ | $lbapu, ubapu \in \mathbb{R}_{\geq 0}$ |
| $retF_m$ | $\sim [0, ubret]$ | $ubret \in \mathbb{R}_{\geq 0}$ |
| $wg_g$ | $fracwg \lceil opc_g \rceil$ | $fracwg \in \mathbb{R}_{\geq 0}$ |

## 3.2 Maximum and minimum capacities of entities

In an attempt to avoid an empty solution space when generating artificial instances and when estimating the installation area of entities, customer demands are taken into account. In line with this, Table 3 presents the parameters related to the capacity of the entities.

The higher the demands, the larger the inventory of products, production, acquisition of raw materials, and return of used products in the reverse flow must be. The storage capacity for a product $m$ in a given entity $i \in I$ ($ic_{mi}^{max}$) can be estimated as a given fraction $icfrac^{max}$ of the maximum total demand for the product ($\sum_{j \in I_c} dmd_{mjt}$) for each period $t \in T$ in the planning horizon. A random variable with uniform distribution which lies between 0 and $\sum_{j \in I_c} dmd_{mjt}$ is summed to the storage capacity estimation of each period. The minimum stock level of a product $m$ ($ic_{mi}^{min}$) is defined as a fraction $icfrac^{min}$ of the maximum inventory capacity.

The maximum technology capacities are estimated by the amount of final products that must be produced in the factories to meet the customer demands. Let us suppose the production is evenly distributed among the factories. In this case, the capacity of a technology $g$ to produce final products $m$ ($pc_g^{max}$) can be estimated as $\lceil \sum_{j \in I_c} \frac{dmd_{mjt}}{|I_f|} \rceil$. Since the production shall not necessarily be evenly distributed among the factories we multiply the resulting value by a random variable with uniform distribution between 2 and 4. The capacity of a technology $g$ is the maximum estimation over all final products in the planning horizon. The minimum use required by a technology $g$ ($pc_g^{min}$) is a fraction $lbpc^{min}$ of the maximum technology capacity. The minimum and maximum installation area of factories and warehouses are selected randomly from uniform distributions in pre-defined intervals.

Each supplier is allowed to provide $\frac{3}{|I_p|}$ of the total amount of raw material necessary to produce the maximum quantity of final products due to the technology capacities. The minimum order of a raw material $m \in M_{rm}$ from a supplier $i \in I_p$ ($sc_{mi}^{min}$) is drawn from a uniform distribution whose range

6

depends on scalars $lbsc$, $ubsc$ and $pc_g^{max}$. The distances between the entities have uniform distribution that falls in the interval $[lbdist, ubdist]$, where $lbdist$ and $ubdist$ are scalars such that $lbdist < ubdist$.

Table 3: Instance generator description: maximum and minimum capacities of entities and technologies.

| Parameter | Value | Variables |
|---|---|---|
| $ic_{mi}^{max}$ | $\max_{t \in T} \left[ icfrac^{max} \sum_{j \in I_c} dmd_{mjt} + Xic^{max} \sim U[0, \sum_{j \in I_c} dmd_{mjt}] \right]$ | $icfrac^{max} \in \mathbb{R}_{\geq 0}$ |
| $ic_{mi}^{min}$ | $icfrac^{min} ic_{mi}^{max}$ | $icfrac^{min} \in \mathbb{R}_{\geq 0}$ |
| $pc_g^{max}$ | $\max_{m \in M_{fp}, t \in T} \left[ (\sim U[2,4]) \left[ \lceil \sum_{j \in I_c} \frac{dmd_{mjt}}{|I_f|} \rceil \right] \right]$ | $icfrac^{max} \in \mathbb{R}_{\geq 0}$ |
| $pc_g^{min}$ | $lbpc^{min} pc_g^{max}$ | $lbpc^{min} \in \mathbb{R}_{\geq 0}$ |
| $ea_i^{max}$ | $\begin{cases} \sim U[lbeaf^{max}, ubeaf^{max}], & \text{if } i \in I_f \\ \sim U[lbeaw^{max}, ubeaw^{max}], & \text{if } i \in I_w \end{cases}$ | $lbeaf^{max}, ubeaf^{max},$ $lbeaw^{max}, ubeaw^{max} \in \mathbb{R}_{\geq 0}$ |
| $ea_i^{min}$ | $\begin{cases} \sim U[lbeaf^{min}, ubeaf^{min}], & \text{if } i \in I_f \\ \sim U[lbeaw^{min}, ubeaw^{min}], & \text{if } i \in I_w \end{cases}$ | $lbeaf^{min}, ubeaf^{min},$ $lbeaw^{min}, ubeaw^{min} \in \mathbb{R}_{\geq 0}$ |
| $sc_{mi}^{max}$ | $3 \lceil \frac{\sum_{(n,g) \in H_{prod}} pc_g^{max} BOM_{mng}^{prod}}{|I_p|} \rceil$ | - |
| $sc_{mi}^{min}$ | $\sim U[lbsc \, sc_{mi}^{max}, ubsc \, sc_{mi}^{max}]$ | $lbsc, ubsc \in \mathbb{R}_{\geq 0}$ |
| $dist_{ij}$ | $\sim U[lbdist, ubdist]$ | $lbdist, ubdist \in \mathbb{R}_{\geq 0}$ |

## 3.3 Parameters related to the objective functions.

In Table 4, one may observe the possible values of the parameters related to the costs associated with the production, operation and inventory of items. The cost of installing a technology $g$ ($tec_g$) is randomly defined with constant probability in an interval that depends on scalar multipliers $lbtech$ and $ubtech$ and on the maximum technology capacity, $pc_g^{max}$. The inventory price of a final product $m$ is defined by the weight of $m$ plus a given scalar $scfrac$. The cost to recover a product $m$ considers the bill ($BOM_{nm}^{rem}$), the price of the resulting final products ($psu_n$) and a random variable drawn from a uniform distribution ranging between two scalars $lbrpc$ and $ubrpc$. Since the cost to recover a product must be substantially smaller than its selling cost, we suggest $ubrpc$ to be 0.5 at most. The average fuel consumption is randomly selected between 14 and 18 liters per $100km$. This interval was chosen in accordance with the case studied by Mota et al. [1].

Table 4: Instance generator description: production, operation and inventory costs.

| Parameter | Value | Variables |
|---|---|---|
| $tec_g$ | $\sim U[lbtech \, pc_g^{max}, ubtech \, pc_g^{max}]$ | $lbtech, ubtech \in \mathbb{R}_{\geq 0}$ |
| $opc_g$ | $\sim U[lbopc, ubopc]$ | $lbopc, ubopc \in \mathbb{R}_{\geq 0}$ |
| $psu_m$ | $\sim U[lbpsu, ubpsu]$ | $lbpsu, ubpsu \in \mathbb{R}_{\geq 0}$ |
| $sc_m$ | $scfrac + pw_m$ | $scfrac \in \mathbb{R}_{\geq 0}$ |
| $rpc_m$ | $BOM_{nm}^{rem} psu_n Xrpc \sim U[lbrpc, ubrpc]$ | $lbrpc, ubrpc \in \mathbb{R}_{\geq 0}$ |
| $rmc_m$ | $\sim U[lbrmc, ubrmc]$ | $lbrmc, ubrmc \in \mathbb{R}_{\geq 0}$ |
| $avc_a$ | $\sim U[14, 18]$ | - |

Table 5 shows the parameters related to social and labor aspects of the supply chain. The construction cost of an entity $i \in I$, $sqmc_i$, is randomly drawn from a uniform distribution whose ranges depend on

7

the average labor cost of the area where the entity is located ($lc_i$) and the maximum installation area of the entity ($ea_i^{max}$). The lower bound for the distribution is half the value of $lc_i ea_i^{max}$ whereas the upper bound also depend on scalars $ubsqmc$ and $sqmcfac$.

In Table 5, $\mu^{gdp}$ is a regional factor, calculated as the inverse of the country Gross Domestic Product (GDP) per capita in euros where an entity $i \in I$ is located (GDPPC). Higher values of $\mu^{gdp}$ correspond to countries with lower GDP per capita. Parameter $pWork_i$ indicates the working conditions in entity $i \in I$. The higher the value of $pWork_i$, the better the working conditions. We set a lower bound of 0.3 for $pWork_i$ since the SSC should avoid considering activities in poor working conditions.

The average labor cost, the regional factor $\mu^{gdp}$ and the working conditions index strongly depend on the area where the entity is located. The GDPPC in which $\mu^{gdp}$ is based on is a strong indicator of the development of a country. In the case study used to define these parameters, more developed countries are often associated with higher average labor costs entities [1]. There are, however, some exceptions. Nonetheless, to avoid biased instances we do not base the parameters related to labor and construction nor the working condition index on the GDPPC.

Table 5: Instance generator description: social and labor aspects.

| Parameter | Value | Variables |
|---|---|---|
| $sqmc_i$ | $\sim U[(0.5 lc_i\ ea_i^{max}), (ubsqmc + sqmcfac\ lc_i\ ea_i^{max})]$ | $ubsqmc, sqmcfac \in \mathbb{R}_{\geq 0}$ |
| $lc_i$ | $\sim U[3.5, 30.4]$ | - |
| $\mu_i^{gdp}$ | $\sim U[0.355, 1.24]$ | - |
| $pWork_i$ | $\sim U[0.3, 1]$ | - |
| $wi_i$ | $\begin{cases} \sim [lbwf, ubwf], & \text{if } i \in I_f \\ \sim [lbww, ubww], & \text{if } i \in I_w \end{cases}$ | $lbwf, ubwf,$ $lbww, ubww \in \mathbb{R}_{\geq 0}$ |
| $wpsq_i$ | $\begin{cases} \sim [lbwpsqf, ubwpsqf], & \text{if } i \in I_f \\ \sim [lbwpsqw, ubwpsqw], & \text{if } i \in I_w \end{cases}$ | $lbwpsqf, ubwpsqf,$ $lbwpsqw, ubwpsqw \in \mathbb{R}_{\geq 0}$ |

The environmental impacts due to the installation of entities, production and remanufacturing of products, and transportation depend on the analysis of the life cycle assessment of the product. Table 6 presents guidelines for generating the remaining parameters. The lower and upper bound values of the uniform distribution intervals from which we base the values were proposed by Mota et al. [1] and are available in [3].

Table 6: Instance generator description: environmental impacts.

| Parameter | Value |
|---|---|
| $ei_c$ | $\sim U[lbenv_c^{ent}, upenv_c^{ent}]$ |
| $es_{ic}$ | $\sim U[lbenv_c^{sup}, upenv_c^{sup}]$ |
| $eg_{mgc}$ | $\sim U[lbenv_c^{prod}, upenv_c^{prod}]$ |
| $et_{ac}$ | $\sim U[lbenv_c^{transp}, upenv_c^{transp}]$ |


 
\end{document}